\documentclass[a4paper,10pt]{amsart}
\usepackage{amsmath,amsfonts,amssymb, amsthm, xypic}
\usepackage{epsfig,psfrag}
\usepackage[all]{xy}

\def\ind{{\lim\limits_{\to}}}

\def\ben{{\mathfrak{b}}}
\def\gen{{\mathfrak{g}}}
\def\hen{{\mathfrak{h}}}
\def\pen{{\mathfrak{p}}}

\def\nen{{\mathfrak{n}}}

\def\Oc{{\mathcal O}}

\def\det{{\operatorname{det}}}
\def\dim{\mathrm{dim}}
\def\Hom{\mathrm{Hom}}
\def\sign{\mathrm{sgn}}
\def\det{\mathrm{det}}

\def\End{\mathrm{End}}
\def\Sym{\mathrm{SYM}}
\def\Ind{\mathrm{Ind}}
\def\KK{\boldsymbol{\mathcal{K}}}
\def\Z{\mathbb{Z}}
\def\Q{\mathbb{Q}}

\def\H{\mathbf{H}}
\def\Cb{\mathbf{C}}
\def\Sb{\mathbf{S}}
\def\SCb{\mathbf{SC}}

\def\Vcb{\pmb{\mathcal V}}
\def\Vc{{\mathcal V}}
\def\B{\mathbb{A}}

\def\f{\mathbf{f}}
\def\Hin{\mathrm{Hilb}_n}
\def\Loc{{\mathrm{Loc}}}
\def\Hi{\mathrm{Hilb}}

\def\U{\boldsymbol{\mathcal{E}}}
\def\SH{\mathbf{S}\ddot{\mathbf{H}}}

\def\N{\mathbb{N}}

\def\a{\alpha}
\def\b{\beta}

\def\C{\mathbb{C}}

\def\Lb{{\boldsymbol{\Lambda}}}
\def\Kb{{\mathbf{K}}}
\def\Zb{{\mathbf{Z}}}
\def\x{\mathbf{x}}
\def\bt{\mathbf{u}}
\def\hb{\mathbf{h}}

\def\y{\mathbf{y}}

\def\z{\mathbf{z}}

\def\t{\tau}
\def\mod{{\text{mod}}}
\def\on{{\ \text{on}\ }}
\def\s{\sigma}
\def\bs{\bar{\sigma}}
\def\X{\mathcal{C}}
\def\hU{\widehat{\boldsymbol{\mathcal{E}}}}

\newtheorem{theo}{\bf{Theorem}}[section]
\newtheorem{lem}[theo]{Lemma}
\newtheorem{claim}[theo]{Claim}
\newtheorem{cor}[theo]{Corollary}
\newtheorem{prop}[theo]{Proposition}
\newtheorem{conj}[theo]{Conjecture}

\theoremstyle{remark}

\newtheorem{rem}[theo]{Remark}

\numberwithin{equation}{section}

\title[Elliptic Hall algebra and Hilbert schemes]{The elliptic Hall algebra and the K-theory of the Hilbert scheme
of $\mathbb{A}^2$}
\author{O. Schiffmann, E. Vasserot}

\begin{document}
\maketitle

\begin{abstract} In this paper we compute the convolution algebra in the equivariant $K$-theory of the Hilbert
scheme of $\mathbb{A}^2$.
We show that it is isomorphic to the elliptic Hall algebra, and hence to the spherical DAHA of $GL_\infty$.
We explain this coincidence via the geometric Langlands
correspondence for elliptic curves, by relating it also to the convolution algebra in the equivariant $K$-theory
of the commuting variety. We also obtain a few other related results ( action of the elliptic Hall algebra on the
$K$-theory of the moduli
space of framed torsion free sheaves over $\mathbb{P}^2$, virtual fundamental classes, shuffle algebras,...).

\end{abstract}

\tableofcontents

\section{Introduction and notation}

\vspace{.1in}

\paragraph{\textbf{0.1}}  Let $\text{Hilb}_n$ denote the Hilbert scheme of $n$ points in $\mathbb{C}^2$. The dimension of the homology groups $H^*(\text{Hilb}_n, \mathbb{Q})$, which were first determined by Ellingsr\"ud and Str\"omme \cite{ES}, are given by
$$\sum_{n \geqslant 0} \sum_{i=0}^{2n} \dim\;H^i(\text{Hilb}_n,\Q)p^iq^n=
\prod_{k \geqslant 1} \frac{1}{1-p^{2k}q^k}.$$
In a groundbreaking work \cite{Nak1}, see also \cite{Groj},
Nakajima obtained a more precise understanding of the space
$$\mathbf{V}=\bigoplus_n H^*(\text{Hilb}_n, \Q)$$
by geometrically realizing it as the Fock space representation of a Heisenberg algebra
$$\mathbb{H}=\C\langle p_{\pm 1}, p_{\pm 2}, \ldots \rangle / ([p_i,p_j]=i \delta_{i+j,0}).$$
More precisely, for $k \in \Z$ and $n, n+k \geqslant 0$ let
$Z_{n+k,n} \subset \text{Hilb}_{n+k} \times \text{Hilb}_n$ stand for
the nested Hilbert scheme and let $[Z_{n+k,n}] \in
H^*(\text{Hilb}_{n+k} \times \text{Hilb}_n)$ be its fundamental
class. The element $[Z_k]=\prod_{n} [Z_{n+k,n}]$ acts on
$\mathbf{V}$ by convolution. Nakajima's theorem may be stated as
follows.

\vspace{.1in}

\begin{theo}[Nakajima] The following hold
\begin{enumerate}
\item[(a)] the operators $\{[Z_k]\;;\; k \in \Z\}$ generate a Heisenberg algebra $\mathbb{H} \subset
End(\mathbf{V})$,
\item[(b)] as an $\mathbb{H}$-module, $\mathbf{V}$ is isomorphic to the Fock space representation, i.e.,
there exists an isomorphism 
$\mathbf{V} \simeq \C[p_1, p_2, \ldots]$
in which the action of the operators $[Z_k]$ is given by
\begin{equation}\label{E:Nak1}
1+\sum_{k \geqslant 1} [Z_k] z^k=
\exp \Big( -\sum_{n \geqslant 1} (-1)^n\frac{p_n}{n}z^n \Big),
\end{equation}
\begin{equation}\label{E:Nak2}
1+\sum_{k \geqslant 1} [Z_{-k}] z^k=
\exp\Big(-\sum_{n \geqslant 1}\frac{1}{n}\frac{\partial}{\partial p_n}z^n\Big).
\end{equation}
\end{enumerate}
\end{theo}

\vspace{.1in}

Note that this result extends to the Hilbert scheme of an arbitrary
smooth quasiprojective surface. One has to replace $\mathbb{H}$ by a
Heisenberg-Clifford algebra modeled over the homology lattice
$H_*(S,\mathbb{Z})$. See also \cite{QinWang1}, \cite{Lehn},
\cite{SL}, \cite{V} for other works in this direction.

In the well-documented analogy between the Hilbert-Chow resolution
$\text{Hilb}_n  \to S^n \C^2$ and the Springer resolution $T^*
\mathcal{B} \to \mathcal{N}$ of the nilpotent cone of a simple Lie
algebra $\mathfrak{g}$, Nakajima's construction corresponds to the
Springer representation of the Weyl group $W$ of $\mathfrak{g}$ in
the homology of the Springer fibers.

One aim of this paper is to generalize Nakajima's work for
equivariant $K$-theory. The torus $T = (\C^*)^2$ acts on $\C^2$ and
on $\text{Hilb}_n$ for any $n \geqslant 1$. Let $K^T(\text{Hilb}_n)$
denote the (algebraic) equivariant $K$-theory group and set
$$\mathbf{L}_R=\bigoplus_{n} K^T(\text{Hilb}_n),$$ a
$R=R_T=\C[q^{\pm 1}, t^{\pm 1}]$-module.
We are interested in a natural $K$-theoretic analog of the
convolution algebra considered by Nakajima. It is not a priori clear
how to define such a \textit{natural} convolution algebra because
the nested Hilbert schemes $Z_{n+k,n}$ are in general singular when
$k \not\in \{-1,0,1\}$ and the classes $\mathcal{O}_{Z_{n+k,n}}$ are
typically badly behaved.

We get around this difficulty by considering only the smooth nested
Hilbert schemes $Z_{n+1,n},$ $Z_{n,n}$,
$Z_{n,n+1}$ and their respective tautological bundles
$\boldsymbol{\tau}_{n+1,n}, \boldsymbol{\tau}_{n,n},
\boldsymbol{\tau}_{n,n+1}$. Namely, we let $\Kb=\C(q^{1/2},t^{1/2})$
and let $$\mathbf{H}_\Kb \subset \bigoplus_{k} \prod_n
K^T(\text{Hilb}_{n+k} \times \text{Hilb}_n) \otimes_{R} \Kb$$ be
the subalgebra generated by the elements
$$\mathbf{f}_{-1,l}=\prod_n \boldsymbol{\tau}^l_{n,n+1}, \qquad
\mathbf{f}_{1,l}=\prod_n \boldsymbol{\tau}^l_{n+1,n}, \qquad l \in
\Z$$ and the Adams operations
$$\mathbf{f}_{0,l}=\prod_n \Psi_l(\boldsymbol{\tau}_{n,n}), 
\qquad
\mathbf{f}_{0,-l}=\prod_n \Psi_{-l}(\boldsymbol{\tau}_{n,n}), 
\qquad
 \Psi_{-l}(\boldsymbol{\tau}_{n,n})=\Psi_l(\boldsymbol{\tau}_{n,n}^*),
\qquad l \geqslant 1.$$
Our first main result identifies $\mathbf{H}_\Kb$ and its action on
$\mathbf{L}_\Kb=\mathbf{L}_R
\otimes_{R} \Kb$ (see Theorem~\ref{T:Main}).

\vspace{.1in}

\begin{theo}\label{T:Intro1} The following hold
\begin{enumerate}
\item[(a)] there is an isomorphism $\Omega:\U_c \simeq \mathbf{H}_\Kb $,
where $\U_c$ is a certain one-dimensional central extension of the
spherical Double Affine Hecke Algebra $\mathbf{S}
\ddot{\mathbf{H}}_{\infty}$ of type $GL_\infty$,
\item[(b)] as an $\U_c$-module, $\mathbf{L}_\Kb$ is isomorphic to the standard representation
on the space of symmetric polynomials $\boldsymbol{\Lambda}_\Kb =\Kb[x_1, x_2, \ldots]^{\mathfrak{S}_{\infty}}$.
\end{enumerate}
\end{theo}

\vspace{.1in}

The algebras $\U_c$ and $\mathbf{S} \ddot{\mathbf{H}}_{\infty}$ were introduced in \cite{BS}, \cite{SV}.
The standard representation on  $\boldsymbol{\Lambda}_\Kb$, see Section~4.7,
involves operators of multiplication by symmetric polynomials, differential operators as well as Macdonald's
difference operators. It may be viewed, in a certain sense, as a  limit as $n \to \infty$ of the
polynomial representations of the spherical DAHAs of  type $GL_n$. The group $SL(2,\Z)$ acts
by automorphisms on $\mathbf{S} \ddot{\mathbf{H}}_{\infty}$  (but not on $\U_c$).
It would be interesting to find a geometric interpretation for this 
symmetry from the point of view of the Hilbert scheme.

The occurrence of Macdonald operators and polynomials comes as no
surprise here~: these are ubiquitous in Haiman's work on the Hilbert
scheme, see \cite{Haiman}. In fact the intertwining operator
$\mathbf{L}_\Kb \simeq \boldsymbol{\Lambda}_\Kb$ in
Theorem~\ref{T:Intro1} coincides with Haiman's identification of
$\mathbf{L}_\Kb$ with the space of symmetric polynomials.

While we were finishing this paper B. Feigin and A. Tsymbaliuk sent
us a copy of their preprint \cite{FT} where the equivariant
$K$-theory of the Hilbert schemes is studied by a different approach
based on shuffle algebras, see Section~9.

What about operators associated to higher rank nested Hilbert
schemes ? We do not know whether the classes
$[\mathcal{O}_{Z_k}]=\prod_n [\mathcal{O}_{Z_{n+k,n}}]$ belong to
$\mathbf{H}_\Kb$ or not for $k \not\in \{-1,0,1\}$. However, the
\textit{virtual classes} $\Lambda \boldsymbol{\mathcal{V}}_k=\prod_n
\Lambda \boldsymbol{\mathcal{V}}_{n+k,n}$ introduced by Carlsson and
Okounkov in \cite{CO} do indeed belong to $\mathbf{H}_\Kb$ and seem
to be much better behaved than the $\mathcal{O}_{Z_k}$. In
Theorem~\ref{T:Virt} and Corollary~\ref{C:Virt} we prove the
following.

\vspace{.1in}

\begin{theo}\label{T:Intro2} The following hold
\begin{enumerate}
\item[(a)] for any $k \in \Z$ the virtual class $\Lambda \boldsymbol{\mathcal{V}}_k$ belongs to
$\mathbf{H}_\Kb$,
\item[(b)] under the isomorphism $\mathbf{L}_\Kb \simeq \boldsymbol{\Lambda}_\Kb$ the
action of the operators $\Lambda \boldsymbol{\mathcal{V}}_k$ is given by
\begin{equation}
1+ \sum_{n \geqslant 1} \boldsymbol{\tau}^*_n \otimes \Lambda(\Vcb_n)z^n=
\exp\bigg( -\sum_{n \geqslant 1} (-1)^n \frac{1-t^nq^n}{1-q^n} p_n
\frac{z^n}{n}\bigg),
\end{equation}
\begin{equation}1+ \sum_{n \geqslant 1} \Lambda(qt\Vcb^*_{-n})z^n=
\exp\bigg( -\sum_{n \geqslant 1} \frac{1-t^nq^n}{1-t^n}
\frac{\partial}{\partial p_n} \frac{z^{-n}}{n}\bigg).
\end{equation}
\end{enumerate}
So the elements $\boldsymbol{\tau}^*_n \otimes \Lambda(\Vcb_n)$,
$\Lambda(qt\Vcb^*_{-n})$, $n \geqslant 0$, generate a Heisenberg
subalgebra of $\mathbf{H}_\Kb$.
\end{theo}

\vspace{.1in}

It would be very interesting to generalize the above
Theorem~\ref{T:Intro1} and Theorem~\ref{T:Intro2} to the case of the
Hilbert scheme of an arbitrary smooth quasiprojective surface $S$.
This would presumably involve a ``global'' version of
$\mathbf{S}\ddot{\mathbf{H}}_{\infty}$ living over the base $K(S)$.
The approach used in the present paper, based on localization to
$T$-fixed points, seems, unfortunately, restricted to toric
surfaces.

The Hilbert scheme $\text{Hilb}_n$ may be interpreted as the moduli
space of framed rank one torsion free-sheaves over
$\mathbb{P}^2(\C)$ of chern class $c_2=n$. For $r \geqslant 1$, let
$M_{r,n}$ be the moduli space of framed torsion free-sheaves over
$\mathbb{P}^2(\C)$ of chern class $c_2=n$ and rank $r$. One may
define ``nested schemes'' $M_{r,n+k,n}$ and tautological bundles
over them, and one may consider the convolution algebra
$\mathbf{H}^{(r)}_K$. The Theorem 1 has an analogue in this setting,
see Theorem~\ref{HR:Main}.

\vspace{.2in}

\paragraph{\textbf{0.2.}} 
Our original motivation for studying the convolution algebra $\mathbf{H}_K$
stems not so much from a desire to understand the geometry of
$\text{Hilb}_n$ as from its interpretation in the framework of the
geometric Langlands program for elliptic curves, which we now
describe. Recall that Beilinson and Drinfeld's geometric Langlands
program, for $GL_r$, predicts the existence, for any fixed smooth
projective curve $\X$ over $\C$, of an equivalence of derived
categories
\begin{equation}\label{E:geom}
Coh \big({\mathrm{Loc}_{r}\X}\big) \;\overset{D}{\simeq} \;
D\text{-}mod\big(\mathrm{Bun}_{r}\X\big).
\end{equation}
Here $\mathrm{Loc}_{r}\X$ is the moduli
stack parametrizing local systems over $\X$ of rank $r$ and
$\mathrm{Bun}_{r}\X$ is the moduli stack of vector bundles of rank
$r$ over $\X$. This correspondence should intertwine the natural
action of $Rep\;GL_r$ on $Coh ({\mathrm{Loc}_{r}\X})$ by tensor
product and the action of $Rep\;GL_r$ on
$D\text{-}mod(\mathrm{Bun}_{r}\X)$ by means of Hecke operators. The
image under this correspondence of the skyscraper sheaves at points
of $\mathrm{Loc}_{r}\X$ corresponding to irreducible local systems
has been determined in \cite{FGV}. The case which is relevant to us
is rather orthogonal. Namely, we are interested in coherent sheaves
supported in the formal neighborhood of the trivial local system
$\C^r_\X$ over $\X$.

Let us specialize from now on $\X$ to be an elliptic curve. Then
$$\mathrm{Loc}_{r}\X = \{\rho: \pi_1(\X)  \to GL_r\}/GL_r=\{(A,B)\in GL_r^2\;;\; AB=BA\}/GL_r$$
and the formal neighborhood of the trivial local system is the formal neighborhood of $(0,0)$ in
$$\widehat\Loc_r\X =\big\{(a,b) \in \mathfrak{gl}_r^2\;;\; [a,b]=0\big\}/GL_r=C_{\mathfrak{gl}_r}/GL_r$$
where $C_r \subset \mathfrak{gl}_r \times \mathfrak{gl}_r$ is the
commuting variety. Observe that this formal neighborhood is
independent of the choice of $\X$. Moreover, the torus $T =
(\C^*)^2$ naturally acts on $\widehat\Loc_r\X$ by 
$$(z_1, z_2)\cdot
(a,b)=(z_1a,z_2b).$$
The relation with the Hilbert schemes is given by the following simple observation. We have
$$Z_{n+r,n}=\big\{ (I,J)\;;\; I \subset J \subset
\C[x,y] \;\text{ideals},\; \text{codim}\; I=n+r, \,\text{codim}\;
J=n\big\}.$$ Hence if $(I,J) \in Z_{n+r,n}$ then $(x|_{J/I},
y|_{J/I})$ define a point in $C_r/GL_r$. One may try to use this map
$Z_{n+r,n} \to C_r/GL_r$ to lift classes from $\bar\Cb^r_R=K^{GL_r
\times T}(C_r)$ to $K^T(Z_{n+r,n})$. In this spirit we define a
convolution algebra structure on the direct sum
$$\bar\Cb_R=\bigoplus_{r \geqslant 0}\bar\Cb^r_R=
\bigoplus_{r \geqslant 0} K^T(\widehat\Loc_r\X)$$ as well as an
action of $\bar\Cb_\Kb=\bar\Cb_R \otimes_R \Kb$ on
$\mathbf{L}_\Kb$. Let $\H^>_\Kb\subset\H_\Kb$ be the subalgebra
generated by the classes $\f_{1,l}$ for $l \in \Z$, and let
$\overline\SCb_\Kb \subset \bar\Cb_\Kb$ be the subalgebra
generated by $\bar\Cb^1_R$.  The \emph{torsion free part} of
$\overline\SCb_\Kb$ is
$$\SCb_\Kb=\overline\SCb_\Kb/torsion.$$

\vspace{.1in}

\begin{theo}\label{T:Intro3} There is an algebra isomorphism
$\rho~: \SCb_\Kb \to \H^>_\Kb.$
It is compatible with the actions of $\overline\SCb_\Kb$ and 
$\H^>_\Kb$ on $\mathbf{L}_\Kb$.
\end{theo}

We conjecture that $\overline\SCb_K$ is actually torsion free (see Theorem~\ref{T:Isoktheory} and Conjecture~\ref{C:torsion}).
On the other hand, it is proved in \cite{SV} that the positive part
$\U^>\simeq \mathbf{S}\ddot{\H}^>_{\infty}$ is isomorphic to the
universal, spherical elliptic Hall algebra. In other words, for
any elliptic curve $\mathcal{C}_k$ defined over some finite field
$k=\mathbb{F}_l$ with Frobenius eigenvalues $\{\sigma,
\overline{\sigma}\}$ the algebra $\U^>$ specializes at $t=\sigma,
q=\overline{\sigma}$ to the spherical Hall algebra
$\mathbf{H}_{\mathcal{C}_k}^{sph}$ of $\mathcal{C}_k$, more
precisely, to its vector bundle part
$\mathbf{H}_{\mathcal{C}_k}^{sph,vec}$. By
definition $\mathbf{H}_{\X_k}^{sph,vec}$ is a convolution algebra of
functions on the moduli stacks $\bigsqcup_r \mathrm{Bun}_{r}\X_k$.
It is shown in \cite{Scano} that $\H^{sph,vec}_{\X_k}$ is the
Grothendieck group\footnote{For simplicity we ignore in this
introduction the (important) questions concerning completions of
Hall algebras, see \cite{BS}, \cite{Scano}.} of a certain category
$\mathcal{E}is_{\X_k}$ of semisimple perverse sheaves over
$\bigsqcup_r \mathrm{Bun}_{r}\X_k$, the so-called \textit{Eisenstein
sheaves} generated by the constant sheaves over the Picard stacks
$\mathrm{Pic}_{\X_k}^d$ for $d \in \Z$. There is a similar category
$\mathcal{E}is_{\X}$ for our complex elliptic curve $\X$. We expect
that $\mathcal{E}is_{\X}$ carries a natural $\Z^2$-grading and that
the associated graded Grothendieck group\footnote{again, we ignore
problems of completions.} $K_0^{gr}(\mathcal{E}is_{\X})$ is
isomorphic to $\U^>$ after a suitable extension of scalars.
Combining Theorem~1 and Theorem~3, we may draw the diagram
$$\xymatrix{ \bigoplus_r K_0(Coh^T(\widehat\Loc_r\X)) &  & \bigoplus_r K_0(Perv( \mathrm{Bun}_{r}\X))\\
\overline\SCb_\Kb \ar[r]^-{\rho}\ar[u] &
\mathrm{End}_\Kb\big(\mathbf{L}_\Kb\big) & \U^> =
K_0^{gr}(\mathcal{E}is_{\X}) \ar[l]_-{\Omega} \ar@<-2ex>[u]} $$ The
vertical arrows are embeddings. The last equality is conjectural.
The faithful actions of $\overline\SCb_\Kb$ and $\U^>$ on
$\mathbf{L}_\Kb$ identify these two algebras. They are both
naturally $\N$-graded by the rank $r$. For each $r$ the action of
$R_{GL_r \times T}$ on $\overline\SCb^r_\Kb$ coincides with the
action of $R_{GL_r \times T}$ on $\U^>[r]$ by means of Hecke
operators, see Proposition~\ref{P:Heckektheory}. We summarize this
in the following corollary, which is the second main result of this
paper.

\vspace{.1in}

\begin{cor} There is an algebra isomorphism
$\SCb_{\Kb}\simeq \U^>$ intertwining, on the graded
components of degree $r$, the natural actions of $R_{GL_r \times T}$
by tensor product and Hecke operators respectively.
\end{cor}

\vspace{.1in}

This corollary may be seen as a version, at the level of
Grothendieck groups, of a geometric Langlands correspondence
(\ref{E:geom}) for local systems living in the formal neighborhood
of the trivial local system (or, equivalently, for the category of
Eisenstein sheaves generated by constant sheaves over
$\mathrm{Pic}_{\X}$). We point, however, two notable differences
with (\ref{E:geom}). First, we have not just an equality of
Grothendieck groups $\SCb^r_{\Kb} \simeq \U^>[r]$ for
each $r$, but an equality of \textit{convolution algebras}. This
convolution product only exists because we are dealing with the
groups $GL_r$. Secondly, we have considered \textit{$T$-equivariant}
coherent sheaves over $\widehat\Loc_r\X$. There is no corresponding
$\Z^2$ grading on the categories $\mathcal{E}is_{\X}$ in our
picture. See \cite{SLectures2}, Lecture 5, for more in that
direction.

\vspace{.2in}

\paragraph{\textbf{0.3.}} The plan of the paper is as follows. Sections~1 and 2 serve as reminders on the elliptic Hall algebra $\U$
and the Cherednik algebra $\mathbf{S}\ddot\H_{\infty}$ and on the
Hilbert schemes respectively. The construction of the standard
representation of $\mathbf{S}\ddot\H^\geqslant_{\infty}$ and
$\U^\geqslant$ is given in Section~1.4. In Section~3 we consider the
convolution algebra $\H_\Kb$ in the equivariant $K$-theory of
Hilbert schemes and state our first main Theorem (see Theorem~1).
Section~4 is devoted to the proof of that result. In Section~5 we
explain how to relate the virtual classes
$\Lambda(\boldsymbol{\mathcal{V}}_n)$ with $\H_\Kb$. Section~6
contains some preparatory results pertaining to the Hecke action on
the Hall algebra $\U$. The convolution algebras $\bar\Cb_\Kb,$
$\overline\SCb_\Kb$ in the equivariant $K$-theory of the
commuting variety are defined and studied in Section~7, where
$\overline\SCb_\Kb$ is compared with the positive part $\H^>_\Kb$
(see the above Theorem~3 and Corollary). Section~8 deals with the
case of moduli spaces of framed torsion sheaves of rank $r > 1$ over
$\mathbb{P}^2(\C)$. Finally, Sections~9 and 10 contain
some remarks on some natural Heisenberg subalgebras of the convolution algebra
$\H^>_\Kb$ as well as a description of
$\H^>_\Kb$ as a shuffle algebra, and a comparison between the
present paper and the recent work \cite{FT}. In an attempt at making
this paper more readable, several technical points have been
postponed to the Appendix. The most important notations have been gathered
together in an index at the very end of the paper.

\vspace{.2in}

\paragraph{\textbf{0.4}}
We will use the standard notation for
partitions. If $\lambda=(\lambda_1\geqslant\lambda_2\dots)$ is a
partition then $\lambda'$ is its conjugate partition, $l(\lambda)$
is its length and $|\lambda|$ is the number of boxes in the Young
diagram associated with $\lambda$. We'll identify $\lambda$ and its
Young diagram, i.e., we write $\lambda=\{(i,j)\in\N^2; 1\leqslant
j\leqslant \lambda_i\}$. If $s \in \lambda$ is a box then we write
$s=(i,j)$. Let
$$a(s)=a_\lambda(s)=\lambda_i-j,\quad l(s)=l_\lambda(s)=\lambda'_j-i,$$
the arm length and the leg length of $s$. We also set
$$i(s)=i_\lambda(s)=i,\quad j(s)=j_\lambda(s)=j.$$
These of course do not depend on $\lambda$. Occasionally, we will write
$$x(s)=i(s)-1, \qquad  y(s)=j(s)-1.$$
So we have the following picture for $\lambda=(10,9^3,6,3^2)$

\vspace{.2in}

\centerline{
\begin{picture}(200,160)
\put(0,20){\line(0,1){70}}
\put(10,90){\line(0,-1){70}}
\put(0,90){\line(1,0){30}}
\put(20,90){\line(0,-1){70}}
\put(30,90){\line(0,-1){70}} \put(0,80){\line(1,0){30}}
\put(40,70){\line(0,-1){50}} \put(0,70){\line(1,0){60}}
\put(50,70){\line(0,-1){50}} \put(60,70){\line(0,-1){50}}
\put(0,60){\line(1,0){90}} \put(70,60){\line(0,-1){40}}
\put(80,60){\line(0,-1){40}} \put(90,60){\line(0,-1){40}}
\put(0,50){\line(1,0){90}} \put(0,40){\line(1,0){90}}
\put(0,30){\line(1,0){100}}
\put(0,20){\line(1,0){100}} \put(100,30){\line(0,-1){10}}
\put(33,43){{$s$}} \put(33,72){{$l$}} \put(34,8){{$j$}}
\put(-14,40){{$i$}} \put(94,40){{$a$}} \put(180,60){{$i(s)=3$}}
\put(180,47){{$j(s)=4$}} \put(180,34){{$a(s)=5$}}
\put(182,21){{$l(s)=2$}} \put(35,52){\line(0,1){16}}
\put(35,68){\line(1,-1){6}} \put(42,45){\line(1,0){46}}
\end{picture}}
\vspace{.05in}
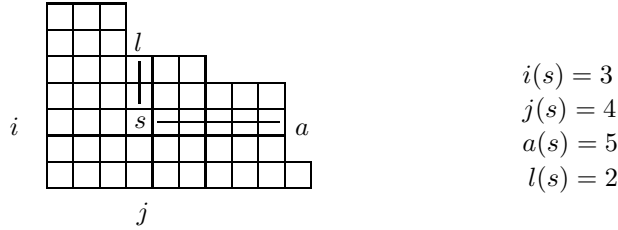
\centerline{\textbf{Figure 1.} Notations for partitions.}

\vspace{.2in}

Denote the set
of partitions by $\Pi$. For $s$ a box of a partition $\lambda$ we
let $C_s$ and $R_s$ be the column and row of $s$.

Given a field $\Kb$ let $\Lb_\Kb=\Kb[x_1, x_2, \ldots
]^{\mathfrak{S}_{\infty}}$ be the $\Kb$-vector space of symmetric
polynomials and $\Lb_\Kb^n=\Kb[x_1, x_2, \ldots, x_n
]^{\mathfrak{S}_n}$  be the space of symmetric polynomials in $n$
variables.

Let $E$ be a complex vector space of dimension $n$. Let $R_{GL(E)}$ denote
the complexified representation ring of $GL(E)$. There is a canonical algebra 
isomorphism $R_{GL(E)}\to\C[x_1^{\pm 1}, \ldots, x_n^{\pm 1}]^{\mathfrak{S}_n}$
 determined by
$$\Lambda^lE \mapsto e_l(x_1, \ldots, x_n)$$
for $1 \leqslant l \leqslant n$. Here $E$ is the tautological representation
of $GL(E)$. Under this isomorphism, the Adams operations
$\Psi_l(E)$ are mapped to the power sum functions
$p_l(x_1, \ldots, x_n)$ for any $l \in \Z^*$.

If $v$ is any variable, we put as usual
$$[n]_v=\frac{v^n-v^{-n}}{v-v^{-1}}, \qquad [n]_v!=[2]_v \cdots [n]_v.$$

\vspace{.2in}

\section{The elliptic Hall algebra, the spherical DAHA and the polynomial representation}

\vspace{.1in}

We first recall some result and definition from \cite{BS} and \cite{SV},
to which the reader is referred for details. Set $\Kb=\C(\s^{1/2}, \bs^{1/2})$,
where $\s,\bs$ are formal variables.

\vspace{.15in}

\paragraph{\textbf{1.1.}} We begin with some recollection on
the elliptic Hall algebra $\U^+$ and its Drinfeld double
$\widehat{\U}$. We set
$$\Zb=\Z^2,\quad
\Zb^*=\Zb \backslash \{(0,0)\},\quad
\Zb^+=\{ (i,j) \in \Zb; i > 0 \;\text{or}\; i=0, j>0\},\quad
\Zb^-=-\Zb^+.$$
For a future use we set also
$$\Zb^>=\{(i,j) \in \Zb; i > 0\},\quad
\Zb^k=\{(i,j) \in \Zb; i = k\},\quad \Zb^\geqslant=\{(i,j) \in \Zb;
i \geqslant 0\},$$ for each $k\in\Z$. For any $\x=(i,j) \in \Zb^*$
let $d(\x)$ denote the greatest common divisor of $i,j$.

We'll also use the following polynomials
$$\alpha_n=\alpha_n(\s,\bs)=(1-(\s\bs)^{-n})(1-\s^n)(1-\bs^n)/n.$$
Note that $\alpha_n=\alpha_{-n}$.

We set $\epsilon_{\x}=1$ if $\x \in \Zb^+$ and $\epsilon_{\x}=-1$ if $\x \in \Zb^-$.
For a pair of non-collinear elements
$\x,\y \in \Zb^*$ we set $\epsilon_{\x,\y}=\sign(\det(\x,\y))$,
an element in $\{\pm 1\}$.
Finally let $\boldsymbol{\Delta}_{\x,\y}$ be
the triangle in $\Zb$ with vertices $\{(0,0), \x, \x+\y\}$.

\vspace{.15in}

\noindent \textbf{Definition.} Let $\widehat{\U}$ be the
$\Kb$-algebra generated by elements $\bt_{\x}$,
$\boldsymbol{\kappa}_{\x}$, $\x \in \Zb^*,$
modulo the following set of relations
\begin{enumerate}
\item[(a)] $\kappa_{\x}$ is central for all $\x$, and $\kappa_{\x+\y}=\kappa_{\x}\kappa_{\y}$,
\item[(b)] if $\x,\y$ belong to the same line in $\Zb$ then
\begin{equation*}\label{E:0001}
[\bt_{\y},\bt_{\x}]=
\frac{\boldsymbol{\kappa}_{\x}-\boldsymbol{\kappa}^{-1}_{\x}}{\alpha_{d(\x)}
}\ \text{if}\ \x=-\y,\quad [\bt_{\y},\bt_{\x}]=0\ \text{else},
\end{equation*}
\item[(c)] if $\x,\y \in \Zb^*$ are such that $d(\x)=1$ and that
$\boldsymbol{\Delta}_{\x,\y}$ has no interior lattice point then
\begin{equation*}\label{E:0002}
[\bt_\y,\bt_{\x}]=\epsilon_{\x,\y}\boldsymbol{\kappa}_{\a(\x,\y)}
\frac{\theta_{\x+\y}}{\alpha_{1}},
\end{equation*}
where
\begin{equation*}
\a(\x,\y)=\begin{cases}
\epsilon_\x(\epsilon_{\x}\x+\epsilon_{\y}\y-\epsilon_{\x+\y}(\x+\y))/2 & \text{if}\; \epsilon_{\x,\y}=1,\\
\epsilon_\y(\epsilon_{\x}\x+\epsilon_{\y}\y-\epsilon_{\x+\y}(\x+\y))/2
& \text{if}\; \epsilon_{\x,\y}=-1,
\end{cases}
\end{equation*}
and where the elements $\theta_{\z}$, $\z \in \Zb$, are given
by
\begin{equation*}\label{E:0003}
\sum_{l\geqslant 0} \theta_{l\x}s^l=\exp\big(\sum_{r \geqslant
1}\alpha_r\bt_{r\x}s^r\big),
\end{equation*}
for any $\x \in \Zb^*$ such that $d(\x)=1$.
\end{enumerate}

\vspace{.15in}

For a future use, note that $\theta_{\z}=\alpha_1\bt_{\z}$ if
$d(\z)=1$. Observe also that the $\Kb$-algebra $\widehat{\U}$
is $\Zb$-graded and has a natural $SL(2,\Z)$-symmetry. We write
$\boldsymbol{\mathcal{K}}$ for the central subalgebra generated by
$\{\boldsymbol{\kappa}_\x; \x \in \Zb\}$. We'll abbreviate
$\mathbf{c}_1=\boldsymbol{\kappa}_{0,1}$ and
$\mathbf{c}_2=\boldsymbol{\kappa}_{1,0}$. We have
$\boldsymbol{\mathcal{K}}=\Kb[\mathbf{c}_1^{\pm 1},
\mathbf{c}_2^{\pm 1}]$. We'll view $\widehat{\U}$ as a
$\KK$-algebra. 

Let $\hU^{\pm}$ be the $\KK$-subalgebra of $\widehat{\U}$
generated by $\{\bt_{\x};\x \in \Zb^{\pm}\}$. It is shown in
\cite{BS}, Section~5, that relations (b), (c) restricted to
$\Zb^{\pm}$ give a presentation of $\U^{\pm}$, and that there is
a biangular decomposition
\begin{equation}\label{E:D1}
\widehat{\U} \simeq \hU^+  \otimes \hU^-.
\end{equation}
For our convenience we put $\bt_{0,0}=1$. It is also proved in
loc.~cit.~  that the $\boldsymbol{\mathcal{K}}$-algebra
$\widehat{\U}$ is generated by the elements $\bt_\x$ with
$\x\in\Zb^1\cup\Zb^0\cup\Zb^{-1}$. It will be helpful to consider a
slight refinement of the isomorphism (\ref{E:D1}). Let ${\hU}^{>}$,
${\hU}^{<}$, ${\hU}^{0}$ be the subalgebras of $\widehat{\U}$
generated by the elements ${\bt}_\x$ with $\x\in\Zb^1$, $\Zb^{-1}$,
$\Zb^0$ respectively. From (\ref{E:D1}) and the defining relations
of $\widehat{\U}$ one deduces the triangular decomposition
\begin{equation}\label{E:D2}
\widehat{\U} \simeq \hU^{>}  \otimes \hU^{0} \otimes \hU^{<}.
\end{equation}
The algebras $\hU^{>}$, $\hU^{<}$ are isomorphic. By
\cite{BS}, Corollary 5.2, the
algebras $\hU^>$, $\hU^<$ contains all the elements $\bt_\x$
for $\x\in\Zb^>$,  $\Zb^<$ respectively.
Let $\hU^{\geqslant}$ be the subalgebra generated
by $\hU^>$, $\hU^0$.
We define $\hU^{\leqslant}$ in the same way. We have
$$\hU^>\subset\hU^+\subset\hU^\geqslant.$$

For any $\omega=(\omega_1,\omega_2) \in (\C^*)^2$, write $\U_\omega$
for the specialization of $\hU$ at $\mathbf{c}_1=\omega_1,
\mathbf{c}_2=\omega_2$. The notations $\U^{\pm}_{\omega}$,
$\U^{>}_{\omega}$, $\U^<_{\omega}$, etc, are clear. We will simply
write $\U$ for $\U_{(1,1)}$. We'll be mostly interested in the
algebra $\U_c$ where
$$c=(1,q^{1/2}t^{1/2}).$$

When $\s$, $\bs$ are specialized to the Frobenius eigenvalues of an elliptic
curve ${X}$ defined over a finite field, $\U^+$ is equal to the
spherical Hall algebra of ${X}$. It is endowed with a natural
coproduct and $\widehat{\U}$ is isomorphic to its Drinfeld double.

\vspace{.2in}

\paragraph{\textbf{1.2.}} The following result will be useful.

\vspace{.1in}

\begin{prop}\label{P:charac} The algebra $\widehat{\U}$ is isomorphic to
the algebra generated by $\hU^<, \hU^0, \hU^>$ modulo the following relations
\begin{equation}\label{E:2}
\begin{split}
[\bt_{0, k},\bt_{1,l}]&=  \,\bt_{1,l+k}, \qquad  k >0,\\
[\bt_{0, k},\bt_{1,l}]&=  -\mathbf{c}_1^k \, \bt_{1,l+k}, \qquad k <0,
\end{split}
\end{equation}
and
\begin{equation}\label{E:205}
\begin{split}
[\bt_{-1, l},\bt_{0,k}]&=  \mathbf{c}_1^{-k}\,\bt_{-1,k+l}, \qquad k >0,\\
[\bt_{-1, l},\bt_{0,k}]&=  - \, \bt_{-1,k+l}, \qquad k <0,
\end{split}
\end{equation}
and finally
\begin{equation}\label{E:5}
\begin{split}
[\bt_{-1,k},\bt_{1,l}]&= \mathbf{c}_2\mathbf{c}_1^{-k} \frac{\theta_{0,k+l}}{\alpha_1},\qquad k+l>0,\\
[\bt_{-1,k},\bt_{1,l}]&=\frac{\mathbf{c_1}^{-k}\mathbf{c}_2-\mathbf{c}_1^k\mathbf{c}_2^{-1}}{\alpha_1},
\qquad k=-l,\\
[\bt_{-1,k},\bt_{1,l}]&=- \mathbf{c}^{-1}_2\mathbf{c}_1^{-l} \frac{\theta_{0,k+l}}{\alpha_1},\qquad k+l<0.
\end{split}
\end{equation}
\end{prop}
\begin{proof} Let $\widetilde{\U}$ momentarily denote the algebra
generated by $\hU^>,\hU^0,\hU^<$ modulo the relations
(\ref{E:2}--\ref{E:5}). We have a natural surjective morphism $\phi:
\widetilde{\U} \to \widehat{\U} $ which is the identity on the
subspace $\hU^> \hU^0 \hU^<$ of $\widetilde\U$. Hence it suffices to
prove that $\hU^> \hU^0 \hU^<=\widetilde{\U}$. The algebras $\hU^>$,
$\hU^<$ are generated by the elements $\bt_\x$ with $\x\in\Zb^1$,
$\Zb^{-1}$ respectively. Hence we have to prove that any monomial of
$\widetilde\U$ in the elements $\bt_{1,l}, \bt_{0,n}, \bt_{-1,k}$
may be straightened to similar monomials in which the $\bt_{1,l},
\bt_{0,n}$ and $\bt_{-1,k}$ appear in that order. It is easy to see
that relations (\ref{E:2})-(\ref{E:5}) precisely enable one to
perform such straightening.\end{proof}

\vspace{.2in}

\paragraph{\textbf{1.3.}} We briefly recall here the relation between $\U$
and the spherical double affine Hecke algebras. See
\cite{SV}, Section~2, for details
\footnote{The conventions used in \cite{SV} differ slightly from the
more standard ones of \cite{C}.}.
Let $\ddot{\H}_n$ be the double affine Hecke algebra (DAHA for
short) of type $GL_n$.
It is a $\C(q^{1/2},t^{1/2})$-algebra generated by elements
$X^{\pm 1}_i, Y^{\pm 1}_i$ for $i=1, \ldots, n$ and $T_j$ for $j=1,
\ldots, n-1$ subject to some relation which we won't write here. Let
${S}$ be the complete idempotent in the finite Hecke algebra $\H_n
\subset \ddot{\H}_n$ generated by the $T_1,T_2,\dots T_{n-1}$. It is
characterized by the property that $T_j S=S T_j=t^{-1/2}S$ for all
$j$. The spherical DAHA is the subalgebra
$\mathbf{S}\ddot{\H}_n={S}\,
 \ddot{\H}_n  {S}$ of $\ddot{\H}_n$. For $l
>0$ we set
$$P^n_{0,l}=S  \sum_i Y_i^l S, \quad P^n_{0,-l}=q^l S \sum_i Y_i^{-l}
S,\quad P^n_{l,0}=q^lS  \sum_i X_i^l S, \quad P^n_{-l,0}= S
\sum_i X_i^{-l} S.$$ These elements generate
$\mathbf{S}\ddot{\H}_n$. More general elements $P^n_{\x}$ for $\x
\in \Zb^*$ are defined in \cite{SV} via the $SL(2,\Z)$ action on
$\mathbf{S}\ddot{\H}_n$. We put also $P^n_{0,0}=1$.

Let $\mathbf{S}\ddot{\mathbf{H}}^{+}_n$ be the subalgebra of
$\mathbf{S}\ddot{\mathbf{H}}_n$ generated by the elements $P_{\x}^n$
for $\x \in \Zb^+$. It is shown in \cite{SV}, Proposition~4.1, that
the assignment $P^n_{\x} \mapsto P^m_{\x}$ extends to a surjective
algebra homomorphism $\mathbf{S}\ddot{\mathbf{H}}^+_n \to
\mathbf{S}\ddot{\mathbf{H}}^+_m$ for $n >m$. This allows one to
consider the limit $ \mathbf{S}\ddot{\H}^+_\infty \subset
\prod_n \mathbf{S}\ddot{\mathbf{H}}^+_n$ of the algebras
$\mathbf{S}\ddot{\mathbf{H}}^+_n$, along with its set of generators
$\{P^{\infty}_{\x};\x \in \Zb^+\}$.

We identify $\Kb$ and $\C(q^{1/2},t^{1/2})$ by means of
\begin{equation}\label{E:qHtH}
\s \mapsto q^{-1}, \qquad \bs \mapsto t^{-1}.
\end{equation}

\vspace{.1in}

\begin{theo}[\cite{SV}]\label{T:SV} For any $n$ the assignment
$$\bt_{\x} \mapsto \frac{1}{q^{d(\x)}-1}P_\x^n,\quad\x \in \Zb^* $$
extends to a surjective algebra homomorphism
$$\Phi_n: \U \to \mathbf{S}\ddot{\H}_n.$$
The restrictions $\Phi^+_n$ of $\Phi_n$ to $\U^+$ form a projective
system of maps and give rise to an algebra isomorphism
\begin{equation*}
\Phi^+_\infty: \U^+ \to
\mathbf{S}\ddot{\H}^+_\infty,\quad \bt_{\x} \mapsto
\frac{1}{q^{d(\x)}-1} P^{\infty}_\x.
\end{equation*}
\end{theo}

Note that it is not clear a priori how to make sense of a stable limit $\mathbf{S}\ddot{\H}_{\infty}$. For this reason, we may take $\U$ as a \textit{definition} of $\mathbf{S}\ddot{\H}_{\infty}$.

\vspace{.2in}

\paragraph{\textbf{1.4.}}
We now define a faithful representation of $\U^+$. First, recall
that there is a faithful representation $\varphi_n$ of $\ddot{\H}_n$
on the space $\Kb[x_1^{\pm 1}, \ldots, x_{n}^{\pm 1}]$ defined by
\begin{align*}
\varphi_n(X_i)&=x_i,\\
\varphi_n(T_i)&=t^{1/2}s_i + \frac{t^{-1/2}-t^{1/2}}{x_i/x_{i+1}-1}(s_i-1),\\
\varphi_n(Y_i)&=\varphi_n(T_i) \cdots \varphi_n(T_{n-1})\omega
\varphi_n(T_1^{-1}) \cdots \varphi_n(T_{i-1}^{-1}).
\end{align*}
Here $s_i$ stands for the transposition $x_i \leftrightarrow
x_{i+1}$, $\partial_i$ is the linear operator defined by $\partial_i
\cdot p(x_j)=p(q^{\delta_{i,j}} x_j)$ and $\omega=s_{n-1} \cdots s_1
\partial_1$. See \cite{C} or \cite{SV}, Section~4.

Let  $\SH^{>}_n, \SH^0_n,\SH^{\geqslant}_n\subset\SH_n$ be the
subalgebras generated by the elements $P^n_\x$ with $\x \in \Zb^>$,
$\Zb^0$, $\Zb^{\geqslant}$ respectively. Note that
$$\SH^>_n\subset\SH^+_n\subset\SH^\geqslant_n.$$
The representation $\varphi_n$ restricts to a representation of
$\SH^{\geqslant}_n$ on the space of symmetric polynomials $\Lb^n_\Kb$. In this representation the
element $P_{0,1}^n$ acts as Macdonald's difference operator
$$\Delta_1^n=\sum_{i=1}^n \left(\prod_{j \neq i} \frac{t^{-1/2}x_i-t^{1/2}x_j}{x_i-x_j} \right)\partial_i.$$
It has as eigenvectors the Macdonald polynomials
$P^{n}_{\lambda}(q,t^{-1})$ where $\lambda$ runs among all partitions
with at most $n$ parts. The $\Delta_1^n$-eigenvalue of
$P^{n}_{\lambda}(q,t^{-1})$ is given by the formula
\begin{equation*}\label{E:eig1}
\beta^1_{\lambda,n}=t^{-\frac{n-1}{2}}\sum_{i=1}^n
q^{\lambda_i}t^{(i-1)},
\end{equation*}
see \cite{Mac}, Section~VI.3. For $l>1$ the element $P^n_{0,l}$ acts
as the linear operator $\Delta^n_l$ whose eigenvalues on the
Macdonald polynomials read
\begin{equation}\label{E:eig2}
\beta^{l}_{\lambda,n}=t^{-l\frac{n-1}{2}}\sum_{i=1}^n
q^{l\lambda_i}t^{l(i-1)},
\end{equation}
see \cite{Mac}, Section~VI.4, (4.15). Similar results also hold for
the elements $P^n_{0,-l}$ with $l \geqslant 1$. Namely, the operator
$\Delta^n_{-l}=\varphi_n(P^n_{0,-l})$ is diagonalizable in the basis
of Macdonald polynomials with eigenvalues given by
\begin{equation}\label{E:eig3}
\beta^{-l}_{\lambda,n}=q^lt^{l\frac{n-1}{2}}\sum_{i=1}^n q^{-l\lambda_i}t^{-l(i-1)}.
\end{equation}
There is a unique algebra homomorphism $\rho_n: \Lb_\Kb^n \to
\Lb_\Kb^{n-1}$ given by $x_i \mapsto x_i$ for $i \leqslant n-1$ and
$x_n \mapsto 0$. It takes $P^n_{\lambda}(q,t^{-1})$ to
$P^{n-1}_{\lambda}(q,t^{-1})$ if $l(\lambda) \leqslant n-1$ and to zero
if $l(\lambda)=n$. It is easy to see from the formulas
(\ref{E:eig2}), (\ref{E:eig3}) that the map $\rho_n$ is not
compatible with the representations $\varphi_n$, i.e., we have
$$\varphi_{n-1} \circ \rho_n \neq \rho_{n} \circ \varphi_n.$$
To remedy this, we need to renormalize the action $\varphi_n$
slightly. Put
$$\widetilde{\Delta}^n_l=t^{l\frac{n-1}{2}}(\Delta_l^n-[n]_{t^{l/2}}).$$
One checks using (\ref{E:eig2}), (\ref{E:eig3}) that we have $\rho_n
\circ \widetilde{\Delta}^n_l = \widetilde{\Delta}_l^{n-1} \circ
\rho_n$ and that for each $l\geqslant 1$ the eigenvalue of
$\widetilde{\Delta}^n_l$ on $P^n_{\lambda}(q,t^{-1})$ is equal to
\begin{equation}\label{E:eig4}
\widetilde{\beta}^{l}_{\lambda,n}=\sum_{i=1}^n
(q^{l\lambda_i}-1)t^{l(i-1)},\quad
\widetilde{\beta}^{-l}_{\lambda,n}=q^l\sum_{i=1}^n
(q^{-l\lambda_i}-1)t^{-l(i-1)}.
\end{equation}
We denote by $\widetilde{\Delta}^\infty_l$ the stable limit of the
operators $\widetilde{\Delta}^n_l$.

\vspace{.1in}

By \cite{BS}, Corollary 5.2, and Theorem~\ref{T:SV} the algebra $\SH^>_n$
is generated by the elements $P^n_\x$ with $\x\in\Zb^1$ and we have
the following relation in $\SH_n^\geqslant$
\begin{equation}\label{E:semidir}
[P^n_{0,l}, P^n_{1,j}]=\epsilon_l(q^l-1)P^n_{1,j+l}.
\end{equation}
Here $\epsilon_l=1$ if $l \geqslant 1$ and $\epsilon_l=-1$ if $l \leqslant
-1$. Let  $\SH_n^> \rtimes \SH_n^0$ be the algebra
generated by $\SH^>_n$ and $\SH^0_n$ modulo the relations
(\ref{E:semidir}).

\begin{lem}\label{L:troy}
The $\Kb$-algebra $\SH^{\geqslant}_n$ is isomorphic to $\SH_n^> \rtimes
\SH_n^0$.
\end{lem}
\begin{proof}
The multiplications in $\SH_n^> \rtimes \SH_n^0$ and
$\SH^{\geqslant}_n$ yield right $\SH_n^0$-module homomorphisms
$$\SH_n^> \otimes_\Kb \SH_n^0\to\SH_n^> \rtimes \SH_n^0 \to \SH^{\geqslant}_n.$$
The first map is surjective by (\ref{E:semidir}). The second one is
also surjective by Theorem~\ref{T:SV} and (\ref{E:D2}). We must
prove that it is injective. It is enough to prove that the
surjective map
$$m:\SH_n^> \otimes_\Kb \SH_n^0\to\SH^{\geqslant}_n$$
given by the multiplication in $\SH^{\geqslant}_n$ is injective.

There is a unique $\Zb$-grading on the algebra $\SH_n^\geqslant$
such that the element $P_\x^n$ has the degree $\x$. For each integer
$i\geqslant 0$ let ${}^i\SH_n^{\geqslant}\subset\SH_n^{\geqslant}$
be subspace spanned by the homogeneous elements whose degree belongs
to $\mathbf Z^i$. It is a right $\SH_n^0$-submodule of
$\SH^{\geqslant}_n$ such that
$\SH^{\geqslant}_n=\bigoplus_{i\geqslant 0}{}^i\SH^\geqslant_n$. The
map $m$ is a surjective right graded $\SH^0_n$-module homomorphism.
Thus it restricts to a surjective right $\SH^0_n$-module
homomorphism
$$m^i:{}^i\SH_n^>\otimes_\Kb \SH_n^0 \to {}^i\SH^\geqslant_n,\quad
{}^i\SH_n^>=\SH_n^>\cap {}^i\SH^\geqslant_n.$$

Let $A$ be the localization of the ring $\C[q^{\pm 1/2},t^{\pm
1/2}]$ with respect to the multiplicative set generated by
$[n]_{t^{1/2}}$. Let $\ddot\H_{n,A}^\geqslant\subset\ddot\H_n$ be
the $A$-subalgebra generated by the elements $X_i, Y_i^{\pm 1},
T_j$, and
let  $\SH^{>}_{n,A}, \SH^0_{n,A}$ be the
$A$-subalgebras generated by the elements $P^n_\x$ with $\x \in
\Zb^>$, $\Zb^0$ respectively.
We put $\SH^{\geqslant}_{n,A}=S\ddot\H_{n,A}^\geqslant S$.
We define also
${}^i\SH^{>}_{n,A}, {}^i\SH^{\geqslant}_{n,A}$ in the obvious way.
The map $m^i$ restricts to an $A$-linear map
$$m_A^i:{}^i\SH_{n,A}^> \otimes_A \SH_{n,A}^0\to{}^i\SH^{\geqslant}_{n,A}.$$
We'll abbreviate $\C=A/(q-1,t-1)$. We have the following
\begin{itemize}
\item ${}^i\SH_{n,A}^>$ is a finitely generated $A$-module such that
${}^i\SH_{n,A}^>\otimes_A\Kb={}^i\SH_n^>$,
\item
$\SH^{\geqslant}_{n,A}$ is a free
$\SH^0_{n,A}$-module, ${}^i\SH^{\geqslant}_{n,A}$ is a direct
summand of finite rank,
and ${}^i\SH^\geqslant_{n,A}\otimes_A\Kb={}^i\SH^{\geqslant}_n$,
\item
$m^i_A\otimes$ id$_\C$ is an isomorphism.
\end{itemize}
The first claim is obvious. It implies that
$$\dim_\Kb({}^i\SH_{n}^>)\leqslant\dim_\C({}^i\SH_{n,A}^>\otimes_A\C).$$
The two other yield the reverse inequality.
Indeed, the third claim implies that
${}^i\SH_{n,A}^\geqslant\otimes_A\C$ is a free
${}^i\SH_{n,A}^0\otimes_A\C$-module of rank
$\dim_\C({}^i\SH_{n,A}^>\otimes_A\C)$.
Thus the second claim implies that
${}^i\SH_{n}^\geqslant$ is a free
${}^i\SH_{n}^0$-module of rank
$\dim_\C({}^i\SH_{n,A}^>\otimes_A\C)$.
Hence the surjectivity of $m^i$ implies that
$$\dim_\Kb({}^i\SH_{n}^>)\geqslant\dim_\C({}^i\SH_{n,A}^>\otimes_A\C).$$
Therefore $m^i$ is a surjective homomorphism of projective
${}^i\SH_{n}^0$-modules of the same (finite) rank.
Hence it is invertible.

The proof of the third claim is immediate.
It is also easy to check that
$\SH_{n,A}^\geqslant$ is a free
$\SH^0_{n,A}$-module, because
it is a direct summand of
$\ddot\H_{n,A}^\geqslant$, the later is free over
$\SH^0_{n,A}$ by the PBW-theorem and the Steinberg-Pittie theorem,
and $\SH^0_{n,A}$ is a polynomial ring.
Therefore we are reduced to check that
$\SH^{\geqslant}_{n}=S\ddot\H_{n}^\geqslant S$.
This is proved as in \cite{SV}, Proposition 2.6.
\end{proof}

\vspace{.1in}

By Lemma~\ref{L:troy} above we may construct a family of
automorphisms of $\SH^{\geqslant}_n$ by the assignment
\begin{equation}\label{E:autom1}
\begin{split}
P^n_{0,l}\mapsto P^n_{0,l} + x_l\ \text{if}\ l \in \Z^*,\quad
P^n_{\x}\mapsto P^n_{\x}\ \text{if}\  \x\in \Zb^>.
\end{split}
\end{equation}
Here $x_l\in\Kb$, $l \in \Z^*$, are arbitrary. The algebra
$\SH^{\geqslant}_n$ is $\Zb$-graded. Thus we can also construct
automorphisms by the assignment
\begin{equation}\label{E:autom2}
P_{i,j}^n \mapsto y_1^iy_2^j P_{i,j}^n
\end{equation}
for all $i,j$. Here $y_1, y_2\in\Kb$ are arbitrary.
Thus there is an automorphism $\theta$ of
$\SH^{\geqslant}_n$ such that
\begin{alignat*}{2}
\theta(P^n_{0,l})=t^{l\frac{n-1}{2}}(P_{0,l}^n-[n]_{t^{l/2}})\
\text{if}\ l \in \Z^*,\quad
\theta(P^n_{i,j})=t^{j\frac{n-1}{2}}P_{i,j}^n\ \text{if}\ (i,j)\in
\Zb^>.
\end{alignat*}
Set $\widetilde{\varphi}_n=\varphi_n \circ \theta$, a representation
of ${\mathbf{S\ddot H}}_n^{\geqslant}$. The discussion above implies
that
$$\widetilde{\varphi}_{n-1}(P_{0,l}^{n-1}) \circ \rho_n = \rho_{n} \circ
\widetilde{\varphi}_n(P^n_{0,l}),
\quad\widetilde{\varphi}_{n-1}(P^{n-1}_{1,0})\circ
\rho_n=p_1=\rho_{n} \circ \widetilde{\varphi}_n(P^n_{1,0}).$$ Since
the algebra $\SH^{\geqslant}_n$ is generated by the elements
$P^n_{0,l}$, $P^n_{1,0}$ we have $\widetilde{\varphi}_{n-1} \circ
\rho_n = \rho_{n} \circ \widetilde{\varphi}_n$. Thus the
representations $\widetilde{\varphi}_n$ yield a representation
$\widetilde{\varphi}_{\infty}$ of $\SH^{\geqslant}_\infty$ on the
$\Kb$-vector space $\Lb_\Kb.$ This representation is faithful
because each of the $\widetilde{\varphi}_n$ is faithful. Further, it
may be characterized as in the following proposition. Let
$P_{\lambda}(q,t^{-1})\in\Lb_\Kb$ be the Macdonald polynomial associated
with the partition $\lambda$.

\vspace{.1in}

\begin{prop} There is a unique representation
$\widetilde{\varphi}_{\infty}$ of $\SH^{\geqslant}_\infty$ on
$\Lb_\Kb$ such that
$\widetilde{\varphi}_{\infty}(P^{\infty}_{1,0})=p_1$ and such that
for any partition $\lambda$ and any integer $l\geqslant 1$ we have
$$\widetilde{\varphi}_{\infty}(P^{\infty}_{0,l}) \cdot P_{\lambda}(q,t^{-1})=
\left(\sum_i (q^{l\lambda_i}-1)t^{l(i-1)} \right) P_{\lambda}(q,t^{-1})$$
and
$$\widetilde{\varphi}_{\infty}(P^{\infty}_{0,-l}) \cdot P_{\lambda}(q,t^{-1})=
q^l\left(\sum_i (q^{-l\lambda_i}-1)t^{-l(i-1)} \right) P_{\lambda}(q,t^{-1}).$$
This representation is faithful.
\end{prop}

\vspace{.1in}

\begin{cor}\label{C:pol} There is a unique representation
$\widetilde{\varphi}$ of $\U^{\geqslant}$ on $\Lb_\Kb$ such that
$\widetilde{\varphi}(\bt_{1,0})=p_1/(q-1)$ and such that for any partition
$\lambda$ and any integer $l \geqslant 1$ we have
$$\widetilde{\varphi}(\bt_{0,l}) \cdot P_{\lambda}(q,t^{-1})=
\left(\sum_i \frac{q^{l\lambda_i}-1}{q^l-1}t^{l(i-1)} \right)
P_{\lambda}(q,t^{-1})$$
and
$$\widetilde{\varphi}(\bt_{0,-l}) \cdot P_{\lambda}(q,t^{-1})=
-\left(\sum_i \frac{q^{-l\lambda_i}-1}{q^{-l}-1}t^{-l(i-1)} \right)
P_{\lambda}(q,t^{-1}).$$
This representation is faithful.
\end{cor}

We call $\widetilde{\varphi}$ the \textit{polynomial}, or \textit{standard}
representation of $\U^{\geqslant}$. In Section~4.7 we will extend this
representation to the whole algebra $\widehat{\U}$.

\vspace{.2in}

\section{Hilbert schemes of $\mathbb{A}^2$}

\vspace{.1in}

This section contains some standard fact on Hilbert schemes.
Most of this may be found in \cite{ES}, \cite{Go}.

\vspace{.2in}

\paragraph{\textbf{2.1.}} Throughout the paper, our ground field will be $\C$.
Let $\Hin$ denote the Hilbert scheme parametrizing length $n$ subschemes of 
$\B^2$.
By Fogarty's theorem it is a smooth irreducible variety of dimension $2n$.
Associating to a closed point of $\Hin$ its ideal sheaf yields a bijection
$$\Hin(\C) =\{I \subset \C[x,y]; I\;\text{is\;an\;ideal\;of\;codimension\;} n\}.$$
We will denote by $S=\C[x,y]$ the ring of regular functions on
$\B^2$. The tangent space $T_I \Hin$ at a closed point $I \in
\Hin(\C)$ is canonically isomorphic to the vector space
$\text{Hom}_S(I, S/I)$.

\vspace{.2in}

\paragraph{\textbf{2.2.}} The torus $T=\C^* \times \C^*$ acts on $\B^2$ via
$(z_1, z_2) \cdot (x,y)=(z_1x, z_2y)$.
There is an induced action on $S$ given by
$(z_1,z_2) \cdot P(x,y)=P(z_1^{-1}x, z_2^{-1}y)$
and one on $\Hin$ such that
$$(z_1, z_2) \cdot I=\{P(z_1^{-1}x, z_2^{-1}y); P(x,y) \in I\},
\quad\forall I\in\Hin(\C).$$
This action has a finite number of isolated fixed points, indexed by the set of partitions of the integer $n$.
To such a partition $\lambda \vdash n$ corresponds the fixed point $I_{\lambda}$ where
$$I_{\lambda}=\bigoplus_{s\not\in \lambda} \C x^{i(s)-1}y^{j(s)-1}.$$
When $I=I_{\lambda}$ is a $T$-fixed
point, there is an induced $T$-action on $T_I \Hin$.
In order to describe this action, we fix a few notations concerning
$T$. Let $R=R_T$ denote the complexified representation ring of $T$.
We have $R=\C[q^{\pm 1}, t^{\pm 1}]$ where
$$q(z_1, z_2)=z_1^{-1}, \qquad t(z_1, z_2)=z_2^{-1}.$$
For $V$ a $T$-module let $[V]$ be its class in $R$. We
abbreviate  $T_{\lambda}=[T_{I_{\lambda}}\Hin]$. It is given by
\begin{equation}\label{E:21}
T_{{\lambda}}=\sum_{s \in \lambda}(t^{l(s)}q^{-a(s)-1}+t^{-l(s)-1}q^{a(s)}).
\end{equation}

\vspace{.2in}

\paragraph{\textbf{2.3.}}
Let $\Theta_n \subset \Hin \times \B^2$ be the universal family and
let $p: \Hin \times\B^2 \to \Hin$ be the projection. The
\textit{tautological bundle} of $\Hin$ is the locally free sheaf
$\boldsymbol{\tau}_{n}=p_*(\mathcal{O}_{\Theta_n})$. The fiber of
$\boldsymbol{\tau}_n$ at a point $I \in \Hin(\C)$ is $S/I$. The
character of the $T$-action on its fiber at the fixed point
$I_{\lambda}$ is
\begin{equation}\label{E:23}
\boldsymbol{\tau}_\lambda=[\boldsymbol{\tau}_{{n}}|_{I_{\lambda}}]=
\sum_{s \in \lambda} t^{y(s)}q^{x(s)}.
\end{equation}

\vspace{.2in}

\paragraph{\textbf{2.4.}} Let $k \geqslant 0$.
The nested Hilbert scheme $Z_{n,n+k}$ is the reduced closed
subscheme of $\Hin \times \Hi_{n+k}$ parametrizing pairs of ideals
$(I,J)$ where $J \subset I$. One defines the nested Hilbert scheme
$\Hi_{n+k,n}$ in a similar fashion. Of course $Z_{n,n}$ is simply
the diagonal of $\Hin \times \Hin$. The schemes $Z_{n,n+k}$ are
singular in general for $k \geqslant 2$, but they are smooth if
$k=0$ or $k=1$, see \cite{Cheah}. The tangent space at a point
$(I,J) \in Z_{n,n+k}$ is the kernel of the natural map
\begin{equation}\label{E:24}
\psi: \text{Hom}_S(I,S/I) \oplus \text{Hom}_S(J,S/J) \to \text{Hom}_S(J,S/I).
\end{equation}
When $k=1$ the map $\psi$ is surjective.
The diagonal $T$-action on $\Hin \times \Hi_{n+k}$ preserves
$Z_{n,n+k}$. The fixed points contained in $Z_{n,n+k}$ are those
pairs $I_{\mu,\lambda}=(I_{\mu},I_{\lambda})$ for which $\mu \subset
\lambda$. When $k=1$ this may be used to give a cell decomposition
of $\Hi_{n,n+1}$, but we won't need this here. We abbreviate
$T_{\mu,\lambda} =[T_{I_{\mu,\lambda}}\Hi_{n,n+1}]$. We may use
(\ref{E:24}) to write a formula for $T_{\mu,\lambda}$
\begin{equation}\label{E:25}
\begin{split}
T_{{\mu},{\lambda}}&= [\text{Hom}_S(I_\mu,S/I_\mu)] +
[\text{Hom}_S(I_\lambda,S/I_\lambda)] -
[\text{Hom}_S(I_\lambda,S/I_\mu)]\\
&=T_{{\mu}}+T_{{\lambda}}-N_{\mu,\lambda},
\end{split}
\end{equation}
where the character of the fiber of the normal bundle
$N_{\mu,\lambda}=[N_{I_{\mu,\lambda}}\Hi_{n,n+1}]$ is equal to
\begin{equation}\label{E:26}
\begin{split}
N_{{\mu},{\lambda}}=\sum_{s \in \mu} \left(
t^{l_{\mu}(s)}q^{-a_{\lambda}(s)-1}+t^{-l_{\lambda}(s)-1}q^{a_{\mu}(s)}\right).
\end{split}
\end{equation}
Of course, similar formulas hold for the nested Hilbert scheme
$Z_{n+1,n}$.

\vspace{.2in}

\paragraph{\textbf{2.5.}}
Let $\pi_1,\pi_2$ be the natural projections of $\Hin \times
\text{Hilb}_{n+1}$ to $\Hin$ and $\text{Hilb}_{n+1}$ respectively.
Over $Z_{n,n+1}$ there is a natural surjective map $\phi: \pi_2^*(
\boldsymbol{\tau}_{n+1}) \to \pi_1^*( \boldsymbol{\tau}_{n})$. Over
the point $(I,J) \in Z_{n,n+1}$ it specializes to the map $S/J \to
S/I$. The kernel sheaf $\mathcal{K}er(\phi)$ is  a line bundle,
which we call the \textit{tautological bundle} of $Z_{n,n+1}$ and
which we denote by $\boldsymbol{\tau}_{n,n+1}$. Over a $T$-fixed
point $I_{\mu,\lambda}$ its character is
\begin{equation}\label{E:27}
{\boldsymbol\tau}_{\mu,\lambda}=
[\boldsymbol{\tau}_{n,n+1}|_{I_{\mu,\lambda}}]=t^{y(s)}q^{x(s)}
\end{equation}
where $s=\lambda \backslash \mu$ is the unique box of $\lambda$ not
contained in $\mu$. For $l>0$ we write
$$\boldsymbol{\tau}^l_{n,n+1}=
[(\boldsymbol{\tau}_{n,n+1})^{\otimes l}],
\qquad
\boldsymbol{\tau}^{-l}_{n,n+1}=
[(\boldsymbol{\tau}^*_{n,n+1})^{\otimes l}].$$

\vspace{.2in}

\paragraph{\textbf{2.6.}}
Let $\pi_1,\pi_2$ be the natural projections of $\Hin \times \Hin$
to $\Hin$. Over $Z_{n,n}$ we have the vector bundle
$\boldsymbol{\tau}_{n,n} =\pi_2^*(\boldsymbol{\tau}_{n})=\pi_1^*(
\boldsymbol{\tau}_{n})$. We call it the \textit{tautological bundle}
of $\Hi_{n,n}$. Over a $T$-fixed point $I_{\lambda,\lambda}$ its
character is
${\boldsymbol\tau}_{\lambda,\lambda}={\boldsymbol\tau}_{\lambda}.$
As above, we abbreviate
$\boldsymbol{\tau}_{n,n}=[\boldsymbol{\tau}_{n,n}].$

\vspace{.2in}

\section{Convolution algebras in K-theory}

\vspace{.2in}

\paragraph{\textbf{3.1.}} Let $G$ be a complex linear algebraic
group. By a $G$-variety we'll mean a quasi-projective complex
variety with an action of the algebraic group $G$. Let $K^G(X)$ be
the complexified Grothendieck group of the category of
$G$-equivariant coherent sheaves on $X$. It is a $R_G$-module, where
$R_G$ is the representation ring of $G$. We will usually denote the
class of a $G$-equivariant sheaf $\mathcal{F}$ by the symbol
$[\mathcal F]$. If $X' \subset X$ is a $G$-stable closed subvariety
let $[X'] \in K^G(X)$ stand for the class of the structure sheaf of
$X'$. If $X'=\{x\}$, a closed point of $X$, we abbreviate
$[x]=[\{x\}]$.

\vspace{.1in}

Let $X_1,X_2,X_3$ be smooth quasi-projective $G$-varieties admitting
\textit{proper} maps to a (possibly singular) quasi-projective
$G$-variety $Y$. Let $p_{12}, p_{13}, p_{23}$ be the projections
from the triple fiber product $X_1 \times_Y X_2 \times_Y X_3$ along
the factor not named. There is a natural map
\begin{equation*}
\begin{split}
\star~: K^G(X_1\times_Y X_2) \otimes K^G(X_2\times_Y X_3) \to
K^G(X_1\times_Y X_3),\ (w,z) \mapsto Rp_{13*}\left( p_{12}^*(w)
\otimes^{L} p_{23}^*(z)\right).
 \end{split}
 \end{equation*}

If $X_1=X_2=X_3=X$ this map endows $K^G(X \times_Y X)$ with the
structure of an associative $R_G$-algebra with unit given by
$[\Delta_X]$ where $\Delta_X \subseteq X \times_Y X$ is the
diagonal.  Now assume in addition that $Y=\{pt\}$. The map $\star$
with $X_1=X_2=X$ and $X_3=\{pt\}$ endows $K^G(X)$ with the structure
of a $K^G(X \times X)$-module.

\vspace{.2in}

\paragraph{\textbf{3.2.}} We apply this formalism to
$X$ equal to the union of Hilbert schemes $\Hi=\bigsqcup_{n\geqslant
0}\Hin$, and $G=T$, $Y=\{pt\}$. This is not directly possible since
$\Hin$ is not proper if $n>0$. We may get around this difficult
however by first localizing to the $T$-fixed loci, which consists
(for each $n$) of a finite number of points.
For any $R$-module $M=M_R$ we write $M_{\Kb}=M_R \otimes_R \Kb$.
Recall that $\Kb$ is identified with $\C(q^{1/2}, t^{1/2})$ by (\ref{E:qHtH}).
The direct image provides us with an isomorphism
\begin{equation}\label{E:32}
i_*: K^T(\Hin^T) \otimes_{R} \Kb = \bigoplus_{\lambda \vdash n}
\Kb[I_{\lambda}] \to K^T(\Hin)
\otimes_{R}  {\Kb}
\end{equation}
where $i: \Hin^T \to \Hin$ is the embedding.
We have also
$$K^T(\Hin \times \text{Hilb}_{m} )\otimes_{R}\Kb=
\bigoplus_{\substack{\lambda \vdash n\\ \mu \vdash m}}
\Kb[I_{\lambda,\mu}].$$ By (\ref{E:32}) each element in
${K}^T(\Hin)_\Kb$ or ${K}^T(\Hin \times \text{Hilb}_{m})_\Kb$ is a
linear combination of classes of coherent sheaves with proper
support. This allows us to define convolution operations

$$\star~:{K}^T(\Hin \times \text{Hilb}_m)_\Kb \otimes_\Kb
{K}^T(\text{Hilb}_m \times \text{Hilb}_k)_\Kb \to
{K}^T(\text{Hilb}_n \times \text{Hilb}_k)_\Kb$$ and
$$\star~:{K}^T(\Hin \times \text{Hilb}_m)_\Kb
\otimes_\Kb {K}^T(\text{Hilb}_m)_\Kb \to {K}^T(
\text{Hilb}_n)_\Kb.$$ We consider the following associative
$\Kb$-algebra
\begin{equation}\label{E:33}
{\mathbf{E}}_\Kb=\bigoplus_{k \in \Z}\prod_n {K}^T(\Hi_{n+k} \times
\Hin)_\Kb
\end{equation}
where the product ranges over all integers $n$ for which $n
\geqslant 0$ and $n+k \geqslant 0$. It acts on the $\Kb$-vector space
\begin{equation}\label{E:34}
\mathbf{L}_\Kb=\bigoplus_{n \geqslant 0}\mathbf L^n_\Kb,\qquad
\mathbf L^n_\Kb={K}^T(\Hin)_\Kb.
\end{equation}
We will denote this representation as $\psi$. The integer $k$ in
(\ref{E:33}) yields a $\Z$-grading on the $\Kb$-algebra
${\mathbf{E}}_\Kb$ such that the action of
${\mathbf{E}}_\Kb$ on $\mathbf{L}_\Kb$ is compatible with the
gradings. The representation $\psi$ is faithful. Note that we write also
\begin{equation}
\mathbf{L}_R=\bigoplus_{n \geqslant 0}\mathbf L^n_R,\qquad
\mathbf L^n_R={K}^T(\Hin).
\end{equation}

\vspace{.2in}

\paragraph{\textbf{3.3.}}
The $\Kb$-vector space $\mathbf{L}_\Kb$ is spanned by the elements
$[I_{\lambda}]$, $\lambda \in \Pi$. Following \cite{Haiman} we may identify
$\mathbf{L}_\Kb$ with $\Lb_\Kb$ via the map
\begin{align}\label{E:identi}
\mathbf{L}_\Kb \to \Lb_\Kb, \
[I_{\lambda}] \mapsto \widetilde{H}_{\lambda}(q,t).
\end{align}
Here $\widetilde{H}_{\lambda}(q,t)$ is the \textit{cocharge}
Macdonald polynomial. It is a renormalization of Macdonald's
original polynomials $P_{\lambda}(q,t)$. The precise relation is as
follows. Let $\gamma_{t}: \Lb_\Kb \to \Lb_\Kb$ be the unique $\Kb$-algebra
homomorphism such that
$\gamma_{t}(p_r)=(1-t^{r})p_r$.
Then
\begin{equation}\label{E:Macdonald}
t^{n(\lambda)}c_\lambda(q,t^{-1}) \cdot  P_{\lambda}(q,t^{-1})=  \gamma_t (\widetilde{H}_{\lambda}(q,t))
\end{equation}
where
$$n(\lambda)=\sum_i (i-1)\lambda_i, \qquad c_{\lambda}(q,t)=\prod_{s \in \lambda} (1-q^{a(s)}t^{l(s)+1}).$$
Under the identification (\ref{E:identi}) elements of
${\mathbf{E}}_\Kb$ appear as $q,t$-operators on symmetric functions.
For any partitions $\mu, \lambda$ and any operator $z$ on
$\mathbf{L}_\Kb$ we denote by $\langle \mu, z \cdot \lambda\rangle$
the coefficient of $\widetilde{H}_{\mu}(q,t)$ in
$z(\widetilde{H}_{\lambda}(q,t))$.

\vspace{.2in}

\paragraph{\textbf{3.4.}} We would ideally like to understand
the action on $\mathbf{L}_\Kb$ of the classes of the nested Hilbert
schemes $Z_{n,n+k}$ and of their tautological bundles, and
abstractly describe the algebra which these generate. It is not
clear, however, how to deal with the classes of $Z_{n,n+k}$ when
$|k|>1$ since these varieties are singular. Consider the
$\Kb$-subalgebra $\H_\Kb\subset{\mathbf{E}}_\Kb$ generated by
$$\begin{gathered}
\f_{-1,l}=\prod_n\boldsymbol{\tau}^l_{n,n+1},
\quad\f_{1,l}=\prod_n\boldsymbol{\tau}^l_{n+1,n},\quad l\in\Z,\cr
\mathbf{e}_{0,l}=\prod_n \Lambda^l \boldsymbol{\tau}_{n,n},
\quad \mathbf{e}_{0,-l}=\prod_n \Lambda^{l}
\boldsymbol{\tau}^*_{n,n},\quad l\in\Z_{>0}.
\end{gathered}$$
It acts on the $\Kb$-vector space $\mathbf L_\Kb\simeq\Lb_\Kb$ in the
obvious way. Note that here we have abbreviated
$\Lambda^l\boldsymbol{\tau}_{n,n}=[\Lambda^l\boldsymbol{\tau}_{n,n}]$
and
$\Lambda^l\boldsymbol{\tau}_{n,n}^*=[\Lambda^l\boldsymbol{\tau}_{n,n}^*]$.
Here the brackets denote the classes in the Grothendieck groups
oh $\Hin$, $\Hi_{n+1}\times\Hin$ and
$\Hin\times\Hi_{n+1}$ respectively. For a
future use, we define another class of elements $\f_{0,l}\in\mathbf
E_\Kb$ for $l \in \Z^*$ through the relations
$$\sum_{k \geqslant 1}\f_{0,\pm k} s^{k-1}=-\frac{d}{ds}\log ( E_{\pm}(s)),
\quad
E_{\pm}(s)=1+\sum_{k \geqslant 1}(-1)^{k} \mathbf{e}_{0,\pm k}s^k.$$
So $\f_{0,l}$ is obtained from
the classes of the tautological bundles $\boldsymbol{\tau}_{n,n}$ by
the Adams operations
$$\f_{0,k}=\prod_n\Psi_k(\boldsymbol{\tau}_{n,n}),\qquad
\f_{0,-k}=\prod_n\Psi_k(\boldsymbol{\tau}_{n,n}^*),\qquad
\Psi_{-k}(\boldsymbol{\tau}_{n,n})=\Psi_k(\boldsymbol{\tau}_{n,n}^*),
\qquad k\geqslant 1.$$
We can now state our first result. Recall the algebra $\U_c$
associated with the central charge $$c=(1, q^{1/2}t^{1/2}).$$ For
$l\in\Z$ and $n\in\N$ we consider the following elements of ${\U}_c$
\begin{equation}\label{E:hrd1}
\hb_{1,l}=t^{1/2}\f_{1, l-1}, \qquad \hb_{-1,k}=-q^{1/2}\f_{-1, k},
\end{equation}
\begin{equation}\label{E:hrd2}
\hb_{0,k}= \f_{0,k}-\frac{1}{(1-q^k)(1-t^k)},
\qquad
\hb_{0,-k}=-\f_{0,-k}+\frac{1}{(1-q^{-k})(1-t^{-k})},\qquad k\geqslant 1.
\end{equation}
The following  is proved in the next section.

\begin{theo}\label{T:Main}  There is an isomorphism of $\Kb$-algebras
\begin{equation}
\Omega~:\ {\U}_c\to\H_\Kb,\ \bt_{i,l}\mapsto\hb_{i,l},\ \forall
i=-1,0,1,\,l \in \Z.
\end{equation}
\end{theo}

\vspace{.2in}

\section{Construction of the isomorphism}

\vspace{.1in}

Our proof of Theorem~\ref{T:Main} is based on the characterization
of $\widehat{\U}$ given in Proposition~\ref{P:charac} and on the
faithful polynomial representation $\varphi$ of $\U^{\geqslant}$
given in Section~1.4.

\vspace{.2in}

\paragraph{\textbf{4.1.}} We begin by computing the action of the class
$[I_{\lambda,\mu}]$ on $K^T(\text{Hilb}_m)$. For any $T$-equivariant
vector bundle $\mathcal{V}$ on a smooth $T$-variety $X$ we set
$$\Lambda[\mathcal{V}]=\sum_{i\geqslant 0}
(-1)^i[\Lambda^i(\mathcal{V})],$$ where $\Lambda^i$ is the usual wedge
power. This operation descends to a well-defined morphism
$$\Lambda: K^T(X) \to K^T(X)$$
which satisfies $\Lambda(x +y)=\Lambda(x) \cdot \Lambda(y)$. Recall
that if $\mathcal{V}$ is a $T$-equivariant vector bundle on $X$ and
$x \in X$ a $T$-fixed point then $[\mathcal{V}|_{x}]\in R$ is the
character of the fiber of $\mathcal{V}$ over $x$, as a $T$-module.
Again, this descends to the Grothendieck group, yielding the
pullback morphism $K^T(X) \to K^T(\{x\})=R$. The following lemma is
well-known.

\vspace{.1in}

\begin{lem}
(a) For partitions $\lambda \vdash n$ and $\mu,\nu \vdash m$ we have
$$[I_{\lambda,\mu}] \star [I_{\nu}]=\delta_{\mu,\nu} \cdot
\Lambda{(T^*_\nu)}\cdot [I_{\lambda}].$$

(b) For any $T$-equivariant vector bundle $\mathcal{V}$ on $\Hin$ we
have in $K^T(\Hin)$
\begin{equation}
[\mathcal{V}]=\sum_\lambda \Lambda(T^*_\lambda)^{-1}\cdot
[\mathcal{V}|_{I_\lambda}]\cdot[I_\lambda].
\end{equation}
\end{lem}

\vspace{.1in}

For any partitions $\lambda,\mu$ we abbreviate
$\Lambda_{\lambda}=\Lambda(T^*_{\lambda})$ and
$\Lambda_{\lambda,\mu}=\Lambda(T^*_{\lambda})\boxtimes \Lambda(T^*_\mu)$.
As a corollary of the above lemma, we get the following

\vspace{.1in}

\begin{cor} For any $r \in \Z$ it holds
$$[\boldsymbol{\tau}^r_{n,n+1}]=\sum_{\substack{\mu \subset \lambda}} \boldsymbol{\tau}^r_{\mu,\lambda}
\cdot \Lambda(N^*_{\mu,\lambda}) \cdot \Lambda^{-1}_{\mu,\lambda}
\cdot [I_{\mu,\lambda}],$$
$$[\boldsymbol{\tau}^r_{n+1,n}]=\sum_{\substack{\mu \subset \lambda}} \boldsymbol{\tau}^r_{\lambda,\mu}
\cdot \Lambda(N^*_{\lambda,\mu}) \cdot \Lambda^{-1}_{\lambda,\mu}
\cdot [I_{\lambda,\mu}],$$
$$[\boldsymbol{\tau}^r_{n,n}]=\sum_{\lambda \vdash n}\boldsymbol{\tau}^r_{\lambda,\lambda}
\cdot \Lambda^{-1}_{\lambda}\cdot [I_{\lambda,\lambda}],$$ where the
sums range over all pairs of partitions $\lambda \subset \mu$ such
that $|\mu|=n$, $|\lambda|=n+1$.
\end{cor}

\vspace{.2in}

\paragraph{\textbf{4.2.}}
We now write explicit formulas for the action of the elements
$\f_{\pm 1,d}$, $\f_{0,d}$ on 
$\mathbf \Lambda_\Kb$. By (\ref{E:21}) we have
$$\Lambda_{\lambda}=
\prod_{s \in \lambda} (1-t^{-l(s)}q^{a(s)+1})(1-t^{l(s)+1}q^{-a(s)}).$$
For $\mu\subset\lambda$, by (\ref{E:26}) we have
$$\Lambda(N^*_{\mu,\lambda})=\prod_{s \in \mu} (1-t^{-l_{\mu}(s)}q^{a_{\lambda}(s)+1})
(1-t^{l_{\lambda}(s)+1}q^{-a_{\mu}(s)})
=\Lambda(N^*_{\lambda,\mu}).$$ We have also
$$\boldsymbol{\tau}^r_{\mu,\lambda}=\boldsymbol{\tau}^r_{\lambda,\mu}=t^{r\cdot y(\lambda\backslash \mu)}q^{r \cdot x(\lambda\backslash
\mu)},\quad \boldsymbol{\tau}^r_{\mu,\mu}=\big(\sum_{s \in
\mu}t^{ y(s)}q^{ x(s)}\big)^r.$$ Using the identification
of Section~3.3, we deduce after straightforward computations

\vspace{.1in}

\begin{lem}\label{L:gyt} (a) For any $r \in \Z$ and any partition $\nu$ we have
\begin{equation*}
\f_{1,r} \cdot \widetilde{H}_{\nu}(q,t)=\frac{1}{(1-q)(1-t)}\sum_{\mu \supset \nu}
q^{(r+1) x(\mu \backslash \nu)}t^{(r+1) y(\mu \backslash \nu)}
L_{\nu,\mu}(q,t)\widetilde{H}_{\mu}(q,t),
\end{equation*}
where the sum ranges over all $\mu \supset \nu$ with
$|\mu|=|\nu|+1$, and
\begin{equation*}
L_{\nu,\mu}(q,t)=\prod_{s \in C_{\mu \backslash
\nu}}\frac{
t^{l_{\nu}(s)}-q^{a_{\nu}(s)+1}}{t^{l_{\nu}(s)+1}-q^{a_{\nu}(s)+1}}
\cdot \prod_{s \in R_{\mu \backslash \nu}}\frac{
t^{l_{\nu}(s)+1}-q^{a_{\nu}(s)}}{t^{l_{\nu}(s)+1}-q^{a_{\nu}(s)+1}}.
\end{equation*}
(b) For any $r \in \Z$ and any partition $\nu$ we have
\begin{equation*}
\f_{-1,r} \cdot \widetilde{H}_{\nu}(q,t)=\sum_{\lambda \subset \nu}
q^{r x(\nu \backslash \lambda)}t^{r y(\nu \backslash
\lambda)}L_{\nu,\lambda}(q,t)\widetilde{H}_{\lambda}(q,t),
\end{equation*}
where the sum ranges over all $\lambda \subset \nu$ with
$|\lambda|=|\nu|-1$, and
\begin{equation*}
L_{\nu,\lambda}(q,t)=\prod_{s \in C_{\nu \backslash \lambda}} \frac{
t^{l_{\nu}(s)+1}-q^{a_{\nu}(s)}}{t^{l_{\nu}(s)}-q^{a_{\nu}(s)}}
\cdot \prod_{s \in R_{\nu \backslash \lambda}} \frac{
t^{l_{\nu}(s)}-q^{a_{\nu}(s)+1}}{t^{l_{\nu}(s)}-q^{a_{\nu}(s)}}.
\end{equation*}
\end{lem}

\vspace{.1in}

In particular we have the following formulas, see \cite{Garsia}, Theorem~1.4.

\vspace{.1in}

\begin{cor}\label{C:frtyue} The following relations holds in
$\mathrm{End}_\Kb(\mathbf\Lambda_\Kb)$
$$\f_{1,-1}=\frac{1}{(1-q)(1-t)} p_1, \qquad \f_{-1,0}=\frac{\partial}{\partial p_1}.$$

\end{cor}

\vspace{.1in}

The action of $\f_{0,r}$ is given by the following formula,
compare \cite{Mac}, I. (2.10').

\vspace{.1in}

\begin{lem}\label{L:hytrfdx} For any $r \in \Z^*$ and any partition
$\nu$ we have
\begin{equation*}
\f_{0,r} \cdot \widetilde{H}_{\nu}(q,t)=\bigl(\sum_{s \in \nu} q^{r
x(s)}t^{r y(s)}\bigr) \widetilde{H}_{\nu}(q,t).
\end{equation*}
\end{lem}

\vspace{.2in}

\paragraph{\textbf{4.3.}} Next, we check
that the generators $\f_{\pm 1, r}$, $\f_{0,l}$ of $\H_\Kb$ satisfy
relations (\ref{E:2}), (\ref{E:205}).

\vspace{.1in}

\begin{prop}\label{P:klq1} For any $l, k \in \Z^*$ we have $[\f_{0,l},\f_{0,k}] =0$.
\end{prop}
\begin{proof} This is obvious since the convolution product
of two classes supported on the diagonal in $\Hi\times\Hi$ is their
tensor product.\end{proof}

\vspace{.1in}

\begin{prop}\label{P:klq2} For any $l  \in \Z^*$, $k \in \Z$ we have
$[\f_{0,l}, \f_{\pm 1, k}]=\pm \f_{\pm 1, k+l}.$
\end{prop}
\begin{proof} For any partition $\lambda$ we set
$B_{\lambda}^l(q,t)=\sum_{s \in \lambda} q^{l x(s)}t^{l y(s)}.$ Then
we have
$$\f_{0,l} \f_{1,k} \cdot \widetilde{H}_{\nu}(q,t)=\frac{1}{(1-q)(1-t)}\sum_{\mu \supset \nu}
B^l_{\mu}(q,t)q^{(r+1) x(\mu \backslash \nu)}t^{(r+1) y(\mu \backslash \nu)}
L_{\nu,\mu}(q,t)\widetilde{H}_{\mu}(q,t),$$
while
$$\f_{1,k} \f_{0,l} \cdot \widetilde{H}_{\nu}(q,t)=\frac{1}{(1-q)(1-t)}\sum_{\mu \supset \nu}
B^l_{\nu}(q,t)q^{(r+1) x(\mu \backslash \nu)}t^{(r+1) y(\mu \backslash \nu)}
L_{\nu,\mu}(q,t)\widetilde{H}_{\mu}(q,t).$$
Substituting the relation
$$B^l_{\mu}(q,t)-B^l_{\nu}(q,t)=q^{lx(s)}t^{ly(s)}$$
for $s=\mu\backslash \nu$ we obtain
$[\f_{0,l},\f_{1,k}]=\f_{1,l+k}$. The second part of the proposition
is identical. \end{proof}

\vspace{.2in}

\paragraph{\textbf{4.4.}} Using Corollary~\ref{C:frtyue}
and Lemma~\ref{L:hytrfdx} we may now compare the representation
$$\psi:\H_\Kb\to\End_\Kb(\mathbf{L}_\Kb)$$
of Section 3.2 with the representation
$$\widetilde{\varphi}: \U^{\geqslant}\to\End_\Kb(\Lb_\Kb)$$
of Corollary 1.5, under the isomorphism  $\mathbf
L_\Kb\simeq\mathbf\Lambda_\Kb$ in (\ref{E:identi}). Recall the
plethystic substitution $\gamma_t \in \mathrm{End}_{\Kb}(\Lb_{\Kb})$
of Section~3.3. By (\ref{E:Macdonald}) we have
$$\gamma_t(\widetilde{H}_{\lambda}(q,t)) =
u_{\lambda}(q,t)
P_{\lambda}(q,t^{-1})$$ for some scalar
$u_{\lambda}(q,t)$. Using
Corollary~\ref{C:pol} and Lemma~\ref{L:hytrfdx} we get, for each $n
\geqslant 1$
\begin{equation*}
\begin{split}
\widetilde{\varphi}(\bt_{0,n}) \gamma_t(\widetilde{H}_{\lambda}(q,t)) &=
\big( \sum_{s \in \lambda}
q^{nx(s)}t^{ny(s)}\big)
u_{\lambda}(q,t) P_{\lambda}(q,t^{-1})\\
&=\big( \sum_{s \in \lambda} q^{nx(s)}t^{ny(s)}\big) \gamma_t(\widetilde{H}_{\lambda}(q,t))\\
&=\gamma_t(\psi(\f_{0,n})\widetilde{H}_{\lambda}(q,t))
\end{split}
\end{equation*}
and
\begin{equation*}
\begin{split}
\widetilde{\varphi}(\bt_{0,-n}) \gamma_t(\widetilde{H}_{\lambda}(q,t))
&=-\big( \sum_{s \in \lambda} q^{-nx(s)}t^{-ny(s)}\big) u_{\lambda}(q,t) P_{\lambda}(q,t^{-1})\\
&=-\big( \sum_{s \in \lambda} q^{-nx(s)}t^{-ny(s)}\big) \gamma_t(\widetilde{H}_{\lambda}(q,t))\\
&=-\gamma_t (\psi(\f_{0-,n}) \widetilde{H}_{\lambda}(q,t)).
\end{split}
\end{equation*}
In addition, we have
$$\widetilde{\varphi}(\bt_{1,0}) \circ \gamma_t = 
\frac{p_1}{(q-1)} \circ \gamma_t = 
\gamma_t \circ \frac{p_1}{(q-1)(1-t)}=-\gamma_t \circ \psi(\f_{1,-1}).$$
Now, let $\H^{\geqslant}_\Kb\subset\H_\Kb$ be the subalgebra generated by
$\f_{1,-1}$ and the elements $\f_{0,n}$ for $n \in \Z^*$. By
Proposition~\ref{P:klq2} the elements $\f_{1,n}$ belong to
$\H^{\geqslant}_\Kb$ for all $n \in \Z$.
Because the representations $\psi$, $\widetilde{\varphi}$ are both faithful and
because the $\Kb$-algebras
$\H^{\geqslant}_\Kb$ and $\U^{\geqslant}$ are respectively
generated by $\f_{1,-1}$, $\f_{0,l}$, $l \in \Z^*$ and $\bt_{1,0}$,
$\bt_{0,l}$, $l \in \Z^*$ we deduce from the above formulas that the assignment
$$\bt_{1,0} \mapsto -\f_{1,-1}, \qquad \bt_{0,n} \mapsto \f_{0,n},
\qquad \bt_{0,-n} \mapsto -\f_{0,-n}, \quad n \geqslant 1$$
extends to an isomorphism $\U^{\geqslant} \to\H^{\geqslant}_\Kb$.
Twisting by automorphisms
(\ref{E:autom1}), (\ref{E:autom2}) we see that the assignment
$$\bt_{1,0} \mapsto t^{1/2} \f_{1,-1}, \qquad
\bt_{0,n} \mapsto \f_{0,n} -\frac{1}{(1-q^n)(1-t^n)}, \qquad
\bt_{0,-n} \mapsto -\f_{0,-n}+\frac{1}{(1-q^{-n})(1-t^{-n})}$$
also extends to an isomorphism of algebras
$\U^{\geqslant}\to\H^{\geqslant}_\Kb$.
In other words, restricting
the map $\Omega$ in Theorem~\ref{T:Main} we get an isomorphism
$\U^{\geqslant}\to\H^{\geqslant}_\Kb$.
In the same way we prove that
$\Omega$ restricts also to an isomorphism
$\U^{\leqslant}
\to\H^{\leqslant}_\Kb$.

\vspace{.2in}

\paragraph{\textbf{4.5.}}
We have just proved that $\Omega$ restricts to $\Kb$-algebra
isomorphisms 
$$\U^{\geqslant}\to\H^{\geqslant}_\Kb,\qquad
\U^{\leqslant}\to\H^{\leqslant}_\Kb.$$
By Proposition~\ref{P:charac}
these two morphisms extend to the whole of $\U_c$ if and only if
relation (\ref{E:5}) holds with the appropriate specialization of
the center. The proof of this fact, which requires somewhat involved
computations, is given in Appendix A.1.

\vspace{.2in}

\paragraph{\textbf{4.6.}}
Let $\H_\Kb^>, \H^0_{\Kb},\H_\Kb^<\subset\H_\Kb$ be the subalgebras
generated by the elements $\f_\x$ with $\x\in\Zb^1, \Zb^0, \Zb^{-1}$
respectively.  By Sections~4.4 and 4.5 there is a well-defined
surjective algebra homomorphism $\Omega: \U_c\to \H_\Kb$ which
restricts to isomorphisms 
$$\U^>\to\H_\Kb^>,\qquad
\U^0\to\H_\Kb^0,\qquad
\U^< \to \H_\Kb^<.$$ Recall the triangular decomposition $\U_c\simeq
\U^> \otimes \U^0 \otimes \U^< $. To prove that the map
$\Omega$ is injective it is enough to prove the following result.

\vspace{.1in}

\begin{prop}\label{T:2} The algebra $\H_\Kb$ has a triangular decomposition,
i.e., the multiplication map induces an isomorphism
$m~: \H_\Kb^> \otimes \H_\Kb^0 \otimes \H_\Kb^< \to \H_\Kb.$
\end{prop}

\begin{proof} Since $\Omega$ is surjective and $\U_c$ has a
triangular decomposition, the multiplication map
$m$ is onto.
We will now prove that it is also injective.
The basic idea is to mimic the construction of a
coproduct on $\H_\Kb$ \footnote{The existence of this coproduct is a
consequence of the identification $\H_\Kb\simeq\U$.}.
Let us argue by contradiction and let $$x=\sum_i P_i \otimes R_i
\otimes Q_i$$ be a nonzero homogeneous element in $\mathrm{Ker}(m)$. We
may assume that the $R_i=R_i(\f_{0,\pm 1}, \f_{0,\pm 2}, \ldots)$
are linearly independent polynomials in the variables $\f_{0,l}$, $l
\in \Z$, and that $P_i$ and $Q_i$ are nonzero. Multiplying by an
element of $\H_\Kb^>$ or $\H_\Kb^<$ if necessary, we may also assume
that $\deg(x)=0$. By definition, we have
\begin{equation}\label{E:PP4}
\sum_i P_i \circ R_i \circ Q_i \cdot \widetilde{H}_{\lambda}(q,t)=0
\end{equation}
for all partitions $\lambda$. We will apply (\ref{E:PP4}) to a certain (asymptotic) kind of partition.
Given partitions $\lambda_1, \lambda_2, \ldots, \lambda_k$, and given an integer $n \gg |\lambda_1|,\ldots, |\lambda_k|$ we let 
$\lambda_1\circledast\ldots\circledast\lambda_k$ stand for the following partition

\centerline{
\begin{picture}(200,200)
\put(20,170){{$\lambda_1$}}
\put(75,115){{$\lambda_2$}}
\put(185,17){{$\lambda_k$}}
\put(57,157){{\small{$(n,(k-1)n)$}}}
\put(112,102){{\small{$(2n,(k-2)n)$}}}
\put(167,57){{\small{$((k-1)n,n)$}}}
\multiput(0,155)(5,0){5}{\line(1,0){3}}
\put(0,0){\line(0,1){180}}
\put(0,0){\line(1,0){200}}
\put(0,180){\line(1,0){5}}
\put(5,180){\line(0,-1){5}}
\put(5,175){\line(1,0){5}}
\put(10,175){\line(0,-1){5}}
\put(10,170){\line(1,0){5}}
\put(15,170){\line(0,-1){10}}
\put(15,160){\line(1,0){10}}
\put(25,160){\line(0,-1){5}}
\put(25,155){\line(1,0){30}}
\multiput(55,100)(0,5){4}{\line(0,1){3}}
\multiput(55,100)(5,0){5}{\line(1,0){3}}
\put(55,155){\line(0,-1){35}}
\put(55,120){\line(1,0){5}}
\put(60,120){\line(0,-1){5}}
\put(60,115){\line(1,0){5}}
\put(65,115){\line(0,-1){5}}
\put(65,110){\line(1,0){10}}
\put(75,110){\line(0,-1){5}}
\multiput(165,0)(0,5){4}{\line(0,1){3}}
\put(75,105){\line(1,0){5}}
\put(80,105){\line(0,-1){5}}
\put(80,100){\line(1,0){30}}
\multiput(112,78)(5,-5){5}{\circle*{1}}
\put(110,100){\line(0,-1){20}}
\put(135,55){\line(1,0){30}}
\put(165,55){\line(0,-1){35}}
\put(165,20){\line(1,0){5}}
\put(170,20){\line(0,-1){5}}
\put(170,15){\line(1,0){10}}
\put(180,15){\line(0,-1){5}}
\put(180,10){\line(1,0){5}}
\put(185,10){\line(0,-1){5}}
\put(185,5){\line(1,0){15}}
\put(200,5){\line(0,-1){5}}
\end{picture}}
\vspace{.05in}
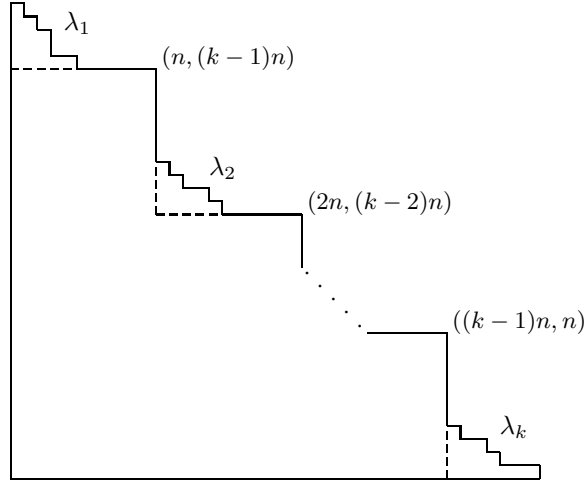
\centerline{\textbf{Figure 2.} An asymptotic partition $\lambda_1\circledast\ldots\circledast\lambda_k$}

\vspace{.2in}

Note that $\lambda_1\circledast\ldots\circledast\lambda_k$ is well-defined as soon as $n > \text{sup}_i(l(\lambda_i),l(\lambda_i'))$.
Put
$$r=\text{sup}_i (\deg(P_i)) =\text{sup}_i(-\deg(Q_i)).$$
Recall that for partitions $\nu, \gamma$ and $z$ an operator on
$\mathbf{L}_\Kb$ we denote by $\langle \nu, z \cdot \gamma\rangle$ the
coefficient of $\widetilde{H}_{\nu}(q,t)$ in $z(\widetilde{H}_{\gamma}(q,t))$.
For $n\gg 0$ we consider the coefficients $$\langle
{\lambda}^{\#}\circledast\gamma\circledast\mu^{\#}, P_i R_i Q_i\cdot
\lambda\circledast\gamma\circledast\mu\rangle,\quad \lambda^{\#} \subset
\lambda,\quad \mu \subset \mu^{\#},\quad |\lambda \backslash
\lambda^{\#}|=|\mu^{\#}\backslash \mu|=r.$$
Since the $Q_i$ are annihilation operators while the $P_i$ are
creation operators, and because $r$ is maximal, the only way to
obtain $\lambda^{\#}\circledast\gamma\circledast\mu^{\#}$ from
$\lambda\circledast\gamma\circledast\mu$ is to use all of $Q_i$ to reduce $\lambda$
to $\lambda^{\#}$ and to use all of $P_i$ to increase $\mu$ to
$\mu^{\#}$. Therefore we have
\begin{equation}\label{E:PP5}
\begin{gathered}
\langle \lambda^{\#}\circledast\gamma\circledast\mu^{\#}, P_i R_i Q_i\cdot (\lambda\circledast\gamma\circledast\mu)\rangle=\cr
=\langle \lambda^{\#}\circledast\gamma\circledast\mu^{\#}, P_i\cdot (\lambda^{\#}\circledast\gamma\circledast\mu)\rangle \,
\langle \lambda^{\#}\circledast\gamma\circledast\mu, R_i \cdot(\lambda^{\#}\circledast\gamma\circledast\mu)\rangle\cr
\langle \lambda^{\#}\circledast\gamma\circledast\mu, Q_i
\cdot(\lambda\circledast\gamma\circledast\mu)\rangle.
\end{gathered}
\end{equation}
Note that (\ref{E:PP5}) is equal to zero unless $\deg(P_i)=-\deg(Q_i)= r$.

\vspace{.05in}

Next, we define automorphisms $\tau \in \text{Aut}(\H_\Kb^<)$, $\rho \in
\text{Aut}(\H_\Kb^>)$ by
$$\tau(\f_{-1,k})=t^k \f_{-1,k},\quad
\rho(\f_{1,l})=q^l\f_{1,l},\quad\forall k,l \in \Z.$$ The existence of
$\tau, \rho$ is a consequence of the isomorphisms
$$\H_\Kb^> \simeq \U^>,\quad\H_\Kb^< \simeq \U^<.$$

\begin{lem}\label{L:PK1} There are constants
$c$ and
$d$ such that
\begin{equation}\label{E:tauy}
\langle \lambda^{\#}\circledast\gamma\circledast\mu, Q\cdot
(\lambda\circledast\gamma\circledast\mu)\rangle=c \;  \langle \lambda^{\#},
\tau^{2n}(Q) \cdot\lambda \rangle,
\end{equation}
\begin{equation}\label{E:rhoy}
\langle \lambda^{\#}\circledast\gamma\circledast\mu^{\#}, P\cdot
(\lambda^{\#}\circledast\gamma\circledast\mu)\rangle=d \; \langle \mu^{\#},
\rho^{2n}(P) \cdot\mu \rangle
\end{equation}
for any operator $P\in \H_\Kb^>[r]$ and $Q \in \H_\Kb^<[-r]$.
\end{lem}
\begin{proof} We prove the first statement only, the second one is identical.
If $Q=\f_{-1,k_r} \cdots \f_{-1,k_1}$ then
$$\langle \lambda^{\#}\circledast\gamma\circledast\mu, Q\cdot  (\lambda\circledast\gamma\circledast\mu)\rangle=
\sum \biggl(t^{\sum k_j y(s_j)}q^{\sum k_j x(s_j)} \cdot
\prod_{j=1}^r L_{\lambda_j\circledast\gamma\circledast\mu,
\lambda_{j+1}\circledast\gamma\circledast\mu}(q,t)\biggr).$$ Here we have
\begin{itemize}
\item
the sum runs over all sequences
$\lambda =\lambda_1 \supsetneq \lambda_2 \cdots
\supsetneq \lambda_{r+1}=\lambda^{\#}$
and $s_i=\lambda_i\backslash \lambda_{i+1}$,
\item
for partitions $\alpha \supset \beta$ with $|\alpha|=|\beta|+1$ we have
$$L_{\alpha,\beta}(q,t)=\prod_{s \in C_{\alpha \backslash \beta}} \frac{
t^{l_{\alpha}(s)+1}-q^{a_{\alpha}(s)}}{t^{l_{\alpha}(s)}-q^{a_{\alpha}(s)}}
\cdot \prod_{s \in R_{\alpha \backslash \beta}} \frac{
t^{l_{\alpha}(s)}-q^{a_{\alpha}(s)+1}}{t^{l_{\alpha}(s)}-q^{a_{\alpha}(s)}}.
$$
\end{itemize}
For $s_j$ a box in $\lambda_j$ we have
$x(s_j)=x_{\lambda}(s_j)$ and $y(s_j)=y_{\lambda}(s_j)+2n$,
where $x_{\lambda}$ and $y_{\lambda}$ denote the coordinate values
when we place the origin at the
bottom left corner of $\lambda$, i.e., at the point $(0,2n)$,
as opposed to the coordinate values when the origin is
at the bottom left corner of $\lambda\circledast\gamma\circledast\mu$.
Similarly we have, for the row and columns of a box $s_j$ in $\lambda_j$
\begin{align*}
R(s_j)=R_{\lambda_j}(s_j),\qquad
C(s_j)=C_{\lambda_j}(s_j) \sqcup C'(s_j)
\end{align*}
where $C'(s_j)=\{(x(s_j),0), \ldots, (x(s_j),2n-1)\}$.
Finally, observe that the armlength $a(u)$ or the
leglength $l(u)$ of a box $u \in \lambda_j$ are
the same whether we consider $u$ as belonging to $\lambda_j$
or to $\lambda_j\circledast\gamma\circledast\mu$.
From the above formulae we deduce that
\begin{equation*}
\begin{split}
&t^{\sum k_j y(s_j)}q^{\sum k_j x(s_j)} \cdot
\prod_{j=1}^r L_{\lambda_j\circledast\gamma\circledast\mu, \lambda_{j+1}\circledast\gamma\circledast\mu}(q,t)=\\
&=t^{\sum k_j y_{\lambda}(s_j)}q^{\sum k_j x_{\lambda}(s_j)} \cdot
\prod_{j=1}^r L_{\lambda_j,\lambda_{j+1}}(q,t) \cdot t^{2n\sum k_j}
\prod_{j=1}^r\prod_{u \in C'(s_j)}
\frac{t^{l_{\sigma_j}(u)+1}-q^{a_{\sigma_j}(u)}}{t^{l_{\sigma_j}(u)}-q^{a_{\sigma_j}(u)}}
\end{split}
\end{equation*}
where we have set $\sigma_j=\lambda_j\circledast\gamma\circledast\mu$. It remains to note that
$$\prod_{j=1}^r
\prod_{u \in
C'(s_j)}\frac{t^{l_{\lambda_j}(u)+1}-q^{a_{\lambda_j}(u)}}{t^{l_{\lambda_j}(u)}-q^{a_{\lambda_j}(u)}}
=\prod_{s \in \lambda \backslash \lambda^{\#}}\prod_{u\in C'(s)}
\frac{t^{y(s)-y(u)+1}-q^{a_{\sigma}(u)}}{t^{y(s)-y(u)}-q^{a_{\sigma}(u)}},$$
where we have set $\sigma=\lambda\circledast\gamma\circledast\mu$, is independent
of the choice of the chain of subdiagrams $(\lambda_j)_j$, and that
$$\sum\biggl(t^{\sum k_j y_{\lambda}(s_j)}q^{\sum k_j x_{\lambda}(s_j)} \cdot
\prod_{j=1}^r L_{\lambda_j,\lambda_{j+1}}(q,t)\biggr)= \langle
\lambda^{\#}, Q \cdot \lambda \rangle.$$ The lemma is proved.
\end{proof}

\vspace{.1in}

Using (\ref{E:PP5}) together with the above lemma, the linear relation (\ref{E:PP4}) may be rescaled to
\begin{equation}\label{E:PP6}
\sum_i \langle \mu^{\#}, \rho^{2n}(P_i)\cdot \mu\rangle \,
\langle \lambda^{\#}\circledast\gamma\circledast\mu, R_i
\cdot(\lambda^{\#}\circledast\gamma\circledast\mu)\rangle \, \langle
\lambda^{\#},\tau^{2n}( Q_i) \cdot \lambda\rangle=0
\end{equation}
for all $\lambda,\lambda^{\#},\gamma, \mu, \mu^{\#}$ as above, and
all $n \gg 0$. Let us choose some $\mu,\mu^{\#}$ and
$\lambda,\lambda^{\#}$ such that $\langle \mu^{\#},
\tau^{2n}(P_i)\cdot \mu\rangle \neq 0$ and $\langle \lambda^{\#},
\rho^{2n}(Q_i)\cdot \lambda\rangle \neq 0$ for at least one value of
$i$. Let us fix $n \gg 0$ and let us vary $\gamma$. Recall that
$R_i$ is a polynomial in the operators $\f_{0,l}$ and observe that
\begin{equation*}
\begin{split}
\langle &\lambda^{\#}\circledast\gamma\circledast\mu, \f_{0,l} \cdot (\lambda^{\#}\circledast\gamma\circledast\mu)\rangle=\\
&= t^{2nl}\langle \lambda^{\#}, \f_{0,l}\cdot \lambda^{\#}\rangle +
q^{nl}t^{nl} \langle \gamma,\f_{0,l}\cdot \gamma\rangle +q^{2nl}
\langle \mu, \f_{0,l}\cdot \mu\rangle + \langle
\emptyset\circledast\emptyset\circledast\emptyset, \f_{0,l}\cdot
(\emptyset\circledast\emptyset\circledast\emptyset)\rangle.
\end{split}
\end{equation*}
Setting
\begin{equation*}
R'_i= \langle \mu^{\#}, \rho^{2n}(P_i)\cdot \mu\rangle\, \langle
\lambda^{\#},\tau^{2n}( Q_i) \cdot \lambda\rangle \, R_i(t^{\pm
n}q^{\pm n}\f_{0,\pm 1}+\alpha_{\pm 1}, t^{\pm 2n}q^{\pm 2
n}\f_{0,\pm 2}+\alpha_{\pm 2}, \ldots)
\end{equation*}
where
$$\alpha_l=\langle \emptyset\circledast\emptyset\circledast\emptyset, \f_{0,l}\cdot (\emptyset\circledast\emptyset\circledast\emptyset)\rangle +t^{2nl}
\langle \lambda^{\#}, \f_{0,l}\cdot \lambda^{\#}\rangle +q^{2nl}
\langle \mu, \f_{0,l}\cdot \mu\rangle,\quad l \in \Z^*,$$ we may
rewrite (\ref{E:PP6}) as
\begin{equation}\label{E:PP8}
\sum_i \langle \gamma, R'_i\cdot \gamma \rangle =0.
\end{equation}
Since this holds for all $\gamma$ with $l(\gamma), l(\gamma') < n$, taking $n$ large enough we deduce that $\sum_i R'_i  \equiv 0$. Remember that the $R_i$s were chosen to be linearly independent; but then the $R'_i$s are also linearly independent and we arrive at a contradiction. This finishes the proof of
Proposition~\ref{T:2}. Theorem~\ref{T:Main} follows.
\end{proof}

\vspace{.2in}

\paragraph{\textbf{4.7.}}
Theorem~\ref{T:Main} allows us to extend the polynomial
representation $\widetilde{\varphi}$ of $\U^{\geqslant}$ defined in Section~1.4
to the whole algebra $\U_c$. We simply use the isomorphism $\Omega$
to transport the  representation of $\H_\Kb$ on $\Lb_\Kb$ to $\U_c$.
Recall that $\s=q^{-1}$ and $\bs=t^{-1}$.
Recall also the operators $\widetilde{\Delta}^{\infty}_{\pm
l}$ on $\Lb_\Kb$ defined in
Section~1.4 for each $l \geqslant 1$. We set $$\Delta^{\infty}_{\pm
l}={\widetilde{\Delta}^{\infty}_{\pm l}}/({q^{l}-1}).$$ The
eigenvalue of $\Delta^{\infty}_{\pm l}$ on $P_{\lambda}(q,t^{-1})$
is equal to $\pm B^{\pm l}_{\lambda}(q,t)$ by (\ref{E:eig4}).

\vspace{.1in}

\begin{prop}\label{P:repinfty} There is a unique faithful representation
$\varphi$ of $\U_c$ on $\Lb_\Kb$ such that
\begin{align*}
\varphi(\bt_{l,0})&=\frac{t^{l/2}}{1-q^l}\;p_l,\\
\varphi(\bt_{-l,0})&=-l\frac{q^{l/2}}{1-t^l}\;\frac{\partial}{\partial p_l},\\
\varphi(\bt_{0,l})&=\Delta^{\infty}_{ l} - \frac{1}{(1-q^{ l})(1-t^{ l})},\\
\varphi(\bt_{0,-l})&=\Delta^{\infty}_{- l} + \frac{1}{(1-q^{- l})(1-t^{- l})}.
\end{align*}
\end{prop}

\vspace{.1in}

Since the $\Kb$-algebra $\U_c$ is generated by $\bt_{0, \pm
1}$, $\bt_{\pm 1, 0}$ it is enough to specify the action of these
elements to determine the representation $\varphi$. Note that there are no finite rank analogues
of the representation $\varphi$ because the ring $\SH_n$ acts only
on the space of symmetric \textit{Laurent} polynomials $\Kb[x_1^{\pm
1}, \ldots, x_n^{\pm 1}]^{\mathfrak{S}_{\infty}}$.

\vspace{.2in}

\section{Virtual classes and their action on $K^T(\text{Hilb})$}

\vspace{.1in}

We have defined the algebra $\H_\Kb$ as the subalgebra of
$\mathbf{E}_{\Kb}$ generated by the classes of the tautological
bundles on the \textit{smooth} nested Hilbert schemes $Z_{n,m}$,
i.e., when $|n-m|\leqslant 1$. The aim of this section is to show that
$\H_{\Kb}$ contains the \textit{virtual classes} of the more general
(singular) nested Hilbert schemes $Z_{n,m}$ as well as the
tautological bundles over them. We also explicitly describe the action
of these virtual classes in the natural representation
$\boldsymbol{\Lambda}_{\Kb}$.

\vspace{.2in}

\paragraph{\textbf{5.1.}}
Consider the virtual vector bundle $\Vc$ over $\Hi\times\Hi$ with
fiber
$$\Vc|_{(I,J)}=\chi(\mathcal{O})-\chi(I,J).$$
Here $I,J$ are closed points of Hilb which are viewed
as ideal sheaves on $\B^2$ and
$$\chi(\mathcal{F},\mathcal{G})=
\sum_{i=0}^2(-1)^i\text{Ext}^i(\mathcal{F},\mathcal{G})$$
for any coherent sheaves $\mathcal{F}$, $\mathcal{G}$ on $\B^2$. We
abbreviate $V_{\lambda,\mu}=[\Vc|_{I_{\lambda,\mu}}]$, an
element of $R$.

\begin{lem}\label{L:51} For any partitions $\lambda,\mu$ such that
$|\lambda|=n$, $|\mu|=m$ the following hold

(a) $V_{\lambda,\mu}=\sum_{s \in \mu}
t^{l_{\lambda}(s)+1}q^{-a_{\mu}(s)}+\sum_{s\in\lambda}
t^{-l_{\mu}(s)}q^{a_{\lambda}(s)+1}=qtV_{\mu,\lambda}^*,$

(b)  $T_\lambda=V_{\lambda,\lambda}$ if $n=m$,

(c)  $N^*_{\lambda,\mu}=V_{\lambda,\mu}-qt$ if $n=m+1$ and
$\lambda\supset\mu$,

(d)  $N^*_{\lambda,\mu}=qtV^*_{\lambda,\mu}-qt$ if $n=m-1$ and
$\lambda\subset\mu$,

(e) $\Lambda(V_{\lambda,\mu})=0$ unless $\mu\subset\lambda$ and
$\Lambda(qtV^*_{\lambda,\mu})=0$ unless
$\lambda\subset\mu$.

\end{lem}

\begin{proof}
Part $(a)$ is proved in \cite{CO}. Part $(b)$ is obvious. Now,
assume that $n=m+1$ and $\lambda\supset\mu$. Let
$s=\lambda\setminus\mu$. From (\ref{E:24}) and $(a)$ we get the
following formula
$V_{\lambda,\mu}=N^*_{\lambda,\mu}+
t^{l_\lambda(s)+1}q^{-a_\mu(s)}=N^*_{\lambda,\mu}+qt$.
This yields $(c)$. Part $(d)$ follows from $(a)$ and $(c)$.
To prove the first statement in part (e) it suffices to notice
that if $\lambda \not\supset \mu$ then there exists a box $s \in \lambda$
with $l_{\mu}(s)=0, a_{\lambda}(s)=-1$ or a box $s \in \mu$ with
$l_{\lambda}(s)=-1, a_{\mu}(s)=0$ (such a box is located at the intersection
of the right boundaries of the Young diagrams of $\lambda$ and $\mu$).
The second statement of (e) is obtained by duality.
\end{proof}

\paragraph{\textbf{5.2.}}
The associative $R$-algebra
\begin{equation}\label{I:33}
{\mathbf{E}}_R=\prod_{n} \bigoplus_{k \in \Z} {K}^T(\Hi_{n+k} \times
\Hin)
\end{equation}
where the product ranges over all integers $n$ for which $n \geqslant 0$
and $n+k \geqslant 0$, acts on the $R$-module
\begin{equation}\label{I:34}
\mathbf L_R=\bigoplus_{n\geqslant 0}{K}^T(\Hin).
\end{equation}
Abbreviate $\Vcb=[\Vc]$, a class in ${\mathbf E}_R$. Let
$\Vcb_{n,m}$ be the restriction of $\Vcb$ to $K^T(\Hin\times\Hi_m)$ and
consider the elements of $\mathbf E_R$
$$\Vcb_k=\prod_n \Vcb_{k+n,n},\quad \Vcb^*_{-k}=\prod_n
\Vcb_{n,n+k},\quad k>0.$$ By part (e) of the above lemma the classes
$\Lambda(\Vcb_k)$, $\Lambda(qt\Vcb^*_{-k})$ are supported on the union of
nested Hilbert schemes $\bigsqcup_{n,m} Z_{n,m}$. 

\begin{rem} The
class $\Lambda(\Vcb_{n,m})$ for $n >m$ or $\Lambda(qt\Vcb^*_{n,m})$ for
$n<m$ is called the \emph{virtual fundamental class} of
$Z_{n,m}$, see \cite{CO}. The reason for this is the following, see \cite{MO}.
The nested Hilbert scheme $Z_{n,m}$ can be embedded in the moduli space of ideal sheaves
on the three-fold $X=\mathbb A^2\times\mathbb P^1$ with Chern classes
\begin{equation}
c_1=0,\qquad c_2=m[\mathbb P^1],\qquad \chi=m+n.
\end{equation}
This map embeds indeed $Z_{n,m}$ as one of the components of some $\mathbb C^*$-action
on this moduli space. Since the later carries a perfect obstruction theory, by a result of Graber and 
Pandharipande, the $\mathbb C^*$-fixed part of this obstruction defines a perfect obstruction theory on
$Z_{n,m}$. The virtual fundamental class above is exactly the one which comes from this
perfect obstruction theory. We are grateful to Maulik and Okounkov for this remark.
\end{rem}

There is a natural embedding
$\mathbf{E}_R \subset \mathbf{E}_{\Kb}$.
Recall also the subalgebra $\H_{\Kb} \subset \mathbf{E}_{\Kb}$
introduced in Section 3.4.
As we will show in this Section,
$$\Lambda(\Vcb_k)\in\H_{\Kb},\quad\forall k>0.$$
More precisely, let us consider the generating series
$$\Lambda^+(\Vcb)(z)=1+\sum_{k \geqslant 1}
\Lambda(\Vcb_k)z^k \in \mathbf{E}_R[[z]],$$
$$\Lambda^-(\Vcb)(z)=1+\sum_{k \geqslant 1}
\Lambda(qt\Vcb^*_{-k})z^{-k} \in \mathbf{E}_R[[z^{-1}]].$$
Let $\Omega: \U_c \to \H_{\Kb}$ be the isomorphism in 
Theorem~\ref{T:Main} and let the elements $\mathbf{u}_{r,d}$ be 
the standard generators of  $\U_c$.
Set
$$a_{l,0}=\begin{cases} t^{-l/2} \Omega(\textbf{u}_{l,l}) \quad & \text{if} \;l>0,\\ q^{l/2}\Omega(\mathbf{u}_{l,0}) \quad & \text{if}\;l<0.\end{cases}$$

\vspace{.1in}

\begin{theo}\label{T:Virt} We have
\begin{equation}\label{E:Virt1}
\Lambda^+(\Vcb)(z)=
\exp\bigg( -\sum_{n \geqslant 1}(-1)^n (1-t^nq^n)a_{n,0}\frac{z^n}{n}\bigg),
\end{equation}
\begin{equation}\label{E:Virt2}
\Lambda^-(\Vcb)(z)=
\exp\bigg(-\sum_{n \geqslant 1} (1-t^nq^n)a_{-n,0}\frac{z^n}{n}\bigg).
\end{equation}
\end{theo}
\begin{proof} We will deal with (\ref{E:Virt1}) only.
The proof of (\ref{E:Virt2}) goes along similar lines.
In degree one, (\ref{E:Virt1}) reads
$$\Lambda(\Vcb_{n+1,n})=(1-tq) [\mathcal{O}_{Z_{n+1,n}}]=
(1-tq)\Lambda(N^*_{Z_{n+1,n}})$$
which holds by virtue of Lemma~\ref{L:51} (c).
We will prove (\ref{E:Virt1}) by showing that both sides are solution
to a certain recurrence equation. Recall the elements
of $\H_{\Kb}$ given by
$$\mathbf{f}_{0,1}=\prod_n \boldsymbol{\tau}_{n,n}, \qquad
\mathbf{f}_{1,1}=\prod_n \boldsymbol{\tau}_{n+1,n}.$$
For a series $A(z) \in 1 + \mathbf{E}_{\Kb}[[z]]$
we consider the functional equation
\begin{equation}\label{FE}
[\mathbf{f}_{0,1},A(z)]=z \big(A(z)  \mathbf{f}_{1,1}
-tq \mathbf{f}_{1,1} A(z)  \big)
\end{equation}
A solution $A(z)=1+\sum_i x_i z^i$ of (\ref{FE}) is uniquely determined
by its first Fourier coefficient $x_1$. Thus (\ref{E:Virt1}) will be a
consequence of the following two lemmas.

\vspace{.1in}

\begin{lem}\label{L:vert12} The series $\Lambda^+(\Vcb)(z)$ satisfies
(\ref{FE}).\end{lem}

\begin{proof} 
We must prove that for any $n \geqslant 1$ we have
$$[\mathbf{f}_{0,1},\Lambda E_n]=\Lambda E_{n-1} \mathbf{f}_{1,1}-tq \mathbf{f}_{1,1}\Lambda E_{n-1}.$$
This means that for every pair of partitions $\lambda,\mu$ with $\mu \subset \lambda$ and $|\lambda|-|\mu|=n$ it holds
\begin{equation}\label{E:virt3}
[\mathbf{f}_{0,1},\Lambda E_n]|_{I_{\lambda,\mu}}=(\Lambda E_{n-1} \mathbf{f}_{1,1})|_{I_{\lambda,\mu}}-tq (\mathbf{f}_{1,1}\Lambda E_{n-1})|_{I_{\lambda,\mu}}.
\end{equation}
By means of Lemma~\ref{L:51} (a) and Lemma~\ref{L:hytrfdx} we have
\begin{equation*}
[\mathbf{f}_{0,1},\Lambda E_n]|_{I_{\lambda,\mu}}=
\Big(\sum_{s \in \lambda \backslash \mu} q^{x(s)}t^{y(s)}\Big)
\Lambda E_{\lambda,\mu}
\end{equation*}
while
\begin{equation}\label{E:virt5}
(\Lambda E_{n-1} \mathbf{f}_{1,1})|_{I_{\lambda,\mu}}=\sum_s q^{x(s)}t^{y(s)} \cdot \Lambda E_{\lambda,\mu+s} \cdot \Lambda N^*_{\mu+s,\mu}
\cdot \Lambda^{-1}_{\mu+s}
\end{equation}
where the sum ranges over all addable boxes of $\mu$ which lie in $\lambda$, and
\begin{equation}\label{E:virt6}
(\mathbf{f}_{1,1}\Lambda E_{n-1})|_{I_{\lambda,\mu}}=\sum_s q^{x(s)}t^{y(s)} \cdot  \Lambda N^*_{\lambda, \lambda \backslash s}
\cdot\Lambda E_{\lambda \backslash s,\mu} \cdot\Lambda^{-1}_{\lambda \backslash s}
\end{equation}
where the sum now ranges over all removable boxes of $\lambda$ which do not lie in $\mu$.
Quantities \eqref{E:virt5} and \eqref{E:virt6} may be nicely expressed 
via the change of variables introduced by Garsia and Tesler, as in Appendix~A. 
Namely, let $x_1, \ldots, x_r, u_0, \ldots, u_r$ be the variables associated to
$\lambda$ and likewise let $x'_1, \ldots, x'_p, u'_0, \ldots, u'_p$ be those 
associated to $\mu$. Then we have
$$\frac{(\Lambda E_{n-1} \mathbf{f}_{1,1})|_{I_{\lambda,\mu}}}{\Lambda_{\lambda,\mu}}=\frac{1}{(1-q)(1-t)}\sum_{i=1}^r \bigg\{\prod_{j=0}^r (u_j-x_i)\cdot \prod_{\substack{j=1 \\ j \neq i}}^r (x_j-x_i)^{-1} \cdot \prod_{j=1}^p (x'_j-x_i) \cdot \prod_{j=0}^p (u'_j-x_i)^{-1}\bigg\}$$
and
$$\frac{(\mathbf{f}_{1,1}\Lambda E_{n-1})|_{I_{\lambda,\mu}}}{\Lambda_{\lambda,\mu}}=\frac{qt}{(1-q)(1-t)}\sum_{i=0}^p \bigg\{\prod_{j=0}^r (u'_i-u_j)\cdot \prod_{j=1}^r (u'_i-x_j)^{-1} \cdot \prod_{j=1}^p (u'_i-x'_j) \cdot \prod_{\substack{j=0\\j \neq i}}^p (u'_i-u'_j)^{-1}\bigg\}.$$
Note that by Lemma~\ref{L:GT} (a),
$$(1-t^{-1})(1-q^{-1}) \sum_{s \in \lambda \backslash \mu} q^{x(s)}t^{y(s)}=\bigg( \sum_{i=1}^r x_i - \sum_{i=0}^r u_i \bigg) - \bigg( \sum_{i=1}^p x'_i -\sum_{i=0}^p u'_i\bigg)$$
so that (\ref{E:virt3}) reduces to the following combinatorial identity,
proved in Appendix~A.2.

\vspace{.1in}

\noindent
\begin{claim}\label{claim:5.4}
The following identity holds 
(in the field of rational functions in the variables involved)
\begin{equation}\label{E:vert8}
\begin{split}
\bigg( \sum_{i=1}^r& x_i - \sum_{i=0}^r u_i \bigg) - \bigg( \sum_{i=1}^p x'_i -\sum_{i=0}^p u'_i\bigg)=\\
&=\sum_{i=1}^r \bigg(\prod_{j=0}^r (u_j-x_i)\cdot \prod_{\substack{j=1 \\ 
j \neq i}}^r (x_j-x_i)^{-1} \cdot \prod_{j=1}^p (x'_j-x_i) \cdot \prod_{j=0}^p i
(u'_j-x_i)^{-1}\bigg)\\
& \qquad - \sum_{i=0}^p \bigg(\prod_{j=0}^r (u'_i-u_j)\cdot 
\prod_{j=1}^r (u'_i-x_j)^{-1} \cdot \prod_{j=1}^p (u'_i-x'_j) \cdot 
\prod_{\substack{j=0\\j \neq i}}^p (u'_i-u'_j)^{-1}\bigg).
\end{split}
\end{equation}
\end{claim}
\end{proof}

\vspace{.1in}

Next, we abbreviate
$$A(z)=\exp \Big( - \sum_{n \geqslant 1} (-1)^n(1-t^nq^n) a_{n,0}z^n/n \Big).$$

\begin{lem}\label{L:vert11} The series
$A(z)$
satisfies (\ref{FE}).\end{lem}

\begin{proof} This is again a direct computation.
Set $B(z)=\sum_n (-1)^n(t^nq^n-1)a_{n,0}z^n/n$ so that $A(z)=\exp(B(z))$.
Introduce the elements $a_{l,1}=t^{-l/2}\Omega(\mathbf{u}_{l,l+1}).$
We have
$$[\mathbf{f}_{0,1}, A(z)]=\sum_{k=1}^{\infty} \frac{1}{k!}
[\mathbf{f}_{0,1}, B(z)^k]=\sum_{k=1}^{\infty} \frac{1}{k!}
\sum_{j=0}^{k-1} B(z)^j [\mathbf{f}_{0,1}, B(z)] B(z)^{k-1-j}.$$
Using the relation $[\mathbf{f}_{0,1}, a_{n,0}]=a_{n,1}$ we get
$$[\mathbf{f}_{0,1},B(z)]=
\sum_{n \geqslant 1}(-1)^n (t^nq^n-1)a_{n,1} \frac{z^n}{n}$$
and more generally if we set
$$B^{(s)}(z)=\sum_{l_1, \ldots, l_s \geqslant 1} (-1)^{l_1+ \cdots + l_s}
(t^{l_1}q^{l_1}-1) \cdots (t^{l_s}q^{l_s}-1) a_{l_1+ \cdots + l_s,1}
\frac{z^{l_1+ \cdots + l_s}}{l_1 \cdots l_s}$$
then
$$[B^{(s)}(z),B(z)]=B^{(s+1)}(z).$$
We deduce that
\begin{equation}\label{E:vert9}
\begin{split}
[\mathbf{f}_{0,1},A(z)]=&\sum_{k=1}^{\infty} \frac{1}{k!}
\bigg( \sum_{j=0}^{k-1} \sum_{l=0}^{k-1-j} \begin{pmatrix} k-1-j\\
l \end{pmatrix} B^{l+j}(z) B^{(k-j-l)}(z) \bigg)\\
=&\sum_{k=1}^{\infty} \frac{1}{k!} \sum_{u=0}^{k-1}
\begin{pmatrix} k\\ u \end{pmatrix} B^{u}(z) B^{(k-u)}(z).
\end{split}
\end{equation}
In a similar fashion, starting from the relation
$[\mathbf{f}_{1,1}, a_{n,0}]=a_{n+1,1}$ we get
\begin{equation}\label{E:vert10}
z[\mathbf{f}_{1,1},A(z)]=\sum_{k=1}^{\infty} \frac{1}{k!} \sum_{u=0}^{k-1}
\begin{pmatrix} k\\ u \end{pmatrix} B^{u}(z) \widetilde{B}^{(k-u)}(z).
\end{equation}
where
$$\widetilde{B}^{(s)}(z)=\sum_{l_1, \ldots, l_s \geqslant 1}
(-1)^{l_1+ \cdots + l_s}(t^{l_1}q^{l_1}-1) \cdots
(t^{l_s}q^{l_s}-1) a_{l_1+ \cdots + l_s+1,1}
\frac{z^{l_1+ \cdots + l_s+1}}{l_1 \cdots l_s}.$$
From (\ref{E:vert9}) and (\ref{E:vert10}) we see that (\ref{FE})
reduces to the identity
$$\sum_{s} \frac{1}{s!} \sum_{l_1 + \cdots + l_s=n}
\frac{(t^{l_1}q^{l_1}-1) \cdots (t^{l_s}q^{l_s}-1)}{l_1 \cdots l_s}=
tq \sum_{s} \frac{1}{s!} \sum_{l_1 + \cdots + l_s=n-1}
\frac{(t^{l_1}q^{l_1}-1) \cdots (t^{l_s}q^{l_s}-1)}{l_1 \cdots l_s}$$
for all $n$. This identity is in turn a corollary of the following formula,
whose proof is left to the reader
$$\sum_{s} \frac{1}{s!} \sum_{l_1 + \cdots + l_s=n}
\frac{(t^{l_1}q^{l_1}-1) \cdots (t^{l_s}q^{l_s}-1)}{l_1 \cdots l_s}=
t^{n-1}q^{n-1}(tq-1).$$
Lemma~\ref{L:vert11} and Theorem~\ref{T:Virt} are proved.
\end{proof}
\end{proof}

\vspace{.2in}

\paragraph{\textbf{5.3.}} Once Theorem~\ref{T:Virt} is established,
it is a simple task to describe the action of the virtual classes
$\Lambda\Vcb_n$ and $\Lambda\Vcb^*_{-n}$ on $\mathbf L_{\mathbf K}$.
The following corollary of Theorem~\ref{T:Virt} may be viewed as a
K-theoretic analog of Nakajima's formulas in Borel-Moore homology (see (\ref{E:Nak1}), (\ref{E:Nak2})).

\vspace{.1in}

\begin{cor}\label{C:Virt} As operators in $\mathbf L_{\mathbf K}\simeq \boldsymbol{\Lambda}_\Kb$,
 we have
$$1+ \sum_{n \geqslant 1} \boldsymbol{\tau}^*_n \otimes \Lambda(\Vcb_n)z^n=
\exp\bigg( -\sum_{n \geqslant 1} (-1)^n \frac{1-t^nq^n}{1-q^n} p_n
\frac{z^n}{n}\bigg),$$
$$1+ \sum_{n \geqslant 1} \Lambda(qt\Vcb^*_{-n})z^n=
\exp\bigg( -\sum_{n \geqslant 1} \frac{1-t^nq^n}{1-t^n}
\frac{\partial}{\partial p_n} \frac{z^{-n}}{n}\bigg).$$ Here  we
have set $$\boldsymbol{\tau}^*_n \otimes \Lambda(\Vcb_n)=\prod_k
\boldsymbol{\tau}^*_{n+k,k} \otimes \Lambda(\Vcb_{n+k,k}).$$
\end{cor}

As mentioned to us by A. Oblomkov, a very similar result appears in \cite{MO}.

\vspace{.2in}

\section{Hecke operators}

\vspace{.1in}

This section is devoted to the action of the Hecke operators on $\U^{>}$.
The definition of the action of $\U^0$ on $\U^>$ is given in Section 6.1.
In Proposition 6.2 this action is expressed in terms of the tensor product by
the tautological bundle over $\Hin$.
In Proposition 6.3 we prove that this action equips $\U^>$ with the structure
of a torsion-free module over the $\Kb$-algebra of symmetric polynomials.
This technical result will be used later, in Theorem 7.14, to compare the
elliptic Hall algebra with the equivariant $K$-theory of the commuting variety.

\vspace{.2in}

\paragraph{\textbf{6.1.}} Recall that $\U^0$ is a commutative polynomial
algebra generated by elements $\{\bt_{0,l}\;;\; l \in \Z^*\}$.
Let $J \subset \U^0$ be the augmentation ideal.
We denote by $\sigma$ the automorphism of $\U^0$ determined by
\begin{equation}\sigma(\bt_{0,l})=\sign(l) \bt_{0,l}.\end{equation} 
We define the map
$$\bullet~: \U^0 \otimes \U^> \to \U^>,\quad
u \bullet v= \omega(\sigma(u)v),$$
where $\omega : \U^{\geqslant } \to \U^>$ is the projection
along $ \U^>\otimes J$. It is clear that
$$(u u') \bullet v=u \bullet (u' \bullet v),$$ so that
$\bullet$ equips $\U^>$ with the structure of a $\U^0$-module.
This action restricts to the graded pieces
$$\U^>[r]=\bigoplus_{d \in \Z} \U^>[r,d],\quad r \geqslant 1.$$
From \cite{SV}, Theorem~6.3, we have
\begin{equation}
\bt_{0,l} \bullet v=\sign(l)\,[\bt_{0,l}, v],\qquad
l \in \Z,\qquad v \in \U^>.
\end{equation}
For any $r \geqslant 1$ we define also a projection
\begin{equation*}
\pi_r~: \U^0 \to \Kb[z_1^{\pm 1}, \ldots, z_r^{\pm 1}]^{\mathfrak{S}_r},\quad
\bt_{0,l} \mapsto p_l(z_1, \ldots, z_r)
\end{equation*}
where $p_l(z_1, \ldots, z_r)=\sum_i z_i^l$ is the power sum function.

\vspace{.1in}

\begin{rem}
When $\U^>$ is interpreted as the Hall algebra
of an elliptic curve  $\bullet$ is identified with the action of
Hecke operators, see \cite{SV}, Section~6. Indeed $\U^0\cap \U^+$ is the Hall
algebra of the category of torsion sheaves, $\U^>$ is the Hall algebra
of vector bundles, and the map $\omega: \U^{\geqslant} \to \U^>$ is the
restriction map from the Hall algebra of all coherent sheaves
to that of vector bundles.
\end{rem}

\vspace{.2in}

\paragraph{\textbf{6.2.}} Fix an integer $r \geqslant 1$. 
The nested Hilbert scheme $Z_{r+k,k} \subset \text{Hilb}_{r+k} \times \text{Hilb}_{k}$
carries a tautological bundle $\boldsymbol{\tau}_{r+k,k}$ whose fiber over a point $(I,J)$ is equal to $J/I$. The 
character of $\boldsymbol{\tau}_{r+k,k}$ over the $T$-fixed point $I_{\lambda,\mu}$ is given by the following 
expression
\begin{equation}\label{E:chartaurk}
[\boldsymbol{\tau}_{r+k,k}{|}_{I_{\lambda,\mu}}]=\sum_{s \in \lambda \backslash \mu} t^{y(s)}q^{x(s)}.
\end{equation}
Let us put
$$\boldsymbol{\tau}_{r}=\prod_{k} \boldsymbol{\tau}_{r+k,k}
\in \prod_k K^T(\text{Hilb}_{r+k} \times \text{Hilb}_{r}).$$
We introduce one final piece of notation~:
for any $l \in \Z^*$ let $$p_l(\boldsymbol{\tau}_r)=\Psi_l(\boldsymbol{\tau}_r)
=\prod_k \Psi_l(\boldsymbol{\tau}_{r+k,k})$$ denote the $l$th Adams operation.
If $\lambda=(\lambda_1, \ldots, \lambda_s)$ is a partition then we set
$$p_{\lambda}(\boldsymbol{\tau}_{r})=
p_{\lambda_1}(\boldsymbol{\tau}_{r}) \otimes \cdots
\otimes p_{\lambda_s}(\boldsymbol{\tau}_{r}),$$ the tensor product of sheaves.
We extend this notation to any symmetric function
$\theta=\sum a_{\lambda}p_{\lambda}$ by linearity.
Note that because $\boldsymbol{\tau}_{r}$ is of rank $r$ we have
$$e_{l}(\boldsymbol{\tau}_r)=\Lambda^l \boldsymbol{\tau}_r=0,\qquad l >r.$$
Therefore $\theta(\boldsymbol{\tau}_r)$ makes sense for
$\theta \in \Kb[z_1^{\pm 1}, \ldots, z_r^{\pm 1}]^{\mathfrak{S}_r}$.
In this notation, we have $$1(\tau_r)=\prod_k \mathcal{O}_{Z_{r+k,k}}.$$
The following proposition connects the action of the Hecke
operators with the tautological bundles.
Recall the algebra isomorphism $\Omega~: \U \to \H_{\Kb}$ given in Section 4.6.

\begin{prop}\label{P:Hecke1}
For any $v \in \U^{>}[r]$ and any $u \in \U^0$ we have
\begin{equation}\label{E:Hecketaut}
\Omega(u \bullet v)=\pi_r(u)(\boldsymbol{\tau}_r) \otimes \Omega(v).
\end{equation}
\end{prop}
\begin{proof} Fix $v \in  \U^{>}[r]$. It is enough to prove (\ref{E:Hecketaut}) for a system of generators of $\U^0$, such as
$\{\bt_{0,l}\;;\; l \in \Z^*\}$. We have, by definition
\begin{equation*}
\begin{split}
\Omega(\bt_{0,l} \bullet v)&= \sign(l) \, \Omega([\bt_{0,l},v])\\
&= \sign(l) \, [\Omega(\bt_{0,l}),\Omega(v)]\\
&=[\f_{0,l},\Omega(v)].
\end{split}
\end{equation*}
In the convolution diagram
$$\xymatrix{ \text{Hilb}_{r+k} & \text{Hilb}_{r+k} \times \text{Hilb}_k \ar[l]_-{q_1} \ar[r]^-{q_2} & \text{Hilb}_{k}}$$
there is, over $Z_{r+k,k}$ a short exact sequence
\begin{equation}\label{E:ses1}
\xymatrix{0 \ar[r]& \boldsymbol{\tau}_{r+k,k} \ar[r] & q_1^* \boldsymbol{\tau}_{r+k} \ar[r]& q_2^*\boldsymbol{\tau}_{k} \ar[r]& 0}
\end{equation}
and hence also a short exact sequence
\begin{equation}\label{E:ses2}
\xymatrix{0 \ar[r]& \Psi_l(\boldsymbol{\tau}_{r+k,k}) \ar[r] & q_1^* \Psi_l(\boldsymbol{\tau}_{r+k}) \ar[r]& q_2^*\Psi_l(\boldsymbol{\tau}_{k}) \ar[r]& 0}
\end{equation}
Equation (\ref{E:Hecketaut}) for $u=\bt_{0,l}$ follows from (\ref{E:ses2}).
\end{proof}

\vspace{.2in}

\paragraph{\textbf{6.3.}} Put $J_r=\mathrm{Ker}(\pi_r)\subset \U^0$.
By Proposition~\ref{P:Hecke1} above and because $\Omega~: \U \to \H_{\Kb}$
is an isomorphism, the Hecke action of $\U^0$ on $\U^{>}[r]$ factors through
$\pi_r$ and yields an action
\begin{equation}\label{E:Hecke2}
\bullet~: \Kb[z_1^{\pm 1}, \ldots, z_r^{\pm 1}]^{\mathfrak{S}_r} \otimes
\U^>[r] \to \U^>[r].
\end{equation}
Observe in particular that the element $z_1 \cdots z_r$ acts on $\U^>[r]$ by an automorphism which we denote by $\xi_r \in \text{Aut}( \U^>[r])$.

\vspace{.1in}

\begin{prop}\label{P:Hecke2} Under this  action $\U^>[r]$ is a torsion-free 
$\Kb[z_1^{\pm 1}, \ldots, z_r^{\pm 1}]^{\mathfrak{S}_r}$-module.
\end{prop}

\begin{proof} We will use some results from \cite{BS}, to which we refer for
definitions. The algebra $\U^+ \subset \U^{\geqslant}$ is equipped with a
(topological) comultiplication
$$\Delta~: \U^+ \to \U^+ \widehat{\otimes} \U^+$$
which satisfies
\begin{equation}\label{E:coproduno}
\Delta(\bt_{0,l})=\bt_{0,l} \otimes 1 + 1 \otimes \bt_{0,l},
\end{equation}
\begin{equation}\label{E:coprodduo}
\Delta(\bt_{1,l})=\bt_{1,l} \otimes 1 + \sum_{d \geqslant 0}  \theta_{0,d}\otimes \bt_{1,l-d}.
\end{equation}
Here $\widehat\otimes$ is the completed tensor product defined in \cite{BS}.
When we need to specify graded components, we write
$\Delta_{r,r'}: \U^+[r+r'] \to \U^+[r] \widehat{\otimes} \U^+[r']$.
We begin with a couple of easy lemmas.

\vspace{.1in}

\begin{lem}\label{L:coprodinj} For any $r \geqslant 1$, the iterated coproduct map $\Delta_{1,\ldots, 1}: \U^+[r] \to \U^+[1] \widehat{\otimes} \cdots \widehat{\otimes} \U^+[1]$ is injective. Moreover, the map $\omega^{\otimes r} \circ \Delta_{1, \ldots, 1}: \U^>[r] \to \U^>[1] \otimes \cdots \otimes \U^>[1]$ is injective.
\end{lem}
\begin{proof} Let $(\;,\;)$ be Green's scalar product on $\U^+$. It is nondegenerate and satisfies the Hopf property, i.e. $(xy,z)=(x \otimes y, \Delta(z))$ for any $x,y,z \in \U^+_{\Kb}$. The first statement of the Lemma follows from the fact that $\bigoplus_{r \geqslant 1}\U^+[r]$ is generated by elements of $\U^+[1]$. In a similar vein, the second statement follows from the fact that $\bigoplus_{r \geqslant 1}\U^>[r]$ is generated by elements of $\U^>[1]$.
\end{proof}

\vspace{.1in}

\begin{lem}\label{L:Heckecoprod} For 
$u \in \Kb[\bt_{0,1}, \bt_{0,2}, \ldots]$ and $v \in \U^>$ we have
$$\omega \otimes \omega \big( \Delta( u \bullet v) \big)= 
\Delta(u) \bullet \Big( \omega \otimes \omega \big( \Delta(v)\big)\Big).$$
\end{lem}

\begin{proof} Using the formula 
$(u_1u_2)\bullet v= u_1 \bullet (u_2 \bullet v)$ we see that it is enough 
to prove the Lemma for $u=\bt_{0,l}$. In that case, we have
\begin{equation*}
\begin{split}
\omega \otimes \omega \big(\Delta(\bt_{0,l} \bullet v)\big)=& \omega \otimes \omega \big( \Delta( [\bt_{0,l},v])\big)\\
=&\omega \otimes \omega \big( [\bt_{0,l} \otimes 1 + 1 \otimes \bt_{0,l},\Delta(v)]\big)\\
=& (\bt_{0,l} \otimes 1 + 1 \otimes \bt_{0,l} )\bullet 
\Big(\omega \otimes \omega (\Delta(v))\Big)\\
=& \Delta(\bt_{0,l}) \bullet 
\Big(\omega\otimes\omega \big( \Delta(v)\big)\Big).
\end{split}
\end{equation*}
\end{proof}

\vspace{.1in}

We may now proceed with the proof of Proposition~\ref{P:Hecke2}.
Let $u \in \U^0$, $v \in \U^>[r]$, $v \neq 0$,
and suppose that $u \bullet v=0$. Our aim is to show that $u \in J_r$.
Twisting by the automorphism $\xi_r$ if necessary, we may assume that $u \in \Kb[\bt_{0,1}, \bt_{0,2}, \ldots]$.
Let us consider the subspace
$$N_r=\sum_{s=1}^r (\U^0) ^{\otimes (s-1)} \otimes J_1
\otimes (\U^0)^{\otimes (r-s)}\subset(\U^0)^{\otimes r}.$$
We claim that it suffices to prove that $\Delta^{(r)}(u) \in N_r$.
Indeed, using the bialgebra isomorphism
$$\pi~: \Kb[\bt_{0,1}, \bt_{0,2}, \ldots] \to 
\Kb[z_1, z_2, \ldots]^{\mathfrak{S}_{\infty}},
\qquad \bt_{0,l} \mapsto p_l$$
we have $J_r \cap \Kb[\bt_{0,1}, \bt_{0,2}, \ldots]=
\pi^{-1} \Kb[e_{r+1}, e_{r+2}, \ldots]$
and standard symmetric functions arguments show that
$(\Delta^{(r)})^{-1}(\pi(N_r)) =\Kb[e_{r+1}, e_{r+2}, \ldots]$.

Now, by Lemma~\ref{L:Heckecoprod} we have
$$0=\Delta_{1, \ldots, 1}(u \bullet v)=\Delta^{(r)} \bullet \omega^{\otimes r} \big(\Delta_{1,\ldots, 1}(v)\big).$$
By Lemma~\ref{L:coprodinj}, we have $\omega^{\otimes r}\big(\Delta_{1, \ldots, 1}(v)\big) \neq 0$, and we may write
\begin{equation}\label{E:Hecke3}
\omega^{\otimes r}\big( \Delta_{1, \ldots, 1}(v)\big)=
\sum_{\underline{d}} a_{\underline{d}}\, \bt_{1,d_1} \otimes \cdots \otimes
\bt_{1,d_r}, \qquad a_{\underline{d}}\neq 0.
\end{equation}
The above sum is in general infinite. However, it follows from \cite{BS},
Proposition 2.1, that~: $d_r$ is bounded from above,
for $d_r$ fixed $d_{r-1}$ is bounded above,
and more generally, for fixed $d_r, d_{r-1}, \ldots, d_{r-i+1}$
the possible values of $d_{r-i}$ are bounded above.
In particular, there exists a unique element
$\underline{d}=(d_1, \ldots, d_r)$ appearing in (\ref{E:Hecke3}) which is maximal for the (right) lexicographic order.
Similarly, modulo $N_r$ we may write
\begin{equation}\label{E:Hecke4}
\Delta^{(r)}(u)=\sum_{\underline{l}} b_{\underline{l}}\, \bt_{0,1}^{l_1}
\otimes \cdots \otimes \bt_{0,1}^{l_r}, \qquad b_{\underline{l}}\neq 0.
\end{equation}
This sum is finite and may be empty. If $\Delta^{(r)}(u) \not\in N_r$
then there exists a unique $\underline{l}=(l_1, \ldots, l_r)$ which
is maximal for the (right) lexicographic order. But then, considering the $\mathbf{Z}$-graded component
\begin{equation*}
\begin{split}
\Delta_{(1,l_1+d_1), \ldots, (1, l_r+d_r)}(u \bullet v)&=
a_{\underline{d}}b_{\underline{l}} \big(\bt_{1,l_1} \otimes \cdots \otimes \bt_{1,l_r}\big) \bullet \big(
\bt_{1,d_1} \otimes \cdots \otimes \bt_{1,d_r}\big)\\
&=a_{\underline{d}}b_{\underline{l}} \big(
\bt_{1,d_1+l_1} \otimes \cdots \otimes \bt_{1,d_r+l_r}\big)\neq 0
\end{split}
\end{equation*}
we reach a contradiction. Thus $\Delta^{(r)}(u) \in N_r$ and $u \in J_r$. We are done.
\end{proof}

\vspace{.2in}

\paragraph{\textbf{6.4.}} By ``transport de structure" and using
Theorem~\ref{T:Main}, we get a Hecke action
$$\bullet~: \H^0_{\Kb} \otimes \H^>_{\Kb} \to \H^>_{\Kb}$$
which is determined by
\begin{equation}\label{E:defHeckeaction}
\begin{split}
(u_1u_2) \bullet v= u_1 \bullet ( u_2 \bullet v),\quad
\f_{0,l} \bullet v=[\f_{0,l},v], \quad \forall\; l \in \Z.
\end{split}
\end{equation}
This action is torsion free in the same sense as Proposition~\ref{P:Hecke2}. Note that we may also use formulas (\ref{E:defHeckeaction}) to extend this to a Hecke action
\begin{equation}\label{E:defHeckeaction2}
\bullet ~:\H^0_{\Kb} \otimes \text{End}_\Kb(\mathbf{L}_{\Kb}) \to
\text{End}_\Kb(\mathbf{L}_{\Kb}).
\end{equation}
Using this notation we have
\begin{equation}\label{E:compatHecke}
h \bullet \psi(u)=\psi(h \bullet u)
\end{equation}
for any $h \in \H^0_{\Kb}$ and $u \in \H_{\Kb}^>$.

\vspace{.2in}

\section{K-theory of the commuting variety}

\vspace{.2in}

The aim of this section is to introduce a ring structure on the
equivariant Grothendieck group of the commuting variety and to
compare this ring with the positive part of the elliptic Hall
algebra. The ring $\bar\Cb_R$ is given in Proposition 7.5. Then we
prove in Proposition 7.9 that it acts on the $R$-module $\mathbf
L_R$ equal to the $K$-theory of the Hilbert scheme. Next, in
Proposition 7.10 we compare the Hecke action on the algebra $\mathbf
H_{\mathbf K}^>$, which is the positive part of the elliptic Hall
algebra by Theorem~\ref{T:Main}, with the natural action of the
representation ring on the equivariant $K$-theory of the commuting
variety. Finally, in Theorem 7.14, we compare $\bar\Cb_{\mathbf
K}$ with $\mathbf H_{\mathbf K}^>$.

\vspace{.2in}

\paragraph{\textbf{7.1.}}
First let us recall a few general facts. Let $G$ be a complex linear
algebraic group. By a variety we'll always mean a quasi-projective
complex variety. We call $G$-variety a variety with a rational
action of $G$.

Let $P\subset G$ a parabolic subgroup and $H\subset P$ a Levi
subgroup. Fix a $H$-variety $Y$.  The group $P$ acts on $Y$ through
the obvious group homomorphism $P\to H$. Let $X=G\times_PY$ be the
induced $G$-variety.

Now assume that $Y$ is smooth. Given a smooth subscheme $O\subset Y$
let $T^*_OY\subset T^*Y$ be the conormal bundle to $O$. It is
well-known that the induced $H$-action on $T^*Y$ is Hamiltonian and
that the zero set of the moment map is the closed $H$-subvariety
$$T^*_HY=\bigcup_{O}T^*_OY\subset T^*Y,$$
where $O$ runs over the set of $G$-orbits. The following lemma is
left to the reader.

\begin{lem}
We have $T^*X=T^*_P(G\times Y)/P$ and $T^*_GX=G\times_PT^*_HY$. The
induction yields a canonical isomorphism $K^H(T^*_HY)=K^G(T^*_GX)$.
\end{lem}

We'll call fibration a smooth morphism which is locally trivial in
the Zariski topology. Let $X'$ be a smooth $G$-variety and $V$ be a
smooth $H$-variety. Assume that we are given $H$-equivariant
homomorphisms $p:V\to Y$ and $q:V\to X'$ which are a fibration and a
closed embedding respectively. Set
$W=G\times_PV$ and consider the following maps
$$\aligned
&g:\ W\to X',\quad (g,v)\, \mod\, P\mapsto gq(v),\hfill\cr &f:\ W\to
X,\quad (g,v)\, \mod\, P\mapsto (g,p(v))\, \mod\, P.
\hfill\endaligned$$ The following properties are immediate.

\begin{lem}
The map $f$ is a $G$-equivariant fibration, the map $g$ is a
$G$-equivariant proper morphism, and the map $(f,g)$ is a closed
embedding $W\subset X\times X'$. The varieties $V$, $W$, $X$, $X'$
are smooth.
\end{lem}

We'll identify $W$ with its image in $X\times X'$. Let
$Z=T^*_W(X\times X')$ the conormal bundle. It is again a smooth
$G$-variety. The obvious projections yield $G$-equivariant maps
$$\phi:Z\to T^*X,\quad\psi:Z\to T^*X'.$$ Consider the $G$-variety
$$Z_G=Z\cap(T^*_GX\times T^*_GX').$$
Recall that a morphism of varieties $S\to T$ is called regular if it
is the composition of a regular immersion $S\subset S'$, i.e., an
immersion which is locally defined by a regular sequence, and of a
smooth map $S'\to T$. Note that a regular map has a finite
tor-dimension and that a morphism $S\to T$ is regular whenever $S$
and $T$ are smooth.

\begin{lem}
(a) The map $\psi$ is proper and regular, the map $\phi$ is regular.

(b) We have $\phi^{-1}(T^*_GX)=Z_G$
and $\psi(Z_G)\subset T^*_GX'$.
\end{lem}

\begin{proof}
To prove that $\psi$ is a proper morphism we may assume that the
fibration $f$ is indeed the obvious projection $$f:\ W=X\times U\to
X,$$ where $U$ is a smooth $G$-variety. Since the map $g$ is proper,
it is enough to check that for all $w=(x,u)\in W$ we have
$$\{\xi\in T^*_w(X\times X');\, \xi(T_wW)=\xi(T_{x'}X')=0\}=\{0\}.$$
This is obvious because we have
$$T_wW+T_{u}X'=T_xX+T_{u}U+T_{u}X'=T_xX+T_{u}X'.$$ The maps $\phi$, $\psi$
are regular because $Z$ is smooth. For instance $\phi$ is the
composition of the projection $T^*(X\times X')\to T^*X$, which is
smooth, and of the obvious inclusion $Z\subset T^*(X\times X')$,
which is regular. Claim $(a)$ is proved.

Now let us concentrate on $(b)$. The second claim is obvious. Let us
prove the first one. Since the set $W\subset X\times X'$ is
preserved by the diagonal action of $G$ we have
$$T_{x'}(Gx')\subset T_wW+T_x(Gx),\quad
 \forall w=(x,x')\in W.$$ Thus we have also
$$Z\cap(T^*_GX\times T^*X')\subset T^*_GX\times T^*_GX'.$$

\end{proof}

\vskip3mm

To avoid confusions we may abbreviate
$$\phi_G=\phi|_{Z_G}:Z_G\to T^*_GX,\quad \psi_G=\psi|_{Z_G}:Z_G\to
T^*_GX'.$$

\vskip3mm

\paragraph{\textbf{7.2.}}
Recall that for any $G$-variety $M$ and any closed $G$-stable
subvariety $N\subset M$ the direct image by the obvious inclusion
$N\to M$ identifies $K^ G(N)$ with the complexified Grothendieck
group $K^G(M\on N)$ of the category of $G$-equivariant coherent
sheaves on $M$ supported on $N$.
Since the map $\psi$ is a proper morphism the derived direct image
yields a map
$$R\psi_*:K^G(Z)\to K^G(T^*X').$$
By Lemma 7.3$(b)$ we have also a map
$$R\psi_*:K^G(Z_G)=K^G(Z\on Z_G)\to K^G(T^*X'\on T^*_GX')=K^G(T^*_GX').$$
Since the map $\phi$ has a finite tor-dimension the derived
pull-back yields a map
$$L\phi^*:K^G(T^*X)\to K^G(Z).$$
By definition $L\phi^*$ is the composition of the pull-back by the
projection $T^*X\times T^*X'\to T^*X$ and the derived pull-back by
the regular immersion $Z\subset T^*X\times T^*X'$. By Lemma 7.3$(b)$
we have also a map
$$L\phi^*:K^G(T^*_GX)=K^G(T^*X\on T^*_GX)\to K^G(Z\on Z_G)=K^G(Z_G).$$
Composing $R\psi_*$ and $L\phi^*$ we get a map
$$R\psi_*\circ L\phi^*:K^G(T^*_GX)\to K^G(T^*_GX').$$
By Lemma 7.1 the induction yields also an isomorphism
$$K^H(T^*_HY)=K^G(T^*_GX).$$ Composing it by $R\psi_*\circ L\phi^*$
we obtain a map
\begin{equation}\label{CO:0}K^H(T^*_HY)\to K^G(T^*_GX').\end{equation}

\vskip3mm

\paragraph{\textbf{7.3.}}
Now, we apply the general construction recalled above to the
particular case of the commuting variety. First, let us fix some
notation. Let $E$ be a finite dimensional $\C$-vector space. We'll
abbreviate $G_E=GL(E)$ and $\gen_E=\End(E)$.
Set
$$C_E=\{(a,b)\in\gen^2; [a,b]=0\}.$$
If no confusion is
possible we'll write $C=C_E$, $G=G_E$ and $\gen=\gen_E$.
We equip $C$ with the diagonal $G$-action and the $T$-action
such that
\begin{equation}\label{E:taction1}
(z_1,z_2)\cdot(a,b)=(z_1a,z_2b).
\end{equation}

\vspace{.2in}

\paragraph{\textbf{7.4.}}
Next we fix a subspace $E_1\subset E$ and we set $E_2=E/E_1$.
We may write $G_i=G_{E_i}$, $\gen_i=\gen_{E_i}$,
$C_i=C_{E_i}$, etc, for $i=1,2$.
Set $$H=G_{1}\times G_{2},\quad P=\{g\in
G;g(E_1)=E_1\}.$$  Let $\gen$, $\pen$ and $\hen$ be the
corresponding Lie algebras. Put
$$Y=\hen,\quad X'=\gen,\quad V=\pen.$$
The $G$-action on $X'$ and the $H$-action on $Y$ are the adjoint
ones.
Put $$C_\gen=C,\quad
C_\hen=C_{1}\times C_{2},\quad C_\pen=\pen^2\cap C.$$
For each $a\in\pen$ let $a_\hen\in\hen$
be the graded linear map associated with $a$.
We apply the
general construction in Section 7.2 with the map $p:V\to Y$,
$a\mapsto a_\hen$ and the obvious inclusion $q:V\to X'$. By the
canonical isomorphisms $\gen^*=\gen$, $\hen^*=\hen$ we'll always
mean the isomorphisms given by the trace.

\begin{lem}
(a) We have
$X=G\times_P\hen$ and $W=G\times_P\pen$.
The maps $W\to X$, $W\to X'$ are given by
$
(g,a)\ \mod \ P\mapsto(g,a_\hen)\ \mod\ P$
and $(g,a)\ \mod\
P\mapsto gag^{-1}.
$

(b) The canonical isomorphisms $\gen^*=\gen$,
$(\gen\times\hen)^*=\gen\times\hen$ yield isomorphisms of
$G$-varieties
$$\begin{aligned}
T^*X'=\gen^2,\quad
T^*X=G\times_P\{(c,a,b)\in\pen\times\hen\times\hen;c_\hen=[a,b]\},\quad
Z=G\times_P\pen^2.\end{aligned}$$
For each $a,b\in\pen$ we have
$$\begin{gathered}
\phi((g,a,b)\ \mod\ P)=(g,[a,b],a_\hen,b_\hen)\ \mod\ P,\quad
\psi((g,a,b)\ \mod\ P)=(gag^{-1},gbg^{-1}). \end{gathered}$$

(c) We have canonical isomorphisms of $G$-varieties
$$\begin{aligned}
T^*_GX=G\times_PC_\hen,\quad T^*_GX'=C,\quad
Z_G=G\times_PC_\pen.
\end{aligned}$$
The maps $\phi_G$, $\psi_G$ are the obvious ones.

\end{lem}

\begin{proof}
We'll frequently use the following isomorphisms without
mentioning them explicitly
$$\begin{gathered}
G\times_P(\hen\times\gen)\to X\times X',\quad
(g,a,b)\ \mod\ P\mapsto ((g,a)\ \mod\ P, gbg^{-1}),\cr
G\times\gen^*\to T^*G,\quad(g,f)\mapsto gf.
\end{gathered}$$
By Lemma 7.1 we have $$\begin{aligned}
T^*(X\times X')
&=T^*_P(G\times\hen\times\gen)/P\\
&=G\times_P\{(f,a)\in(\gen\times\hen\times\gen)^*\times(\hen\times\gen);\,
f(-b,[\delta b,a])=0,\forall b\in \pen\}.
\end{aligned}$$
Let  $\delta$ is the linear map
$$\delta:\pen\to\hen\times\gen,\quad a\mapsto(a_\hen,a).$$
We have
$$Z=T^*_W(X\times X'),\quad
W= G\times_P\delta\pen,\quad T^*W=T^*_P(G\times\pen)/P.$$
Let $\delta\pen^\perp\subset(\hen\times\gen)^*$ be the orthogonal of
$\delta\pen$.
We have also
$$\begin{aligned}
&T^*(X\times
X')|_W=G\times_P\{(f,a)\in(\gen\times\hen\times\gen)^*\times\pen;\,
f(-b,\delta[b,a])=0,\forall b\in \pen\},\\
&T^*W=G\times_P\{(f,a)\in(\gen\times\pen)^*\times\pen;\,
f(-b,[b,a])=0,\forall b\in \pen\},\\
&Z=G\times_P(\delta\pen^\perp\times\pen).
\end{aligned}$$
Here the inclusion $Z\subset T^*(X\times X')$ is given
by the inclusion $\delta :\pen\to\hen\times\gen$ and the inclusion
$$\delta\pen^\perp=\{0\}\times\delta\pen^\perp\subset
\{0\}\times(\hen\times\gen)^*\subset(\gen\times\hen\times\gen)^*.$$
The canonical isomorphism
$(\hen\times\gen)^*\to\hen\times\gen$ identifies $\delta\pen^\perp$
with $\delta'\pen\simeq\pen$, where
$$\delta':\pen\to\hen\times\gen,\ a\mapsto(-a_\hen,a).$$
This yields an isomorphism
$$Z=G\times_P\pen^2.$$ By Lemma 7.1 we have
$$\begin{aligned}
T^*X
&=T^*_P(G\times\hen)/P\\
&=G\times_P\{(f,a)\in(\gen\times\hen)^*\times\hen;\,
f(-b,[b_\hen,a])=0,\forall b\in \pen\},\\
T^*_GX&=G\times_PT^*_H\hen,\\
&=G\times_P C_\hen,\\
&=G\times_P\{(f,a)\in(\gen\times\hen)^*\times\hen;\,
f(-b,[b_\hen,a])=f(c,0)=0,\forall b\in \pen, c\in\gen\}.
\end{aligned}$$
Here the inclusion $\hen^*\times\hen\subset(\gen\times\hen)^*\times\hen$
is given by the map
$$\hen^*=\{0\}\times\hen^*\subset(\gen\times\hen)^*.$$
The map $\phi$ is the composition of the
chain of maps
$$Z\subset T^*(X\times X')\to T^*X.$$
Fix $a,b\in\pen$. Consider the element $\xi=(g,a,b)\
\mod\ P$ of $Z$. We may identify $a$ with $\delta'a$, which
can be regarded as an element in
$$(\hen\times\gen)^*=
\{0\}\times(\hen\times\gen)^*\subset(\gen\times\hen\times\gen)^*=
\gen^*\times\hen^*\times\gen^*,$$
and $b$ with $\delta b$, which is an element of $\hen\times\gen$.
So $\xi$ can be viewed as an element in $T^*(X\times X')$, see above. Set
$\delta'a=
(0,f_\hen,f).$
A short computation yields
$$\phi(\xi)= (g,-[b_\hen,f],f_\hen,b_\hen)\ \mod\ P$$
where the bracket is the coadjoint action. Now, observe that the
canonical map identifies $f_\hen,f$ with $-a_\hen,a$ respectively.
This yields the formula for $\phi$ in part $(b)$. Finally Lemma
7.3$(b)$ yields
$$Z_G=\phi^{-1}(T^*_GX).$$ Therefore we have
$Z_G=G\times_PC_\pen.$ The other claims are left to the reader.
\end{proof}

\vspace{.2in}

Next, we set $\bar\Cb_{E,R}=K^{G\times T}(C_E)$. A vector space
isomorphism $E\simeq E'$ yields a $R$-module isomorphism $\mathbf
C_{E,R}\simeq \bar\Cb_{E',R}.$ Let
$$\bar\Cb_R=\ind_E \bar\Cb_{E,R},$$
where the limit runs over the groupoid formed by all finite
dimensional vector spaces with their isomorphisms.
There is a $G\times T$-action on $T^*X$ and $T^*X'$ is given by
$$\begin{gathered}
(z_1,z_2)\cdot(g,c,a,b)\ \mod\ P=(g,z_1z_2c,z_1a,z_2b)\ \mod\ P,\cr
(z_1,z_2)\cdot(g,a,b)\ \mod\ P=(g,z_1a,z_2b)\ \mod\
P.\end{gathered}$$
We define as in (\ref{CO:0}) a $R$-linear homomorphism
$$K^{H\times T}(C_\hen)\to K^{G\times T}(C)
=\bar\Cb_{E,R}.$$ We'll abbreviate $\bar\Cb_i=\bar\Cb_{E_i,R}$
for $i=1,2$. By the Kunneth formula, see \cite{CG}, Chapter~5.6, it
can be viewed as a map
\begin{equation}\label{CO:3}\bar\Cb_{1}\otimes_R \bar\Cb_{2}\to \mathbf
C_{E,R}.\end{equation}

\vspace{.1in}

\begin{prop}  The map (\ref{CO:3}) equips $\bar\Cb_R$
with the structure of an associative unital $R$-algebra.
\end{prop}

\begin{proof}
Fix a flag $E_1\subset E_{2}\subset E$.
First, we define the following varieties
\begin{itemize}
\item
$X_1$ is the set of tuples $(F_1,F_2,a)$
where $F_1\subset F_2\subset E$ is a
flag such that $F_1\simeq E_1$, $F_2\simeq E_2$
and $a$ is an endomorphism of the graded vector space
$F_1\oplus(F_2/F_1)\oplus (E/F_2)$,
\item
$X_2$ is the set of pairs $(F_1,a)$ where $F_1\subset E$ is a
vector subspace isomorphic to $E_1$ and $a$ is an endomorphism
of the graded vector space $F_1\oplus(E/F_1)$,
\item $X_3=\gen$.
\end{itemize}
Next, we define the following ones
\begin{itemize}
\item
$W_1$ is the set of pairs $(F_1,a)$ where $F_1\subset E$ is a
vector subspace isomorphic to $E_1$ and $a\in\gen$ preserves $F_1$,
\item
$W_2$ is the set of tuples $(F_1,F_2,a)$
where $F_1\subset F_2\subset E$ is a
flag such that $F_1\simeq E_1$, $F_2\simeq E_2$
and $a\in\gen$ preserves $F_1$ and $F_2$,
\item
$W_3$ is the set of tuples $(F_1,F_2,a)$
where $F_1\subset F_2\subset E$ is a
flag such that $F_1\simeq E_1$, $F_2\simeq E_2$
and $a$ is an endomorphism
of the graded vector space $F_1\oplus(E/F_1)$
which preserves the subspace $\{0\}\oplus(F_2/F_1)$.
\end{itemize}

There are obvious inclusions
$W_1\subset X_3\times X_2$, $W_2\subset X_3\times X_1$
and $W_3\subset X_2\times X_1$.
These inclusions factor through an isomorphism
$$W_2=W_1\times_{X_2}W_3.$$
Now, we consider the smooth varieties
$$
Z_1=T^*_{W_{1}}(X_3\times X_2),\quad
Z_2=T^*_{W_{2}}(X_3\times X_1),\quad
Z_3=T^*_{W_{3}}(X_2\times X_1).
$$
The intersection of $(W_1\times X_1)\cap(X_3\times W_2)$ is
transverse in $X_3\times X_2\times X_1$.
Thus, by \cite{CG}, Theorem~2.7.26, the obvious projection
$T(X_3\times X_2\times X_1)\to T^*(X_3\times X_1)$ factors to an isomorphism
$$Z_1\times_{T^*X_2}Z_3\to Z_2.$$
Therefore, we have the following diagram with a Cartesian square
\begin{equation}\label{CO:4}
\xymatrix{
T^*X_3&Z_1\ar[l]_{\psi_1}\ar[r]^-{\phi_1}&T^*X_2\\
&Z_2\ar[ul]^{\psi_2}\ar[u]_\a\ar[r]^-{\b}\ar[dr]_{\phi_2}
&Z_3\ar[u]_{\psi_3}\ar[d]^{\phi_3}\\
&&T^*X_1.}
\end{equation}
The maps are all regular (all the varieties are smooth),
and $\psi_{1}$, $\psi_{3}$ are proper.
Put $\psi_{2}=\psi_{1}\circ\alpha$ and
$\phi_{2}=\phi_{3}\circ\b$. Then $\psi_2$ is also proper.
Therefore we can define the following maps
\begin{equation}\label{CO:5}\gathered
I_{1}=R(\psi_{1})_*\circ L\phi_{1}^*:K^{G\times T}(T^*X_2\on
T^*_GX_2)\to K^{G\times T}(T^*X_3\on T^*_GX_3),\cr
I_{3}=R(\psi_{3})_*\circ L\phi_{3}^*:K^{G\times T}(T^*X_1\on
T^*_GX_1)\to K^{G\times T}(T^*X_2\on T^*_GX_2),\cr
I_{2}=R(\psi_{2})_*\circ L\phi_{2}^*: K^{G\times T}(T^*X_1\on
T^*_GX_1)\to K^{G\times T}(T^*X_3\on
T^*_GX_3).\endgathered\end{equation} Note that $Z_1$, $Z_2$, $Z_3$,
$T^*X_2$ are smooth with $\dim Z_1+\dim Z_3=\dim Z_2+\dim T^*X_2$,
that $\psi_3$ is proper and that $\a\times\b$ is a closed embedding
$Z_2\subset Z_1\times Z_3$. Therefore, by base change we have
$I_{2}=I_{1}\circ I_{3}$, see Proposition \ref{D:1}. Now, set
$$\begin{gathered}
P_{1}=\{g\in G;g(E_1)=E_1\},\quad P=\{g\in
P_{1};g(E_2)=E_2\}.\end{gathered}$$ The Lie algebras of $P_1$, $P$
and their Levi factors are denoted by $\pen_1$, $\pen$, $\hen_1$, $\hen$.
Lemma 7.4 yields
$$T^*_GX_3=C,\quad T^*_GX_2=G\times_{P_{1}}C_{\hen_{1}},\quad
T^*_GX_1=G\times_{P}C_\hen.$$ Thus $I_{1}$, $I_{3}$ yield maps
$$\gathered
I_{1}:K^{G_{E_1}\times T}(C_{1})\otimes_R K^{G_{E/E_1}\times
T}(C_{E/E_1})\to K^{G\times T}(C),\cr I_{3}:K^{H\times T}(C_\hen)\to
K^{G_{E_1}\times T}(C_{1})\otimes_R K^{G_{E/E_1}\times T}(C_{E/E_1}).
\endgathered
$$
The equality $I_{2}=I_{1}\circ I_{3}$ implies that $\bar\Cb_R$
is associative in the same way as in \cite{L}, Lemma~3.4.

\end{proof}

\vspace{.1in}

\begin{rem}
Although $Z$ is smooth the variety $Z_G$ is not locally a complete
intersection in general. This explains why we used $L\phi^*$ rather than
$L(\phi_G)^*$ (which may not be well defined)
in the definition of the map (\ref{CO:0}). Note also that the map
$\phi$ is not flat in general. This explains why we used
the derived functor $L\phi^*$.
\end{rem}

\vspace{.2in}

\paragraph{\textbf{7.5.}}
Next, we modify slightly the construction in Section 7.3 in order to
get a $\bar\Cb_R$-module. First, let us recall the relation between
the commuting variety and the Hilbert scheme. Set
$$\begin{aligned}
N=N_{E}=\gen^2\times E^*\times E, \quad M=M_{E}=\{(a,b,\varphi,v)\in
N; [a,b]+v\circ\varphi=0\}.
\end{aligned}$$
Here $v\circ\varphi$ is regarded as an element of $\gen$. Let
$N^s\subset N$ be the set of tuples $(a,b,v,\varphi)$ such that the
elements $a^{i_1}b^{i_2}\cdots(v)$ span $E$ as $(i_1,i_2,\dots)$
runs over the set of all tuples of integers $\geqslant 0$. Put
$M^s=M\cap N^s$. The group $G \times T$ acts on $M^s$ as
\begin{equation}\label{E:taction2}
(g,z_1,z_2) \cdot (a,b,\varphi,v)=(z_1 g a g^{-1}, z_2 g b g^{-1}, \varphi g^{-1}, z_1z_2 g v).
\end{equation}

We set $\mathbf M_{E,R}=K^{G\times T}(M^s)$. A vector space
isomorphism $E\simeq E'$ yields an $R$-module isomorphisms $\mathbf
M_{E,R}\simeq \mathbf M_{E',R}.$ Let
$$\mathbf M_R=\ind_E \mathbf M_{E,R},$$
where the limit runs over the groupoid formed by all finite
dimensional vector spaces with their isomorphisms. If $n=\dim \;E$
there is a $G$-torsor
$$M^s\to\Hin,\quad
(a,b,v,\varphi)\mapsto\{p(x,y)\in\mathbb C[x,y];p(a,b)v=0\}.$$
Recall the $R$-module $\mathbf L_R$ from (\ref{I:34}).
The following is now obvious.

\begin{lem}
If $\dim\; E=n$ there are canonical $R$-module isomorphisms $\mathbf
M_{E,R}=K^T(\Hin)$ and $\mathbf M_R=\mathbf L_R$.
\end{lem}

\vspace{.2in}

\paragraph{\textbf{7.6.}}
Now fix a subspace $E_1\subset E$ and we set $E_2=E/E_1$. Let $\pi:E\to
E_2$ be the obvious projection. Let $H$, $P$, $\hen$ and $\pen$ be
as in Section 7.4. We set
$$X'=\gen\times E,\quad Y=\hen\times E_2,\quad V=\pen\times E,
\quad X=G\times_PY,\quad W=G\times_PV.$$
The $G$-action on $E$ is the obvious one, the $P$ action on $E_2$ is
the composition of the obvious map $P\to H$ and the $H$-action on
$E_2$. The $G$-action on $X'$, the $H$-action on $Y$ and the
$P$-action on $V$ are the diagonal ones. We'll also write
$$\begin{gathered}
N_\gen=N,\quad
N_\hen=\hen^2\times E_2^*\times E_2=(\gen_1)^2\times N_2,\quad
N_\pen=\pen^2\times E^*_2\times E,
\cr
M_\gen=M,\quad M_\hen=C_{1}\times M_{2},\quad
M_\pen=N_\pen\cap M.
\end{gathered}$$
Here we have identified $E_2^*$ with a subspace of $E^*$ via the
transpose of $\pi$. We'll apply the general construction in Section
7.2 with the map $p:V\to Y$, $(a,v)\mapsto (a_\hen,\pi(v))$ and the
obvious inclusion $q:V\to X'$. The following is proved as Lemma 7.4.

\begin{lem}
(a) We have $X=G\times_P(\hen\times E_2)$ and
$W=G\times_P(\pen\times E)$. The maps $W\to X$, $W\to X'$ are given
by $(g,a,v)\ \mod \ P\mapsto(g,a_\hen,\pi(v))\ \mod\ P$ and
$(g,a,v)\ \mod\ P\mapsto (gag^{-1},gv).$

(b) We have canonical isomorphisms of $G$-varieties
$$\begin{aligned}
T^*X'=N,\quad T^*X=G\times_P\{(c,a,b,\varphi,v)\in\pen\times N_\hen;
c_\hen=[a,b]+v\circ\varphi\},\quad Z=G\times_PN_\pen.
\end{aligned}$$
For all $(a,b,\varphi,v)\in N_\pen$ we have
$$\begin{gathered}
\phi((g,a,b,\varphi,v)\ \mod\ P)= (g,[a,b]+v\circ\varphi
,a_\hen,b_\hen,\varphi,\pi(v))\ \mod\ P,\cr \psi((g,a,b,\varphi,v)\
\mod\ P)=(gag^{-1},gbg^{-1},(\varphi\circ\pi)g^{-1},gv).
\end{gathered}$$

(c) We have canonical isomorphisms of $G$-varieties
$$\begin{aligned}
T^*_GX=G\times_P M_\hen,\quad T^*_GX'=M,\quad
Z_G=G\times_PM_\pen.
\end{aligned}$$
The maps $\phi_G$, $\psi_G$ are the obvious ones.

\end{lem}

\vskip0.2in

We can now prove the following.

\begin{prop}
There is a representation $\rho:\bar\Cb_{R}\to\End_R(\mathbf
M_{R})$.
\end{prop}

\begin{proof}
We'll abbreviate $\mathbf M_i=\mathbf M_{E_i,R}$ for $i=1,2$. First
we define a $R$-linear homomorphism
\begin{equation}\label{CO:8}\bar\Cb_{1}\otimes_R \mathbf M_{2}\to
\mathbf M_{E,R}.\end{equation} We set $N_\pen^s=N_\pen\cap N^s$,
$M_\pen^s=M_\pen\cap M^s$ and
$$Z^s=Z\cap\psi^{-1}(N^s)=G\times_PN_\pen^s,\quad
Z^s_G=Z_G\cap Z^s=G\times_PM_\pen^s.$$
Note that $N_\pen^s$, $M_\pen^s$, $Z^s$, $Z^s_G$ are open in
$N_\pen$, $M_\pen$, $Z$, $Z_G$.
The map $\psi$ restricts to proper morphisms
$$\psi_s:Z^s\to N^s,\quad\psi_{G,s}:Z^s_G\to M^s.$$
Taking the derived direct image we get a $R$-linear map
\begin{equation}\label{CO:9}R(\psi_s)_*:
K^{G\times T}(Z^s_G)=K^{G\times T}(Z^s\on Z^s_G)
\to K^{G\times T}(N^s\on M^s)= K^{G\times T}(M^s).
\end{equation}
Next, we set
$N_\hen^s=(\gen_{1})^2\times N_{2}^s$,
$M_\hen^s=C_{1}\times M_{2}^s$ and
$$\begin{gathered}
T^*X^s=T^*X\cap (G\times_P(\pen\times N_\hen^s)), \quad
T^*_GX^s=T^*_GX\cap T^*X^s=G\times_PM_\hen^s.
\end{gathered}$$
We have $\phi(Z^s)\subset T^*X^s$ by Lemma 7.8$(b)$. Hence
the restriction of $\phi$ yields morphisms
$$\phi_s:Z^s\to T^*X^s,\quad
\phi_{G,s}:Z^s_G\to T^*_GX^s.$$
The map $\phi_s$ has a finite tor-dimension because $\phi$ is regular.
Hence the derived pull-back $L\phi_s^*$ is well-defined, and it yields a
$R$-linear homomorphism
\begin{equation}\label{CO:10}
L\phi_s^*:K^{G\times T}(T^*_GX^s) =K^{G\times T}(T^*X^s\on T^*_GX^s) \to
K^{G\times T}(Z^s\on Z^s_G)= K^{G\times T}(Z^s_G).
\end{equation}
Finally, by induction we have a canonical isomorphism of $R$-modules
\begin{equation}\label{CO:11}\bar\Cb_{1}\otimes_R \mathbf M_{2}=
K^{G\times T}(T^*_GX^s).\end{equation}
We define the map (\ref{CO:8}) to be the
composition of (\ref{CO:11}), (\ref{CO:10}) and (\ref{CO:9}).

Now we must prove that the map (\ref{CO:8}) defines a $\mathbf
C_R$-action on $\mathbf M_R$. The proof is the same as the proof of
Proposition 7.5. More precisely, fix a flag $E_1\subset E_2\subset
E$ and define $P_1,P, H_1, H$ and their Lie algebras as in loc.~cit.
We define the following varieties
\begin{itemize}
\item
$X_1$ is the set of tuples $(F_1,F_2,a,v)$
where $F_1\subset F_2\subset E$ is a
flag such that $F_1\simeq E_1$, $F_2\simeq E_2$
and $a$ is an endomorphism of the graded vector space
$F_1\oplus(F_2/F_1)\oplus (E/F_2)$, and $v$ is an element of $E/F_2$,
\item
$X_2$ is the set of pairs $(F_1,a,v)$ where $F_1\subset E$ is a
vector subspace isomorphic to $E_1$, $a$ is an endomorphism
of the graded vector space $F_1\oplus(E/F_1)$, and $v$ is an element of $E/F_1$,
\item $X_3=\gen\times E$.
\end{itemize}
Next, we define the following ones
\begin{itemize}
\item
$W_1$ is the set of tuples $(F_1,a,v)$ where $F_1\subset E$ is a
vector subspace isomorphic to $E_1$, $a\in\gen$ preserves $F_1$,
and $v\in E$,
\item
$W_2$ is the set of tuples $(F_1,F_2,a,v)$
where $F_1\subset F_2\subset E$ is a
flag such that $F_1\simeq E_1$, $F_2\simeq E_2$,
$a\in\gen$ preserves $F_1$ and $F_2$, and $v\in E$,
\item
$W_3$ is the set of tuples $(F_1,F_2,a)$
where $F_1\subset F_2\subset E$ is a
flag such that $F_1\simeq E_1$, $F_2\simeq E_2$,
$a$ is an endomorphism
of the graded vector space $F_1\oplus(E/F_1)$
which preserves the subspace $\{0\}\oplus(F_2/F_1)$,
and $v\in E/F_1$.
\end{itemize}
We have canonical inclusions $W_1\subset X_3\times X_2$,
$W_3\subset X_2\times X_1$ and $W_2\subset X_3\times X_1$.
Set
$$Z_1=T^*_{W_{1}}(X_3\times X_2),\quad
Z_2=T^*_{W_{2}}(X_3\times X_1),\quad Z_3=T^*_{W_{3}}(X_2\times
X_1).$$ We have again $W_2=W_1\times_{X_2}W_3$ and
$Z_2=Z_1\times_{T^*X_2}Z_3$. In particular we have also the
commutative diagram (\ref{CO:4}) with a Cartesian square. Therefore,
defining the maps $I_1$, $I_2$, $I_3$ as in (\ref{CO:5}) by base
change we get $I_2=I_1\circ I_3 $. Further, by Lemma 7.8 we have
$$\gathered
I_{1}:K^{G_{1}\times T}(C_{1})\otimes_R K^{G_{E/E_1}\times
T}(M_{E/E_1})\to K^{G\times T}(M),\cr I_{3}:K^{H\times T}(M_\hen)\to
K^{G_{1}\times T}(C_{1})\otimes_R K^{G_{2}\times T}(M_{E/E_1}),
\endgathered
$$
where $M_\hen=C_1\times C_{E_2/E_1}\times M_{E/E_2}$. This implies
that the ring $\bar\Cb_R$ acts on the Abelian group
$$\ind_E K^{G\times T}(M_E).$$

Finally we must prove that $\bar\Cb_R$ acts also on $\mathbf M_R$.
By Lemma 7.8 we have
$$T^*X_3=N,\quad
T^*X_2\subset G\times_{P_1}(\pen_1\times N_{\hen_1}),\quad
T^*X_1\subset G\times_{P}(\pen\times N_{\hen}),
$$
where $N_{\hen_1}=(\gen_1)^2\times N_{E/E_1}$
and $N_{\hen}=(\gen_1)^2\times (\gen_{E_2/E_1})^2\times N_{E/E_2}$.
We set
$$T^*X_3^s=N^s,\quad
T^*X_2^s=T^*X_2\cap(G\times_{P_1}(\pen_1\times N_{\hen_1}^s)),\quad
T^*X_1^s=T^*X_1\cap(G\times_{P}(\pen\times N_{\hen}^s)).
$$
Set also
$$Z_3^s=\psi_3^{-1}(T^*X_2^s),\quad Z_2^s=\psi_2^{-1}(T^*X_3^s),\quad
Z_1^s=\psi_1^{-1}(T^*X_3^s).$$
We still have a commutative diagram with a Cartesian square
$$
\xymatrix{
T^*X_3^s&Z_1^s\ar[l]_{\psi_1}\ar[r]^-{\phi_1}&T^*X_2^s\\
&Z_2^s\ar[ul]^{\psi_2}\ar[u]_\a\ar[r]^-{\b}\ar[dr]_{\phi_2}
&Z_3^s\ar[u]_{\psi_3}\ar[d]^{\phi_3}\\
&&T^*X_1^s.}
$$
We can check that
$Z_2^s=Z_3^s\times_{T^*X_2^s}Z_1^s$.
Then the rest of the proof is as above.
\end{proof}

\vspace{.2in}

\paragraph{\textbf{7.7.}}
We keep the same notation as in the previous section. Note that we
have $\mathbf M_R=\mathbf L_R$ by Lemma 7.7. We'll denote as usual
by $\bar\Cb_\Kb$, $\mathbf{M}_{\Kb}$, etc, the extensions of
scalars from $R$ to $\Kb$. Our next goal is to compare the
representation $\rho$ with the faithful representation
$$\psi:\mathbf H_\Kb\to\End_\Kb(\mathbf L_\Kb)$$
given in Section 3.4. For each finite dimensional vector space $E$
we have a $R$-submodule $\bar\Cb_{E,R}\subset\bar\Cb_{R}$ which
depends only on the dimension of $E$.
Recall that $\U^0 \simeq \H^0_{\Kb}$ is a free polynomial algebra in
the tautological classes $\f_{0,l}=\prod_n
\Psi_l(\boldsymbol{\tau}_{n,n})$, $l\in\Z$. We define, for each $E$,
a projection map
\begin{equation}\label{E:defpiprime}
\begin{split}
\pi_E~: \H^0_{\Kb} \to  R_{G_E \times T} \otimes_{R} \Kb, \quad
\f_{0,l} \mapsto  \Psi_l(E)
\end{split}
\end{equation}
where we also write $E$ for the standard $G$-module.
The collection of maps $\pi_E$ endows
$\bar\Cb_\Kb$ with the structure
of an $\H^0_\Kb$-module :
$$\H^0_\Kb \otimes_\Kb \bar\Cb_{E,\Kb}\to\bar\Cb_{E,\Kb},\quad
(h,u)\mapsto h \, u=\pi_E(h)u,$$ where the multiplication in the right hand side
is the tensor product of a class in equivariant K-theory by a
representation. In Section~6.4 we have also defined a Hecke action
$$\bullet~: \H^0_\Kb \otimes_\Kb  \End_\Kb(\mathbf{L}_\Kb) \to
\End_\Kb(\mathbf{L}_\Kb).$$

\vspace{.1in}

\begin{prop}\label{P:Heckektheory}
The map $\rho:\bar\Cb_{\Kb}\to\End_\Kb(\mathbf L_{\Kb})$
intertwines the $\H^0_\Kb$-module structure on $\bar\Cb_{\Kb}$
with the Hecke action on $\End_\Kb(\mathbf{L}_\Kb)$, i.e., for $h
\in \H^0_\Kb$ and $u \in \bar\Cb_{\Kb}$ we have $\rho(h \, u)=h
\bullet \rho(u).$
\end{prop}

\begin{proof} Fix a flag $E_1\subset E$, $E_2=E/E_1$ with
$\dim\; E_1=n$, $\dim\; E_2=k$. By (\ref{E:defHeckeaction2}) we must check that
\begin{equation}\label{E:claim}
\rho(\Psi_l(E_1)u)(x)=\Psi_l(\tau_{n+k})\otimes\rho(u)(x)
-\rho(u)(\Psi_l(\tau_k)\otimes x),\end{equation}
for each
$$u\in\bar\Cb_1=K^{G_1\times T}(C_1),\quad x\in \mathbf M_2=K^{G_2\times T}(M_2^s)=
K^T(\Hi_k).$$
Here the tensor product is the tensor product of coherent sheaves
over $\Hi_{n+k}$ and $\Hi_k$ respectively.
Recall the diagram
$$
\xymatrix{
T^*_GX^s\ar[d]_g&Z^s_G\ar[l]_-{\phi_{G,s}}\ar[d]_{g'}\ar[r]^{\psi_{G,s}}&
(T^*_GX')^s\ar[d]\\
T^*X^s&Z^s\ar[l]_-{\phi_s}\ar[r]^{\Psi_s}&(T^*X')^s,
}
$$
where
$T^*_GX^s=G\times_P(C_1\times M_2^s)$, $Z^s=G\times_PN_\pen^s$ and
$(T^*_GX')^s=M^s$.
Recall also the induction
$$\Ind:\bar\Cb_{1}\otimes_R\mathbf M_{2}\to
K^{G\times T}(T^*_GX^s).$$
We have
$$\rho(u)(x)=R(\psi_s)_*L\phi_s^*g_*(I_{u,x}),\quad
I_{u,x}=\Ind(u\otimes x).$$
Therefore we have
$$\rho(u)(\Psi_l(\tau_k)\otimes x)=
R(\psi_s)_*L\phi_s^*g_*(\theta'_{\Psi_l(E_2)}\otimes I_{u,x})),$$
where $\theta'_{\Psi_l(E_2)}\in K^{G\times T}(T^*X^s)$
is the class induced from the class of the
representation $\Psi_l(E_2)$ in $R^{P\times T}$.
Note that the $G_2$-module $E_2$ is regarded as a $P$-module via the obvious map
$P\to H=G_1\times G_2$, and that it is equipped with the trivial $T$-action.
Thus, the projection formula yields
\begin{equation}\label{E:formule2}
\rho(u)(\Psi_l(\tau_k)\otimes x)=
R(\psi_s)_*(\theta_{\Psi_l(E_2)}\otimes L\phi_s^*g_*(I_{u,x})),\end{equation}
where $\theta_{\Psi_l(E_2)}\in K^{G\times T}(Z^s)$
is induced from the class of
$\Psi_l(E_2)$ in $R^{P\times T}$.
For a similar reason we have also
\begin{equation}\label{E:formule3}
\rho(\Psi_l(E_1)u)(x)=
R(\psi_s)_*(\theta_{\Psi_l(E_1)}\otimes L\phi_s^*g_*(I_{u,x})),\end{equation}
Finally, we have
\begin{equation}\label{E:formule4}
\Psi_l(\tau_{n+k})\otimes\rho(u)(x)=
\Psi_l(E)R(\psi_s)_*L\phi_s^*g_*(I_{u,x}),\end{equation}
where $\Psi_l(E)$ is identified with its class in $R^{G\times T}$.
The formula (\ref{E:claim}) is a direct consequence of
(\ref{E:formule2}), (\ref{E:formule3}), (\ref{E:formule4}).
It is enough to observe that (\ref{E:formule4}) is equivalent to
the equality
$$\Psi_l(\tau_{n+k})\otimes\rho(u)(x)=
R(\psi_s)_*(\Psi_l(E)L\phi_s^*g_*(I_{u,x}))$$
and that
$\Psi_l(E)=\theta_{\Psi_l(E_1)}+\theta_{\Psi_l(E_1)}$
in $K^{G\times T}(Z^s)$.

\end{proof}

\vspace{.1in}

If $E$ is one-dimensional then
$$\bar\Cb^1_R=\bar\Cb_{E,R}=R_{\C^*\times
T}[C],\quad R_{\C^*\times T}=R[z,z^{-1}].$$ Recall that the symbol
$[C]$ denotes the class of $\Oc_C$ and that $R=R_T$. Write
\begin{equation}
\theta_l=z^l[C].
\end{equation}
We can now compare the representation
$\psi:\mathbf H_\Kb\to\End_\Kb(\mathbf L_\Kb)$ from Section 3.2 with
the representation $\rho:\bar\Cb_\Kb\to\End_\Kb(\mathbf L_\Kb)$ from
Proposition 7.10.

\begin{prop}\label{P:Heckektheory2}
If $E$ is one-dimensional then we have $\rho(\theta_l)=\psi(\mathbf
f_{1,l})$ for each $l\in\Z$.
\end{prop}

\begin{proof}
Let us change the notation. Set $E_2=E/E_1$ where $E_1\subset E$ is
a line and $\dim\; E=n+1$. Let $x_1=z^l[C_{1}]$ and $x_2\in
K^T(\Hin)$. We view $x_1,x_2$ as elements of
$$\bar\Cb^1_R=K^{G_1\times T}(C_1),\quad
\mathbf M_2=K^{G_2\times T}(M_{2}^s).$$
We have
$$\boldsymbol{\tau}_{n+1,n}^l\star x_2\in K^T(\Hi_{n+1})=\mathbf M_{E,R}.$$
We must check that the image $x_3\in\mathbf M_{E,R}$ of $x_1\otimes
x_2$ by the map (\ref{CO:8}) is equal to
$\boldsymbol{\tau}_{n+1,n}^l\star x_2$. In Section 7.6 we have
defined a map $\phi_s:Z^s\to T^*X^s$  which restricts to the map
$$\phi_{G,s}:Z^s_G=G\times_PM_\pen^s\to
T^*_GX^s=G\times_P(C_{1}\times M_{2}^s).$$
Consider the Cartesian square of smooth connected varieties
$$
\xymatrix{
T^*_GX^s\ar[d]_g&Z^s_G\ar[l]_-{\phi_{G,s}}\ar[d]_{g'}\\
T^*X^s&Z^s\ar[l]_-{\phi_s}
.}
$$
The vertical maps are the canonical closed embeddings.
The induction yields an isomorphism
$$\Ind:\bar\Cb_{1}\otimes_R\mathbf M_{2}\to
K^{G\times T}(T^*_GX^s)$$
such that
$x_3=R(\psi_s)_*(y_3)$, where
$$y_3=L\phi_s^*g_*\Ind(x_1\otimes x_2)\in K^{G\times T}(Z_G^s).$$
By Lemma 7.8$(c)$ the variety $Z^s_G$ is the set of tuples
$(a,b,v,\varphi)\in M^s$ such that $a,b\in\pen$ and
$\varphi(E_1)=0$. It is well-known that the stability condition
implies that $\varphi=0$. Thus $Z^s_G$ is a $G$-torsor over
$\Hi_{n+1,n}$. Hence we have
$$\begin{gathered}
\dim\; Z_G^s=\dim\; G+\dim\;\Hi_{n+1,n}=n^2+4n+3,
\cr
\dim\; Z^s=\dim\; G-\dim\; P+\dim\; N_\pen=2n^2+5n+3,
\cr
\dim\; T^*_GX^s=\dim\; G-\dim\; P+\dim\; M_\hen^s=n^2+3n+2,
\cr
\dim\; T^*X^s=2\dim\; X=2n^2+4n+2.
\end{gathered}
$$
Thus Proposition \ref{D:1} yields
$L\phi_s^*\circ g_*=g'_*\circ L\phi_{G,s}^*$.
Therefore we have
$$x_3=R(\psi_{G,s})_*L\phi_{G,s}^*\Ind(x_1\otimes x_2).$$
Thus we are reduced to observe that the canonical isomorphism
$$K^{G\times
T}(Z^s_G)=K^T(Z_{n+1,n})$$
takes the class
$L\phi_{G,s}^*\Ind(x_1\otimes x_2)$
to $\boldsymbol{\tau}^l_{n+1,n}\otimes \pi_2^*x_2$.
Recall that $\boldsymbol{\tau}_{n+1,n}^l$ is the $l$-th power of the
tautological bundle over $\Hi_{n+1,n}$.
Thus the claim follows from the definition of $x_1$ and the formula for
$\phi_{G,s}$ recalled above.

\end{proof}

\vspace{.2in}

\paragraph{\textbf{7.8.}} The $R$-algebra $\bar\Cb_R$ is naturally
$\N$-graded, with the piece $\bar\Cb_R^n$ of degree
$n$ equal to the limit
$$\bar\Cb_R^n=\ind_E \bar\Cb_{E,R},$$
over the groupoid formed by all $n$-dimensional vector spaces with
their isomorphisms. Consider the $R$-subalgebra
$\overline\SCb_R\subset\bar\Cb_R$ generated by $\bar\Cb_R^1$. We'll abbreviate
$$\overline\SCb_R^n=\bar\Cb^n_R\cap\overline\SCb_R,\qquad
\overline\SCb_\Kb=\overline\SCb_R\otimes_R \Kb.$$
We have defined in Section~4.6 a $\Kb$-subalgebra $\mathbf
H_\Kb^{>}\subset\mathbf H_\Kb$. Proposition~\ref{P:Heckektheory2}
bears the following consequence.

\begin{cor}
For $n=1$ the assignment $\theta_l\mapsto\mathbf f_{1,l}$, $l\in\Z$,
defines uniquely a surjective $\Kb$-algebra homomorphism
$\Xi:\overline\SCb_\Kb\to\mathbf H^{>}_\Kb$ such that $\rho=\psi\circ\Xi$.
\end{cor}

\vspace{.2in}

\paragraph{\textbf{7.9.}}
Our next goal is to prove that the map $\Xi$ is indeed
an isomorphism. For each vector space $E$ of dimension $n$ the
direct image by the inclusion $C\subset\gen\times\gen$ yields a
$R_{G\times T}$-module homomorphism
\begin{equation}\label{CO:14}\bar\Cb_R^n\to K^{G\times
T}(\gen\times\gen).\end{equation} We conjecture that (\ref{CO:14})
is an injective map. By the Thomason concentration theorem this
conjecture is equivalent to the following one.

\begin{conj}\label{C:torsion} The
$R_{G\times T}$-module $\bar\Cb^n_R$ is torsion-free.
\end{conj}

\noindent The image of (\ref{CO:14}) is the quotient of
$\bar\Cb_R^n$ by its torsion $R_{G\times T}$-submodule. Let
$\SCb_R^n$ be the image of $\overline\SCb_R^n$ by (\ref{CO:14}) and let $\SCb_R=\bigoplus_{n\geqslant 0}\SCb^n_R.$

\vspace{.1in}

\begin{theo}\label{T:Isoktheory}
(a) The map (\ref{CO:14}) factors to a surjective algebra
homomorphism $\overline\SCb_R\to\SCb_R$.

(b) The map $\Xi$ factors to a $\Kb$-algebra isomorphism
$\Xi~:\SCb_\Kb\to\mathbf H^>_\Kb.$

\end{theo}

\begin{proof}
The map (\ref{CO:14}) yields a surjective $R$-linear
homomorphism $\overline\SCb_R\to\SCb_R$. It is easy to check that the
multiplication on $\overline\SCb_R$ descends to $\SCb_R$. Now, 
set $n=\dim\; E$.  In Section 5.2 we have defined a class
$\Lambda(\Vcb_n)\in\mathbf H _\Kb^>$. Since the map
$\Xi:\overline\SCb_\Kb\to\mathbf H^{>}_\Kb$ is surjective we may fix an
element $\nu_n\in\overline\SCb^n_\Kb$ which maps to $\Lambda(\Vcb_n)$. By
the Thomason concentration theorem, the direct image by the
inclusion $\{0\}\subset\gen^2$ yields an isomorphism
$$\mathrm{Frac}(R_{G\times T})\to \bar\Cb^n_R\otimes_{R_{G\times
T}}\mathrm{Frac}(R_{G\times T}).$$ Hence, for each
$x\in\overline\SCb^n_\Kb$ there are non-zero elements $\a,\b\in R_{G\times
T}$ such that
$$\a\, x=\b\,\nu_n.$$
Now, assume that $\rho(x)=0$. Then we have using Proposition~\ref{P:Heckektheory} and (\ref{E:compatHecke})
$$
\psi(\beta \bullet \Lambda(\boldsymbol{\mathcal{V}}_n))=\beta
\bullet \psi(\Lambda(\Vcb_n))= \b \bullet
\rho(\nu_n)=\rho(\b\,\nu_n) =\rho(\a\, x)=\a\bullet\rho(x)=0.$$ Thus
$\b \bullet \Lambda(\boldsymbol{\mathcal{V}}_n)=0$ because $\psi$ is
faithful and $\b=0$ by Proposition~\ref{P:Hecke2}. Hence $x$ is a
torsion element of $\overline\SCb^n_\Kb$.
\end{proof}

\vspace{.1in}

\noindent Combining the above result and Theorem~\ref{T:Main} we
obtain the following corollary.

\begin{cor}\label{C:Main} There is a graded algebra isomorphism
$\Gamma~:\SCb_{\Kb} \to \U^>$ such that $\Gamma(
\theta x)=\theta \bullet \Gamma(x)$ for any $x
\in\SCb_{\Kb}$ of degree $n$ and $\theta \in R_{GL_n\times T} $.
\end{cor}

\vspace{.2in}

\paragraph{\textbf{7.10.}}
The associativity of the multiplication, proved in Proposition 7.5,
yields a $R$-linear map \begin{equation}\label{CO:16}(\mathbf
C_R^1)^{\otimes n}\to\bar\Cb_R^n.\end{equation} The left hand side is identified
with $R[z_1^{\pm 1},z_2^{\pm 1},\cdots z_n^{\pm 1}]$, see Section
7.7. For  $\ell=(l_1,l_2,\dots, l_n)\in\Z^n$ we set
$$z^\ell=z_1^{l_1}z_2^{l_2}\cdots z_n^{l_n}\in(\bar\Cb_R^1)^{\otimes n}.$$
Let us fix a $n$-dimensional vector space $E$. We have a canonical
isomorphism $$\bar\Cb^n_R=K^{G\times T}(C).$$ Let $B\subset G$ be
a Borel subgroup and $H\subset B$ be a maximal torus. Let
$\ben,\hen$ be the Lie algebras of $B,H$ and $\nen\subset\ben$ be
the nilpotent radical. Let $\theta_\ell$ be the character of $H$
associated with $\ell$. It yields a 1-dimensional representation of
$B$. Set $Z=G\times_B\ben^2$. For each $r\geqslant 0$ we have the
$G\times T$-equivariant vector bundle over $Z$
$$\Lambda_Z^r(\ell)=G\times_B(\ben^2\times\theta_\ell\otimes\Lambda^r\nen^*).$$
Consider the complex
$$\Lambda_Z(\ell)=\Bigl\{\cdots\to q^{-2}t^{-2}\Lambda^2_Z(\ell)\to
q^{-1}t^{-1}\Lambda^1_Z(\ell)\to\Lambda^0_Z(\ell) \Bigr\},
$$
where the differential is given by the following map
$$\ben^2\times\Lambda^{r+1}\nen^*\to\ben^2\times\Lambda^{r}\nen^*,\quad
(a,b,\omega)\mapsto\iota_{[a,b]}\omega.$$ Here $\iota$ denotes the
contraction. The cohomology sheaves of this complex are supported on
the subset $Z_G=G\times_BC_\ben$. Consider the proper map
$$\psi:\ Z\to \gen^2,\ (g,a,b)\
\mod\ B\mapsto (gag^{-1},gbg^{-1}).$$ We have $\psi(Z_G)\subset C$.
Thus the class $\Lambda(\ell)=[\psi_*\Lambda_Z(\ell)]$ belongs to
$\bar\Cb_R^n$. The $R$-module $\overline\SCb_R^n$ is described
by the following proposition. Since we'll not use it, the proof is left to the
reader.

\begin{prop} The following holds
\begin{enumerate}
\item[(a)] for each $\ell\in\Z^n$
the image of $z^\ell$ by the map (\ref{CO:16}) is equal to
$\Lambda(\ell)$,
\item[(b)] the $R$-module $\overline\SCb_R^n$ is spanned by the elements
$\Lambda(\ell)$ with $\ell\in\Z^n$.
\end{enumerate}
\end{prop}

\vspace{.2in}

\section{Higher rank}

\paragraph{\textbf{8.1.}}
Fix integers $r>0$, $n\geqslant 0$. Let $M_{r,n}$ be the moduli
space of framed torsion-free sheaves on $\mathbb P^2$ with rank $r$
and second Chern class $n$ (over $\mathbb C$). More precisely,
closed points of $M_{r,n}$ are isomorphism classes of pairs
$(\mathcal E,\Phi)$ where $\mathcal E$ is a torsion-free sheaf which
is locally free in a neighborhood of $\ell_\infty$ and
$\Phi:\mathcal E|_{\ell_\infty}\to\Oc_{\ell_\infty}^r$ is a framing
at infinity. Here $\ell_\infty=\{[x:y:0]\in\mathbb P^2\}$ is the
line at infinity. Recall that $M_{r,n}$ is a smooth variety of
dimension $2rn$ which admits the following alternative description.
Let $E$ be a $n$-dimensional vector space. As above we'll abbreviate
$G=G_E$, $\gen=\gen_E$. There is an isomorphism of algebraic variety
$M_{r,n}=M_{r,E}^s/G$ where
\begin{equation}
\gathered M_{r,E}^s=\{(a,b,\varphi,v)\in M_{r,E}; (a,b,\varphi,v)\
\text{is\ stable}\},\cr M_{r,E}=\{(a,b,\varphi,v)\in
N_{r,E};[a,b]+v\circ\varphi=0\},\cr
N_{r,E}=\gen^2\times\Hom(E,\C^r)\times\Hom(\C^r,E).\endgathered
\end{equation}
The $G$-action is given by
$g(a,b,\varphi,v)=(gag^{-1},gbg^{-1},\varphi g^{-1},gv)$ and
$(a,b,\varphi,v)$ is stable iff there is no proper subspace
$E_1\subsetneq E$ which is preserved by $a,b$ and contains
$v(\mathbb C^r)$.

\vspace{.2in}

\paragraph{\textbf{8.2.}}
Consider the tori  $T^r=(\C^*)^{r}$ and $T^2=(\C^*)^2$.
Let $R=R_T$ be the complexified Grothendieck ring of $T=T^r\times
T^2$. We have $R=\C[q^{\pm 1},t^{\pm 1},\chi_1^{\pm 1},\dots,\chi_r^{\pm
1}]$ where
\begin{equation}
q(h,z_1, z_2)=z_1^{-1}, \qquad t(h,z_1, z_2)=z_2^{-1},
\qquad \chi_\a(h,z_1, z_2)=h_\a^{-1},
\end{equation}
 and
$h=(h_1,h_2,\dots h_r)$, $\a=1,2,\dots r$. 
In Section~8, we make the following change of notation
$$\Kb=\C(q^{1/2},t^{1/2},e^{1/2}_1, \ldots, e^{1/2}_r)$$
and we extend the scalars of all the algebras $\U,$ $\SCb, \dots$ 
defined in the previous section to $\mathbf{K}$.
We equip $M_{r,n}$ and
$M_{r,E}^s$ with the $T$-action given by
\begin{equation}
(h,z_1,z_2)\cdot(a,b,v,\varphi)=(z_1a,z_2b,vh^{-1},z_1z_2h\varphi).
\end{equation}
This action has a finite number of
isolated fixed points, indexed by the set of $r$-tuples of
partitions with total weight $n$. To such a tuple
$\lambda=(\lambda^{(1)},\lambda^{(2)},\dots\lambda^{(r)})$
corresponds a fixed point $I_{\lambda}$ such that the class
$T_\lambda=[T_{I_\lambda}M_{r,n}]$ in $R$ is given by
\begin{equation}\label{HR:1}
T_\lambda=\sum_{\a,\b=1}^r \chi_\b^{-1}
\chi_\a\Bigl(\sum_{s\in\lambda^{(\a)}}t^{l_{\lambda^{(\b)}}(s)}q^{-a_{\lambda^{(\a)}}(s)-1}+
\sum_{s\in\lambda^{(\b)}}t^{-l_{\lambda^{(\a)}}(s)-1}q^{a_{\lambda^{(\b)}}(s)}\Bigr),
\end{equation}
see \cite{NY}, Theorem~2.11.

\vspace{.2in}

\paragraph{\textbf{8.3.}}
The \textit{tautological bundle} of $M_{r,n}$ is the $T$-equivariant
locally free sheaf $\boldsymbol{\tau}_{n}$ given by
$$\boldsymbol{\tau}_{n}=M_{r,E}^s\times_G E.$$
The character of the $T$-action on
its fiber at the fixed point $I_{\lambda}$ is
\begin{equation}\label{HR:2}
\boldsymbol{\tau}_\lambda=[\boldsymbol{\tau}_{{n}}|_{I_{\lambda}}]=\sum_\a
\sum_{s \in \lambda^{(\a)}} \chi_\a^{-1}
t^{y(s)}q^{x(s)},
\end{equation}
see \cite{NY}, Theorem~2.11 and \cite{VV}, Lemma~6. The classes of the tangent bundles and the tautological bundles
are related by the following equation. For each $\lambda$ we have
\begin{equation}\label{E:tangentrelation}
T_{\lambda}=-(1-q^{-1})(1-t^{-1})\boldsymbol{\tau}_\lambda \otimes \boldsymbol{\tau}_\lambda^*  + \boldsymbol{\tau}_\lambda \otimes W^* + q^{-1}t^{-1} \boldsymbol{\tau}_\lambda^* \otimes W
\end{equation}
where $W=\chi_1^{-1} + \cdots + \chi_r^{-1}$ 
is the class of the tautological representation of the torus $T^r$.

\vspace{.2in}

\paragraph{\textbf{8.4.}}
We can now define the \textit{Hecke correspondence}. It is the
variety
$$M_{r,n,n+1}=Z_{r,E}^s/G,$$
where $Z_{r,E}^s$ is the variety of all tuples $(a,b,\varphi,v,E_1)$
where $(a,b,\varphi,v)\in M_{r,E}^s$ and $E_1\subset E$ is a line
preserved by $a,b$ such that $\varphi(E_1)=0$. The variety
$Z_{r,E}^s$ is a $G$-torsor over $M_{r,n,n+1}$, a smooth variety of
dimension $2rn+r+1$. Further the assignments
$$(a,b,\varphi,v)\mapsto (\bar a, \bar b,\bar\varphi,\bar v),\,
(a,b,\varphi,v)$$ yield a closed immersion $M_{r,n,n+1}\subset
M_{r,n}\times M_{r,n+1}$. Here $\bar a,\bar b\in\gen_{E_2}$,
$\bar\varphi\in\Hom(E_2,\C^r)$ are the induced linear maps and $\bar
v=\pi\circ v$. As before we have set $E_2=E/E_1$. The fixed points
contained in $M_{r,n,n+1}$ are those pairs
$I_{\mu,\lambda}=(I_\mu,I_\lambda)$ for which $\mu\subset\lambda$,
i.e., we have $\mu^{(\a)}\subset\lambda^{(\a)}$ for all $\a$, and
$\mu$, $\lambda$ have total weight $n$, $n+1$ respectively.  The class in $R$ of the fiber of the normal bundle to $M_{r,n,n+1}$ at the point $I_{\mu,\lambda}$ has the following expression
\begin{equation}\label{E:conormal1}
N_{\mu,\lambda}=-(1-q^{-1})(1-t^{-1})\boldsymbol{\tau}_{\mu} \otimes \boldsymbol{\tau}_{\lambda}^* + \boldsymbol{\tau}_{\mu} \otimes W^* + q^{-1}t^{-1} \boldsymbol{\tau}_{\lambda}^* \otimes W -q^{-1}t^{-1}.
\end{equation}
See Appendix B for details.

\vspace{.2in}

\paragraph{\textbf{8.5.}}
Let $\pi_1$, $\pi_2$ be the projections of $M_{r,n}\times M_{r,n+1}$
to $M_{r,n}$, $M_{r,n+1}.$ Over $M_{r,n,n+1}$ there is a natural
surjective map $\phi: \pi_2^*( \boldsymbol{\tau}_{n+1}) \to \pi_1^*(
\boldsymbol{\tau}_{n})$. The kernel sheaf $\mathcal{K}er(\phi)$ is a
line bundle called the \textit{tautological bundle} of $M_{r,n,n+1}$
which we denote by $\boldsymbol{\tau}_{n,n+1}$. Over a $T$-fixed
point $I_{\mu,\lambda}$ its character is
\begin{equation}\label{HR:5}
{\boldsymbol\tau}_{\mu,\lambda}=[\boldsymbol{\tau}_{n,n+1}|_{I_{\mu,\lambda}}]=
\chi_\a^{-1}t^{y(s)}q^{x(s)},
\end{equation}
where $s=\lambda^{(\a)} \backslash \mu^{(\a)}$ is the unique box of
$\lambda$ not contained in $\mu$. We define the Hecke correspondence
$M_{r,n+1,n}$ and the tautological bundle $\boldsymbol{\tau}_{n+1,n}$ over it in
the obvious way.

\vspace{.2in}

\paragraph{\textbf{8.6.}}
We'll abbreviate $M_r=\bigsqcup_{n\geqslant 0}M_{r,n}$. Now, we
apply the formalism in Section 3.1 to $X=M_r$, $G=T$ and $Y=\{pt\}$.
Note that $M_{r,n}$ is not proper but has a finite number of fixed
points. Hence the direct image provides us with an isomorphism
\begin{equation}
i_*: K^T(M_{r,n}^T) \otimes_{R} \Kb = \bigoplus_{\lambda}
\Kb[I_{\lambda}] \to K^T(M_{r,n})
\otimes_{R} {\Kb}
\end{equation}
where $i: M_{r,n}^T \to M_{r,n}$ is the embedding.
We have also
$$K^T(M_{r,n} \times M_{r,m} )\otimes_{R}\Kb=\bigoplus_{\lambda,\mu}
\Kb[I_{\lambda,\mu}].$$ This allows us to define convolution
operations

$$\begin{gathered}\star~:{K}^T(M_{r,n}\times M_{r,m})_\Kb
\otimes {K}^T(M_{r,m} \times M_{r,k})_\Kb \to {K}^T(M_{r,n} \times
M_{r,k})_\Kb,\cr \star~:{K}^T(M_{r,n}\times M_{r,m})_\Kb \otimes
{K}^T(M_{r,m})_\Kb \to {K}^T(M_{r,n})_\Kb.\end{gathered}$$
Therefore, the associative $\Kb$-algebra
\begin{equation}\label{HR:6}
{\mathbf{E}}_\Kb=\bigoplus_{k \in \Z}\prod_n {K}^T(M_{r,n+k} \times
M_{r,n})_\Kb,
\end{equation}
where the product ranges over all integers $n \geqslant 0$ with $n+k
\geqslant 0$, acts on the $\Kb$-vector space
\begin{equation*}
\mathbf{L}_\Kb=\bigoplus_{n \geqslant 0}\mathbf L^n_\Kb,\qquad
\mathbf L^n_\Kb={K}^T(M_{r,n})_\Kb.
\end{equation*}
The integer $k$ yields a $\Z$-grading on ${\mathbf{E}}_\Kb$.

\vspace{.2in}

\paragraph{\textbf{8.7.}}
We'll write
\begin{equation}
\begin{gathered}
\boldsymbol{\tau}^l_{n,n+1}=[(\boldsymbol{\tau}_{n,n+1})^{\otimes l}],
\quad
\boldsymbol{\tau}^{-l}_{n,n+1}=[(\boldsymbol{\tau}^*_{n,n+1})^{\otimes
l}],\quad  l \in\Z_{>0},\cr
\f_{-1,l}=\prod_n\boldsymbol{\tau}^l_{n,n+1},
\quad\f_{1,l}=\prod_n\boldsymbol{\tau}^l_{n+1,n},\quad l\in\Z,\cr
\mathbf{e}_{0,l}=\prod_n \Lambda^l \boldsymbol{\tau}_{n,n},
\quad \mathbf{e}_{0,-l}=\prod_n \Lambda^{l}
\boldsymbol{\tau}^*_{n,n},\quad l\in\Z_{>0}.
\end{gathered}
\end{equation}
Once again, we define the elements $\f_{0,l}\in\mathbf E_\Kb$ for $l
\in \Z^*$ through the relations
$$\sum_{k \geqslant 1}\f_{0,\pm k} s^{k-1}=-\frac{d}{ds}\log(E_{\pm}(s)),\quad
E_{\pm}(s)=1+\sum_{k \geqslant 1} (-1)^{k} \mathbf{e}_{0,\pm
k}s^k.$$ So $\f_{0,l}$ is obtained from the classes of the
tautological bundles $\boldsymbol{\tau}_{n,n}$ by the Adams
operations
$$\f_{0,k}=\prod_n\Psi_k(\boldsymbol{\tau}_{n,n}),\qquad
\f_{0,-k}=\prod_n\Psi_k(\boldsymbol{\tau}_{n,n}^*).$$
Finally we introduce some elements $\mathbf h_{i,l}$, with $i=-1,0,1$,
$l\in\Z$ and $(i,l)\neq(0,0)$, defined as
\begin{equation}\label{E:hrd12}
\hb_{1,l}=t^{1-r/2}\f_{1, l-r}, \qquad \hb_{-1,l}=
(-1)^r\det(W)q^{1-r/2}\f_{-1, l},
\end{equation}
\begin{equation}\label{E:hrd22}
\hb_{0,k}= \f_{0,k}-\frac{p_k(\chi_1^{-1},\dots,\chi_r^{-1})}{(1-q^k)(1-t^k)},
\qquad
\hb_{0,-k}=-\f_{0,-k}+
\frac{p_k(\chi_1,\dots,\chi_r)}{(1-q^{-k})(1-t^{-k})},\qquad k\geqslant 1,
\end{equation}
where $\det(W)=(\chi_1 \cdots \chi_r)^{-1}$.
Consider the $\Kb$-subalgebra $\H_\Kb\subset{\mathbf{E}}_\Kb$
generated by these elements. The
$\Kb$-subalgebras $\mathbf H^>_\Kb$, $\mathbf H^0_\Kb$ and $\mathbf
H^<_\Kb$ are defined as in Section 4.6. We have a faithful
representation $\psi:\H_\Kb\to\End_\Kb(\mathbf{L}_\Kb)$. Compare
Section 3.2. Consider the central charge $$c^r=(1,q^{r/2}t^{r/2}).$$

\begin{theo}\label{HR:Main}  There is a $\Kb$-algebra isomorphism
$\U_{c^r}\to\H_\Kb$ such that $\bt_{i,l}\mapsto\hb_{i,l}$ for all
$i=-1,0,1$ and $l \in \Z$.
\end{theo}

The proof is similar to the proof of
Theorem 3.1. Let us just recall the main arguments.

\vspace{.2in}

\paragraph{\textbf{8.8.}} 
First we construct a $\Kb$-algebra homomorphism
$\Xi:\SCb_\Kb\to\mathbf H_\Kb^>$ as in Proposition 7.9.
Composing it with the isomorphism $\U^>\simeq\SCb_\Kb$ in
Theorem 3.1 and Theorem 7.14 we get a $\Kb$-algebra homomorphism
\begin{equation}\label{8.8:isom1}\U^>\to\mathbf H^>_\Kb.
\end{equation}
To define $\Xi$ we first construct a representation $\rho$ of $\bar\Cb_\Kb$ on
$\mathbf L_\Kb$. To do that we let $E_1$, $E_2$, $\pi$, $H$, $P$,
$\hen$ and $\pen$ be as in Sections 7.4, 7.6, and we set
$$\gathered
X'=\gen\times \Hom(\C^r,E),\quad Y=\hen\times \Hom(\C^r,E_2),\quad
V=\pen\times \Hom(\C^r,E),\cr X=G\times_PY,\quad
W=G\times_PV.\endgathered$$ The $G$-action on $X'$, the $H$-action
on $Y$ and the $P$-action on $V$ are the obvious ones. Set also
$$\begin{gathered}
N_\hen=\hen^2\times \Hom(E_2,\C^r)\times \Hom(\C^r,E_2),\quad
N_\pen=\pen^2\times \Hom(E_2,\C^r)\times \Hom(\C^r,E), \cr
N_\gen=N_{r,E},\quad M_\gen=M_{r,E},\quad M_\hen=C_{E_1}\times
M_{r,E_2},\quad M_\pen=N_\pen\cap M_\gen.
\end{gathered}$$
Then we apply the general construction in Section 7.2 with the map
$p:V\to Y$, $(a,v)\mapsto (a_\hen,\pi\circ v)$ and the obvious
inclusion $q:V\to X'$. We have the following lemma.

\begin{lem}
(a) We have canonical isomorphisms of $G$-varieties
$$\begin{aligned}
T^*X'=N_\gen,\quad T^*X=G\times_P\{(c,a,b,\varphi,v)\in\pen\times
N_\hen; c_\hen=[a,b]+v\circ\varphi\},\quad Z=G\times_PN_\pen.
\end{aligned}$$
For all $(a,b,\varphi,v)\in N_\pen$ we have
$$\begin{gathered}
\phi((g,a,b,\varphi,v)\ \mod\ P)=
(g,[a,b]+v\circ\varphi,a_\hen,b_\hen,\varphi,\pi\circ v)\ \mod\
P,\cr \psi((g,a,b,\varphi,v)\ \mod\
P)=(gag^{-1},gbg^{-1},(\varphi\circ\pi)g^{-1},gv).
\end{gathered}$$

(b) We have canonical isomorphisms of $G$-varieties
$$\begin{aligned}
T^*_GX=G\times_P M_\hen,\quad T^*_GX'=M_\gen,\quad
Z_G=G\times_PM_\pen.
\end{aligned}$$
The maps $\phi_G$, $\psi_G$ are the obvious ones.

\end{lem}

\noindent Note that $K^T(M_{r,n})=K^{G\times T}(M_\gen^s)$. Thus the
proof of Proposition 7.9 and the formulas above yield the
representation we need
\begin{equation}
\rho:\bar\Cb_\Kb\to\End_\Kb(\mathbf L_\Kb).
\end{equation}
Next, the proof of Proposition 7.11 implies that if $E$ is
one-dimensional then we have 
\begin{equation}
\rho(\theta_l)=\psi(\mathbf f_{1,l}),\qquad l\in\Z.
\end{equation}
Since the representation 
$\psi:\mathbf H_\Kb\to\End_\Kb(\mathbf L_\Kb)$ 
from Section 3.2 is faithful, this
yields a $\Kb$-algebra homomorphism 
\begin{equation}
\label{toto12}
\overline\SCb_\Kb\to\mathbf H_\Kb^>,\qquad \theta_l\mapsto \mathbf f_{1,l}.
\end{equation}
Now, we prove the following.

\begin{lem} The map \eqref{toto12}
factors to an injective $\Kb$-algebra homomorphism
\begin{equation}\Xi:\SCb_\Kb\to\mathbf H_\Kb^>.
\end{equation}
\end{lem}

\begin{proof} We must prove that the map \eqref{toto12}
factors to an algebra homomorphism
$\Xi:\SCb_\Kb\to\mathbf H_\Kb^>$
and that $\Xi$ is injective.
We'll be very brief. The proof is similar to the proof of
\cite{SV5}, Theorem 6.3.
Details are left to the reader.
Set $G=GL_n$. We have
$$R_{G\times T}\otimes_R\Kb=\Kb[z_1^{\pm 1}, \ldots, z_n^{\pm 1}]^{\mathfrak{S}_n}=\Sb^n_\Kb.$$
We'll view $\Lb_\Kb^n$ as a $\Kb$-subalgebra of $\Sb_\Kb^n$ in the obvious way.
Recall that there is a unique $\Kb$-algebra homomorphism $\Lb_\Kb\to\Lb^n_\Kb$ taking
$p_{l}=z_1^l+z_2^l+\cdots$ to $p_{l}(z_1,\dots,z_n)=z_1^l+\cdots+z_n^l$ for each $l>0$.
The tensor product in equivariant K-theory gives an Hecke action
\begin{equation}\Lb_\Kb\otimes_\Kb\bar\Cb_\Kb\to\bar\Cb_\Kb.
\end{equation}
It is not difficult to prove that this action factors to an action on $\overline\SCb_\Kb$, $\SCb_\Kb.$ 
Further, it factors through the canonical quotient $\Lb_\Kb\to\Lb^n_\Kb$
on the degree $n$ pieces $\overline\SCb_\Kb^n,$
$\SCb_\Kb^n.$
Next, there is an Hecke action
\begin{equation}\Lb_\Kb\otimes_\Kb\mathbf H_\Kb^>\to\mathbf H_\Kb^>.
\end{equation}
More precisely,  the action of the  power sum 
$p_{l}$ is given by the formula
\begin{equation}
p_l\bullet u=[p_l,u],\qquad u\in\End_\Kb(\mathbf L_\Kb),
\end{equation}
where $p_l$ acts on $\mathbf L^n_\Kb$ through the tensor product by $\Psi_l(\tau_n)$.
Once again, it is not difficult to prove that this action preserves $\mathbf H_\Kb^>$,  
and that the action on the degree $n$ piece
$\mathbf H_\Kb^>[n]$ factors through the canonical quotient $\Lb_\Kb\to\Lb^n_\Kb$.
Further, we have
\begin{equation}
p_l\bullet\mathbf f_{1,k}=\mathbf f_{1,k+l},\qquad k\in\Z.
\end{equation}

Now, we claim that the map \eqref{toto12} intertwines both Hecke actions and that the action of 
$\Lb^n_\Kb$ on
$\mathbf H_\Kb^>[n]$ is torsion free. This implies that the map $\Xi$ is well defined.
The injectivity of $\Xi$ is a consequence of the surjectivity of $\Xi$, of the commutation of $\Xi$ with
the Hecke actions, and of the localization theorem in K-theory,
see \cite{SV5}, Theorem 6.3, for more details.  Write $\Xi'$ for the map \eqref{toto12}.
The first part of the claim follows from the following relations
\begin{equation}
p_l\bullet\Xi'(\theta_k)=p_l\bullet\mathbf f_{1,k}=\mathbf f_{1,k+l}=\Xi'(\theta_{k+l})=\Xi'(p_l\bullet\theta_k).
\end{equation}
To prove the second part of the claim, we consider the surjective $\Kb$-linear map
\begin{equation}
\iota: \Kb[z_1^{\pm 1}, \ldots, z_n^{\pm 1}]\to\mathbf H_\Kb^>,\qquad
z_1^{k_1}\cdots z_n^{k_n}\mapsto\mathbf f_{1,k_1}\cdots\mathbf f_{1,k_n}.
\end{equation}
A direct computation shows that
\begin{equation}
\iota(P)=\varpi_n(P)\bullet\gamma_n,
\end{equation}
where $\gamma_n=\iota(1)$  and $\varpi_n$ is the \emph{twisted symmetrization map}
\begin{equation}
\varpi_n:\begin{cases}
\Kb[z_1^{\pm 1}, \ldots, z_n^{\pm 1}]\to\Kb(z_1, \ldots, z_n)^{\mathfrak{S}_n}\\
P(z_1,\dots,z_n)\mapsto\text{SYM}_n\big(g(z_1,\dots,z_n)\,P(z_1,\dots,z_n)\big).
\end{cases}
\end{equation}
Here $\text{SYM}_n$ and $g(z_1,\dots,z_n)$ are as in Section 10.1.
Therefore, to check the second part of the claim it is enough to observe that the map
\begin{equation}
\Lb^n_\Kb\to\End_\Kb(\mathbf L_\Kb),\qquad
p\mapsto p\bullet\gamma_n
\end{equation}
is injective. This is left to the reader, compare \cite[prop.~3.12]{SV5}.

\end{proof}

\vspace{.2in}

\paragraph{\textbf{8.9.}} 
In the same way we define an injective $\Kb$-algebra homomorphism
\begin{equation}\label{8.8:isom2}\U^<\to\mathbf H^<_\Kb.
\end{equation}

\vspace{.2in}

\paragraph{\textbf{8.10.}} 
Finally one checks that the class $[\mathbf f_{1,l},\mathbf
f_{-1,k}]$ is supported on the diagonal of $M_{r,n}\times M_{r,n}$
as in \cite{VV}, Lemma~9, and this contribution is computed as in 
\cite{VV}, Section~4.5. More precisely, we have the following.

\begin{prop}
(a) The class $[\mathbf h_{-1,k},\mathbf h_{1,l}]$ 
is supported on the diagonal of $M_{r,n}\times M_{r,n}.$

(b) For $k+l>0$ we have 
$[\mathbf h_{-1,k},\mathbf h_{1,l}]=q^{r/2}t^{r/2}\theta_{0,k+l}/\alpha_1$ with
$$\sum_{l\geqslant 0}\theta_{0,l}\,s^l=
\exp\Big(\sum_{k\geqslant 1}\alpha_k\mathbf h_{0,k}\,s^k\Bigr).$$

(c) For $k+l<0$ we have 
$[\mathbf h_{-1,k},\mathbf h_{1,l}]=-q^{-r/2}t^{-r/2}\theta_{0,k+l}/\alpha_1$ with
$$\sum_{l\geqslant 0}\theta_{0,-l}\,s^l=
\exp\Big(\sum_{k\geqslant 1}\alpha_k\mathbf h_{0,-k}\,s^k\Bigr).$$

(d) For $k+l=0$ we have 
$[\mathbf h_{-1,k},\mathbf h_{1,l}]=(q^{r/2}t^{r/2}-q^{-r/2}t^{-r/2})/\alpha_1$.
\end{prop}

\begin{proof}
First, we have the following formulas, compare Section 4.1,
\begin{equation}\label{Er:31}
\pmb\tau_{n,n+1}^l=\sum_{\mu \subset \lambda} 
\pmb\tau_{\mu,\lambda}^l \cdot
\Lambda(N^*_{\mu,\lambda})\cdot\Lambda_{\mu,\lambda}^{-1}\cdot 
[I_{\mu,\lambda}],
\end{equation}
\begin{equation}\label{Er:32}
\pmb\tau_{n+1,n}^l=\sum_{\mu\subset\lambda} 
\pmb\tau_{\lambda,\mu}^l \cdot
\Lambda(N^*_{\lambda,\mu})\cdot\Lambda_{\lambda,\mu}^{-1}\cdot[I_{\lambda,\mu}],
\end{equation}
\begin{equation}\label{Er:33}
\pmb\tau_{n,n}^l=\sum_{\mu} 
\pmb\tau_{\mu}^l\cdot \Lambda_{\mu}^{-1} \cdot [I_{\mu,\mu}].
\end{equation}
Here $\mu$ and $\lambda$ are $r$-partitions 
of $n$ and $n+1$ respectively.

Now, assume that $\lambda$, $\mu$ are $r$-partitions 
of $n$ and that $\sigma$, $\pi$ are $r$-partitions 
of $n-1$, $n+1$ respectively, with
$\sigma\subset \lambda,$ $\mu\subset\pi$ and
$\lambda\neq\mu$. Then, the $r$-partitions
$\sigma$, $\pi$ are completely determined by $\lambda$, $\mu$ and 
\eqref{HR:2} gives
the following identity
$$\pmb\tau_\lambda+\pmb\tau_\mu=\pmb\tau_\sigma+\pmb\tau_\pi.$$
Therefore, using the identities from Sections 8.2, 8.3 and 8.4,
a short computation gives
\begin{equation}
\label{tutu2}
N_{\lambda,\sigma}+N_{\mu,\sigma}-T_\sigma= 
N_{\pi,\lambda}+N_{\pi,\mu}-T_\pi.
\end{equation}
Recall that we have, see \eqref{E:hrd12},
$$\hb_{1,l}=t^{1-r/2}\f_{1, l-r}, \qquad 
\hb_{-1,k}=(-1)^r\det(W)q^{1-r/2}\f_{-1, k}.$$
Therefore, using \eqref{Er:31}, \eqref{Er:32} and 
\eqref{tutu2} we get 
\begin{equation}
\label{tutu3}
[\mathbf h_{-1,k},\mathbf h_{1,l}]\star[I_\lambda]=c_{\lambda,k+l}\,[I_\lambda]
\end{equation}
for some constant $c_{\lambda,k+l}$ which remains to be computed. 
To do so, observe first that
$$\gathered
\mathbf f_{-1,k}\mathbf f_{1,l-r}\star[I_\lambda]=\sum_{\lambda\subset\pi}
\pmb\tau_{\lambda,\pi}^{k+l-r}\cdot
\Lambda(N_{\lambda,\pi}^*+N_{\pi,\lambda}^*)\cdot
\,\Lambda_{\lambda,\pi}^{-1}\cdot [I_\lambda],\\
\mathbf f_{1,l-r}\mathbf f_{-1,k}\star[I_\lambda]=\sum_{\sigma\subset\lambda}
\pmb\tau_{\sigma,\lambda}^{k+l-r}\cdot
\Lambda(N_{\lambda,\sigma}^*+N_{\sigma,\lambda}^*)\cdot
\,\Lambda_{\lambda,\sigma}^{-1}\cdot[I_\lambda],
\endgathered$$
modulo $[I_\mu]$'s with $\mu\neq\lambda$. 
Next, set $H_\lambda=(1-q)(1-t)\pmb\tau_\lambda-W.$ 
For $\lambda\subset\pi$, we have
$$\aligned
H_\lambda&=H_\pi-(1-q)(1-t)\pmb\tau_{\lambda,\pi},\\
N_{\lambda,\pi}-T_\lambda
&=-q^{-1}t^{-1}\pmb\tau_{\lambda,\pi}^* H_{\lambda}-q^{-1}t^{-1},\\
&=-q^{-1}t^{-1}\pmb\tau_{\lambda,\pi}^* H_{\pi}+1-q^{-1}-t^{-1},\\
N_{\pi,\lambda}-T_\pi
&=\pmb\tau_{\lambda,\pi}H_\pi^*-q^{-1}t^{-1},\\
&=\pmb\tau_{\lambda,\pi}H_{\lambda} ^*+1-q^{-1}-t^{-1}.
\endaligned$$
Now, we consider the following sums
$$B_\lambda=\sum_{\sigma\subset\lambda}\pmb\tau_{\sigma,\lambda},
\quad A_\lambda=\sum_{\lambda\subset\pi}\pmb\tau_{\lambda,\pi}.$$
The proof of the following lemma is left to the reader, 
compare \cite{VV}, Lemma 7.

\begin{lem}
\label{lem:toto1}
For each $r$-partition $\lambda$ of $n$ we have 
$H_\lambda=qtB_\lambda-A_\lambda.$
\end{lem}

\noindent
For $\lambda\subset\pi$, using the lemma above, we get the following formulas 
$$\aligned
N_{\lambda,\pi}^*+N_{\pi,\lambda}^*-T_\lambda^*-T_\pi^*
&=qt(\pmb\tau_{\lambda,\pi}A_{\lambda}^*+\pmb\tau_{\lambda,\pi}^*B_\lambda-1)-
(\pmb\tau_{\lambda,\pi}B_{\lambda}^*+
\pmb\tau_{\lambda,\pi}^*A_\lambda-1)
-q-t,\\
N_{\lambda,\sigma}^*+N_{\sigma,\lambda}^*-T_\lambda^*-T_\sigma^*
&=qt(\pmb\tau_{\sigma,\lambda}A_{\lambda}^*+
\pmb\tau_{\sigma,\lambda}^*B_\lambda-1)-
(\pmb\tau_{\sigma,\lambda}B_{\lambda}^*+\pmb\tau_{\sigma,\lambda}^*A_\lambda-1)-q-t.
\endaligned$$
Therefore, we get also
$$
\gathered
(-1)^r\det(W)^{-1}(1-q)(1-t)(qt)^{r/2-1}\,c_{\lambda,k}=
\sum_{\lambda\subset\pi}\pmb\tau_{\lambda,\pi}^{k-r}\,
\frac{\Lambda(qt\pmb\tau_{\lambda,\pi} A_{\lambda}^*+
qt\pmb\tau_{\lambda,\pi}^*B_\lambda-qt)}
{\Lambda(\pmb\tau_{\lambda,\pi}^*A_{\lambda}+
\pmb\tau_{\lambda,\pi}B_\lambda^*-1)}-\\
-\sum_{\sigma\subset\lambda}\pmb\tau_{\sigma,\lambda}^{k-r}\,
\frac{\Lambda(qt\pmb\tau_{\sigma,\lambda} A_{\lambda}^*+
qt\pmb\tau_{\sigma,\lambda}^*B_\lambda-qt)}
{\Lambda(\pmb\tau_{\sigma,\lambda}^* A_{\lambda}+
\pmb\tau_{\sigma,\lambda}B_\lambda^*-1)}.
\endgathered
$$
Write $a_\lambda=\prod_{i\in I}a_i$ and
$b_\lambda=\prod_{j\in J}b_j$, where
$$\{a_i\;;\;i\in  I\}=\{\pmb\tau_{\lambda,\pi}\;;\;\lambda\subset\pi\},
\qquad
\{b_j\;;\;j\in J\}=\{\pmb\tau_{\sigma,\lambda}\;;\;\sigma\subset\lambda\}.$$
Consider the formal series 
$$C(s)=\det(W)^{-1} (1-q)(1-t)(qt)^{r/2-1}\,
\sum_{k\geqslant 0}c_{\lambda,k}\,s^k.$$ 
A short computation yields
$$\aligned
-a_\lambda\, b_\lambda^{-1}C(s)&=
\sum_{i\in I} \frac{1}{1-a_is}\,
\prod_{k\in I\setminus\{i\}} \frac{1-qta_i/a_k}{1-a_i/a_k}
\prod_{j\in J}\frac{1-qtb_j/a_i}{1-b_j/a_i}-\\
&\qquad-\sum_{j\in J} \frac{1}{1-b_js}\, 
\prod_{i\in I} \frac{1-qtb_j/a_i}{1-b_j/a_i}
\prod_{k \in J\setminus\{j\}} \frac{1-qtb_k/b_j}{1-b_k/b_j}.
\endaligned
$$
Now, we have the following lemma.

\begin{lem}
Given commutative formal variables
$\{a_i\;;\;i\in I\}$ and $\{b_j\;;\;j\in J\}$ we have
$$
\aligned
& \sum_{i\in I} \frac{1}{1-a_is}\, 
\prod_{k\in I\setminus\{i\}}\frac{1-qta_i/a_k}{1-a_i/a_k}
\prod_{j\in J} \frac{1-qtb_j/a_i}{1-b_j/a_i}-
\sum_{j\in J}\frac{1}{1-b_js}\,\prod_{i\in I}\frac{1-qtb_j/a_i}{1-b_j/a_i}
\prod_{k\in J\setminus\{j\}}\frac{1-qtb_k/b_j}{1-b_k/b_j}=\\
&=\frac{1}{qt-1}\prod_{i\in I} \frac{qt-a_is}{1-a_is}
\prod_{j\in J}\frac{1-qtb_js}{1-b_js}-\frac{(qt)^{\sharp J}}{qt-1}.
\endaligned
$$
\end{lem}

\begin{proof} 
It is enough to prove the following identity
\begin{equation*}
\begin{split}
& \sum_{i\in I} \frac{a_is}{1-a_is}\, 
\prod_{k\in I\setminus\{i\}}\frac{1-qta_i/a_k}{1-a_i/a_k}
\prod_{j\in J} \frac{1-qtb_j/a_i}{1-b_j/a_i}-
\sum_{j\in J}\frac{b_js}{1-b_js}\,\prod_{i\in I}\frac{1-qtb_j/a_i}{1-b_j/a_i}
\prod_{k\in J\setminus\{j\}}\frac{1-qtb_k/b_j}{1-b_k/b_j}=\\
&=\frac{1}{qt-1}\prod_{i\in I} \frac{qt-a_is}{1-a_is}
\prod_{j\in J}\frac{1-qtb_js}{1-b_js}-\frac{(qt)^{\sharp I}}{qt-1}.
\end{split}
\end{equation*}
Both sides of the equality are rational functions in $s$ of degree $0$, with at most simple poles. 
One checks that the poles and residues are the same. This implies the equality, 
up to a possible constant. But both sides vanish at $s=0$. 
So this constant is zero.
\end{proof}

Using the lemma above, we get
$$\aligned
-a_\lambda b_\lambda^{-1}C(s)
&=\frac{(qt)^{\sharp I}}{qt-1}\prod_{\lambda\subset\pi} \frac{1-(qt)^{-1}\tau_{\lambda,\pi}s}
{1-\tau_{\lambda,\pi}s}\,
\prod_{\sigma\subset\lambda}\frac{1-qt\tau_{\sigma,\lambda}s}
{1-\tau_{\sigma,\lambda}s}-\frac{(qt)^{\sharp J}}{qt-1}.
\endaligned$$
Now, fix splitting sums of one-dimensional characters
$$\pmb\tau_\lambda^*=\phi_1+\cdots+\phi_n,\qquad
W^*=\chi_1+\cdots+\chi_r.$$ 
By Lemma \ref{lem:toto1}, we have
$$\gathered
H_\lambda=\sum_{\sigma\subset\lambda}qt\pmb\tau_{\sigma,\lambda}-
\sum_{\lambda\subset\pi}\pmb\tau_{\lambda,\pi}=
\sum_i(1-q)(1-t)\phi_i^{-1}-\sum_\alpha \chi_\alpha^{-1},
\qquad
(qt)^{\sharp J}a_\lambda^{-1}b_\lambda=\det(W)^{-1}.
\endgathered$$
Therefore, we get
$$
\aligned
-\det(W)C(s)
&=\frac{(qt)^r}{qt-1}
\prod_{\lambda\subset\pi} \frac{1-(qt)^{-1}\tau_{\lambda,\pi}s}
{1-\tau_{\lambda,\pi}s}\,
\prod_{\sigma\subset\lambda}\frac{1-qt\tau_{\sigma,\lambda}s}
{1-\tau_{\sigma,\lambda}s}
-\frac{1}{qt-1},\\
&=\frac{(qt)^r}{qt-1}
\Bigl(\sum_{k\geqslant 0}(-s)^k\Lambda^kH_\lambda\Bigr)\,
\Bigl(\sum_{k\geqslant 0}(-s/qt)^k\Lambda^kH_\lambda\Bigr)^{-1}
-\frac{1}{qt-1},\\
&=\frac{(qt)^r}{qt-1}
\prod_{i=1}^n
\frac{1-qt\phi_i^{-1}s}{1-(qt)^{-1}\phi_i^{-1}s}\,
\frac{1-q^{-1}\phi_i^{-1}s}{1-q\phi_i^{-1}s}\,
\frac{1-t^{-1}\phi_i^{-1}s}{1-t\phi_i^{-1}s}\,
\,\prod_{\alpha=1}^r\frac{1-(qt\chi_\alpha)^{-1}s}{1-\chi_\alpha^{-1}s}
-\frac{1}{qt-1}.
\endaligned
$$
We abbreviate $p_l(\phi_i^{-1})=p_l(\phi_1^{-1},\dots,\phi_n^{-1})$
and $p_l(\chi_\alpha^{-1})=p_l(\chi_1^{-1},\dots,\chi_r^{-1})$.
Observe that
$$\prod_{i=1}^n\frac{1-q^{-1}\phi_i^{-1}s}{1-q\phi_i^{-1}s}=
\exp\Big(\sum_{l\geqslant 1}(q^{l}-q^{-l})p_l(\phi_i^{-1})s^l/l\Big).$$
We finally get
$$\aligned
\alpha_1\sum_{k\geqslant 0}c_{\lambda,k}\,s^{k}&=(qt)^{r/2}
\exp\Bigl(\sum_{l\geqslant 1}
\alpha_l\,\Big(\frac{p_l(\chi_\alpha^{-1})}{(1-q^l)(1-t^l)}
-p_l(\phi_i^{-1})\Big)\,s^l\Bigr)
-(qt)^{-r/2},\\
\alpha_l=\alpha_{-l}&=
\big((qt)^{l/2}-(qt)^{-l/2}\big)
\big(q^{l/2}-q^{-l/2}\big)
\big(t^{l/2}-t^{-l/2}\big)/l\\
&=\big(q^{-l}-q^l\big)/l+\big(t^{-l}-t^l\big)/l+
\big((qt)^l-(qt)^{-l}\big)/l,\qquad l\geqslant 1.
\endaligned
$$
Now, recall that for $l\geqslant 1$ we have, see \eqref{E:hrd22},
$$\gathered
\hb_{0,l}= \f_{0,l}-\frac{p_l(\chi_\alpha^{-1})}{(1-q^l)(1-t^l)},
\qquad
\hb_{0,-l}=-\f_{0,-l}+\frac{p_l(\chi_\alpha)}{(1-q^{-l})(1-t^{-l})},\\
\mathbf f_{0,\pm l}\star[I_\lambda]=
p_l(\phi_1^{\mp 1},\dots,\phi_n^{\mp 1})\,[I_\lambda].
\endgathered$$ 
Thus, we obtain
$$\aligned
\alpha_1\sum_{k\geqslant 0}c_{\lambda,k}\,s^{k}\,[I_\lambda]
&=(qt)^{r/2}\,
\exp\Bigl(-\sum_{l\geqslant 1}
\alpha_{l}\,\mathbf h_{0,l}\,s^l\Bigr)\star[I_\lambda]-(qt)^{-r/2}[I_\lambda].
\endaligned
$$
Comparing this expression with \eqref{tutu3}, we get parts
$(b)$ and $(d)$ of the proposition. Part $(c)$ is similar and is left to
the reader.
\end{proof}

\vspace{.1in}

Therefore, we have proved the following.

\begin{cor} There is a $\Kb$-algebra homomorphism
$\U_{c^r}\to\H_\Kb$ such that $\bt_{i,l}\mapsto\hb_{i,l}$ for all
$i=-1,0,1$ and $l \in \Z$.
\end{cor}

\vspace{.2in}

\paragraph{\textbf{8.11.}} 
To prove that the $\Kb$-algebra homomorphism
$\U_{c^r}\to\H_\Kb$ above is injective, we are now reduced to prove the 
following.

\vspace{.1in}

\begin{prop} The algebra $\H_\Kb$ has a triangular decomposition,
i.e., the multiplication map induces an isomorphism
$m~: \H_\Kb^> \otimes \H_\Kb^0 \otimes \H_\Kb^< \to \H_\Kb.$
\end{prop}

\begin{proof} The proof is similar to the proof of Proposition \ref{T:2}.
For an operator $f$ on
$\mathbf{L}_\Kb$ and $r$-partitions $\lambda,$ $\mu$
we denote by $\langle \mu , f , \lambda\rangle$ the
coefficient of $[I_{\mu}]$ in $f([I_{\lambda}])$.
Given partitions $\lambda_1, \lambda_2, \ldots, \lambda_k$ and given an 
integer $n \gg |\lambda_1|,\ldots, |\lambda_k|,$ let the symbol
$\lambda_1\circledast\ldots\circledast\lambda_k$ denote the $r$-partition
whose first part is the partition $\lambda_1\circledast\ldots\circledast\lambda_k$
from Section 4.6  and the $r-1$ other partitions are empty.
Given a finite family of elements 
$$P_i\in\H^{>}_\Kb,\qquad R_i\in\H^{0}_\Kb,\qquad Q_i\in\H^{<}_\Kb,$$
we set $x=\sum_iP_iR_iQ_i$.
Assume that $m(x)=0$. We may also assume that $P_iQ_iR_i$ is homogeneous of degree 0
(for the rank grading) for each $i$. Then
$$\sum_i\langle \mu\,,\,P_iR_iQ_i\,,\,\lambda\rangle=0$$
for each $r$-partitions $\lambda$, $\mu$.
For $n$ large enough we consider the coefficients 
\begin{equation}
\langle
\bar\lambda_1\circledast\lambda_2\circledast\bar\lambda_3\, , P_i R_i Q_i\,,
\lambda_1\circledast\lambda_2\circledast\lambda_3\rangle,\quad \bar\lambda_1 \subset
\lambda_1,\quad \lambda_3 \subset \bar\lambda_3,\quad |\lambda_1 \backslash
\bar\lambda_1|=|\bar\lambda_3\backslash \lambda_3|=t
\end{equation}
with 
$t=\text{sup}_i (\deg(P_i)) =\text{sup}_i(-\deg(Q_i)).$
Since  $Q_i$ is an annihilation operator and $P_i$ is a
creation operator, the coefficient
\begin{equation}
\begin{gathered}
\langle \bar\lambda_1\circledast\lambda_2\circledast\bar\lambda_3\,, P_i R_i Q_i\,, 
\lambda_1\circledast\lambda_2\circledast\lambda_3\rangle\end{gathered}
\end{equation}
factorizes as in \eqref{E:PP5}, and it is zero unless $\deg(P_i)=-\deg(Q_i)= t$.
We claim that \eqref{E:tauy}, \eqref{E:rhoy} hold again for some non-zero $c,d\in K_r$.
Then \eqref{E:PP6} hold again
for all $\lambda_1,\bar\lambda_1,\lambda_2, \lambda_3, \bar\lambda_3$ as above and
all large enough $n$. If $x\neq 0$ we may assume that the elements $R_i$
are linearly independent polynomials in the $f_{0,l}$'s and that $P_i,Q_i\neq 0$.
Then the same argument as in the proof of  Proposition \ref{T:2} yields a contradiction.
The proof of the claim is the same as the proof of Lemma \ref{L:PK1}. It is left to the reader.

\end{proof}

\vspace{.1in}

\begin{rem}
In the case of the Hilbert scheme we must set
$r=1$ and $\chi_\alpha=1$ in the formulas above.
\end{rem}

\vspace{.2in}

\section{Heisenberg subalgebras}

\vspace{.15in}

This short section contains a few remarks concerning a natural
family of Heisenberg subalgebras in $\H_{\Kb}$ and their action on
$\mathbf{L}_{\Kb}=\bigoplus_n K^T(\text{Hilb}_n)_{\mathbf{K}}$.

\vspace{.2in}

\paragraph{\textbf{9.1. Heisenberg subalgebras.}} Fix $\mu \in \mathbb{Q} \cup \{\infty\}$. Write
$\mu=d/r$ with $r \geqslant 0$, $d$ and $r$ coprime (and $d=1$ if
$r=0$). Let us set for simplicity
$\mathbf{u}^{\mu}_l=\mathbf{u}_{lr,ld}$ for any $l \in \Z$. Let
$\widehat{\U}^{\mu}$ be the subalgebra of $\widehat{\U}$ generated
by $\boldsymbol{\mathcal{K}}$ and the elements $\{
\mathbf{u}^\mu_{\pm 1}, \mathbf{u}^\mu_{\pm 2}, \ldots\}$. This
algebra is isomorphic to a (quantum) Heisenberg algebra. The
defining relations are
$$[\boldsymbol{\mathcal{K}}, \mathbf{u}^\mu_{l}]=0, \qquad [\mathbf{u}^\mu_{l}, \mathbf{u}^\mu_{n}]=
\delta_{n, -l}
\frac{(\mathbf{c}_1^r\mathbf{c}_2^d)^n-(\mathbf{c}_1^r\mathbf{c}_2^d)^{-n}}{\alpha_n}.$$
As before we have set $\a_n=(1-(\sigma
\overline{\sigma})^{-n})(1-\sigma^n)(1-\overline{\sigma}^{-n})$. We
define subalgebras $\U_c^{\mu}$ in the obvious way. Two cases are of
special interest
\begin{enumerate}
\item[i)] $\mu=\infty$. The subalgebra $\U_c^{\infty}$ acts on $\mathbf{L}_{\Kb}$ by
Macdonald operators (see Section~1.4.). It also gives rise, through
the adjoint action, to the action of the Hecke operators on $\U_c$
(see Section 6). In this case the central charge vanishes and
$\U_c^{\infty}$ is a polynomial algebra,
\item[ii)] $\mu=0$. The subalgebra $\U_c^{0}$ acts on $\mathbf{L}_{\Kb}$ by means of the correspondences induced by the virtual classes $\t^*_n \otimes \Lambda(\boldsymbol{\mathcal{V}}_n)$ and $\Lambda(qt\boldsymbol{\mathcal{V}}^*_{-n})$ for $n \geqslant 1$ (see Section~5).
\end{enumerate}
Note that (when the central charge is set to zero) the various Heisenberg subalgebras
$\widehat{\U}^{\mu}$ are interchanged by the $SL(2,\Z)$-action.
When viewed as the Hall algebra of an elliptic curve,
$\widehat{\U}^{\mu}$ consists of the functions supported on the set of semistable coherent sheaves of slope $\mu$ (see \cite{BS}, Section~5).

\vspace{.2in}

\paragraph{\textbf{9.2. Casimir operators.}} There is a natural nondegenerate pairing on $\widehat{\U}^{+}$,
coming from its realization as a Hall algebra (see \cite{BS},
Section~2 and Section~4.5). It is given by
$$\langle \mathbf{u}^\mu_{l}, \mathbf{u}^\mu_{n}\rangle=\delta_{l,n}/\a_{l}.$$
More generally we have
$$\langle \mathbf{u}^{\mu}_{\lambda}, \mathbf{u}_{\nu}^{\mu}\rangle=\delta_{\lambda,\nu} \prod_i m_i(\lambda)!\, \a_{i}^{-m_i(\lambda)}$$
where $\lambda,\nu$ are partitions, $\lambda=(\lambda_1, \lambda_2,
\ldots)=(1^{m_1(\lambda)}, 2^{m_2(\lambda)}, \ldots)$ and
$\mathbf{u}^\mu_{\lambda}=\prod_i \mathbf{u}^\mu_{\lambda_i}$. We
introduce the canonical Casimir operator
$$C^\mu=\sum_{\lambda} \frac{1}{\langle \mathbf{u}^{\mu}_{\lambda},\mathbf{u}^\mu_{\lambda}\rangle}
\mathbf{u}_{\lambda}^\mu \mathbf{u}^\mu_{-\lambda}$$ where the sum
ranges over all partitions and where we have set
$\mathbf{u}^\mu_{-\lambda}=\prod_i \mathbf{u}^\mu_{-\lambda_i}$.
Although the sum defining it is infinite, the operator $C^\mu$ acts
on $\mathbf{L}_{\Kb}$ since $\mathbf{u}_{-\lambda}^{\mu}$ acts by
zero on $K^T(\text{Hilb}_n)_{\mathbf{K}}$ as soon as $|\lambda| >
n$. We expect that these Casimir operators are relevant to the study
of the monodromy of the so-called \textit{quantum differential
equation} arising in GW/DT theory of the Hilbert scheme of points in
the plane (see \cite{OP}).

\vspace{.1in}

\begin{rem}
i) There is a factorization
\begin{equation*}
\begin{split}
C^\mu&=\bigg(\sum_n \frac{1}{n!} {\a_1^n}(\mathbf{u}^\mu_1)^n (\mathbf{u}^\mu_{-1})^n \bigg)
\bigg(\sum_n \frac{1}{n!} {\a_2^n}(\mathbf{u}^\mu_2)^n (\mathbf{u}^\mu_{-2})^n \bigg)\cdots\\
\end{split}
\end{equation*}
The factors in the expression of $C^\mu$ above are not the
exponential of any natural expressions because the algebra
$\U^\mu_c$ is not commutative. The Heisenberg algebras
$\widehat{\U}^\mu$ have a Hopf algebra structure coming from their
Hall algebra realization. The primitive elements are exactly the
$\mathbf{u}^\mu_l$ with $l \in \Z$. This leads to another variant of
the Casimir operator, given by
$${C'}^\mu=\exp \bigg( {\a_1} \mathbf{u}^\mu_1 \mathbf{u}^\mu_{-1}\bigg)
\exp\bigg( {\a_2} \mathbf{u}^\mu_2 \mathbf{u}^\mu_{-2}\bigg)\cdots =
\exp \bigg( \sum_{l} {\a_l} \mathbf{u}^\mu_l
\mathbf{u}^\mu_{-l}\bigg).$$ For instance, using Section~4.7, we get
${C'}^0=\exp\big( \sum_{l \geqslant 1} ((tq)^{l/2}-(tq)^{-l/2}) p_l
\frac{\partial}{\partial p_l}\big)$. Note that $C^\mu$ and
${C'}^\mu$ coincide if the central charge vanishes, i.e., when
$\mu=\infty$.\\
ii) The group $SL(2,\mathbb{Z})$ does not act on $\U_c$ (because of the central charge), but
its unipotent subgroup $\begin{pmatrix} 1 & \mathbb{Z} \\ 0 & 1 \end{pmatrix} \simeq \mathbb{Z}$ does, via
automorphisms $\rho_n~: \mathbf{u}_{r,d} \mapsto \mathbf{u}_{r, d+nr}$. In particular,
$\rho_n( C^\mu)=C^{\mu+n}$ for any $\mu \in \mathbb{Q}$. Furthermore, the representation
$\mathbf{L}_{\Kb}$ and its twist $\mathbf{L}_{\Kb}^{\rho_1}$ are isomorphic, with intertwining operator 
$$\nabla~: \mathbf{L}_{\Kb} \mapsto \mathbf{L}_{\Kb}^{\rho_1}, \quad [I_{\lambda}] \mapsto t^{n(\lambda)}q^{n(\lambda')} [I_{\lambda}].$$
 As a consequence, the Casimir operators on $\mathbf{L}_{\Kb}$ satisfy
$$C^{\mu+1} = \nabla C^{\mu} \nabla^{-1}$$
for any $\mu \in \mathbb{Q}$.
Note that the operator $\nabla$ restricts, for any given $n$, to the action on the K-theory of the tensor product by the relative ample line bundle
$\mathcal{O}(1)$ on $\text{Hilb}_n$ (see \cite{Haiman}).
\end{rem}

\vspace{.2in}

\section{The shuffle algebra}

\vspace{.15in}

In this last section we provide an alternative algebraic description
of the algebras $\H^>_{\Kb}$ or $\U^>$. This presentation involves
a certain (noncommutative) shuffle product on the algebra of symmetric
functions $\Lb_{\Kb}$. It appeared independently in \cite{FT}.
We also discuss and compare the results of \textit{loc. cit.} with ours.

\vspace{.2in}

\paragraph{\textbf{10.1.}} One may use the faithful action of $\U^>$ on the
representation $\mathbf{L}_{\Kb}$ together with the virtual classes
$\Vcb_{l}$, $l \geqslant 1$, to give a new presentation of the $\U^>$
as some kind of $q,t$-shuffle algebra. The precise formalism,
which we now briefly recall, was introduced by B. Feigin and A. Odesskii in
\cite{FO}. Let $g(z) \in \C(z)$ be any rational function.
For $r \geqslant 1$ we put $g(z_1, \ldots, z_r)=\prod_{i<j} g(z_i/z_j)$.
Let us denote by
\begin{equation}
\begin{split}
\Sym_r: \C(z_1, \ldots, z_r) \to \C(z_1, \ldots,
z_r)^{\mathfrak{S}_r},\ P(z_1, \ldots, z_r) \mapsto \sum_{\sigma \in
\mathfrak{S}_r} P(z_{\sigma(1)}, \ldots, z_{\sigma(r)})
\end{split}
\end{equation}
the \emph{symmetrization map} and let us consider
the weighted symmetrization
\begin{equation*}
\begin{split}
\Psi_r: \C[z_1^{\pm 1}, \ldots, z_r^{\pm r}] \to \C(z_1, \ldots,
z_r)^{\mathfrak{S}_r},\quad P(z_1, \ldots, z_r) \mapsto \Sym_r \big(
g(z_1, \ldots, z_r) P(z_1, \ldots, z_r)\big).
\end{split}
\end{equation*}
We denote by $\mathbf{S}_r$ the image of $\Psi_r$.
There is a unique map $m_{r,r'}: \mathbf{S}_r \otimes \mathbf{S}_{r'}
\to \mathbf{S}_{r+r'}$ which makes the following diagram commute
\begin{equation}\label{E:shufflediagram}
\xymatrix{ \C[z_1^{\pm 1}, \ldots, z_r^{\pm 1}] \otimes \C[z_1^{\pm 1}, \ldots, z_{r'}^{\pm 1}] 
\ar[r]^-{\Psi_r \otimes \Psi_{r'}} \ar[d]_-{i_{r,r'}} & \mathbf{S}_r \otimes \mathbf{S}_{r'} \ar[d]_-{m_{r,r'}} \\ 
\C[z_1^{\pm 1}, \ldots, z_{r+r'}^{\pm 1}] \ar[r]^-{\Psi_{r+r'}} & \mathbf{S}_{r+r'}}
\end{equation}
where $i_{r,r'}\big( P(z_1, \ldots, z_r) \otimes Q(z_1, \ldots, z_{r'})\big)=
P(z_1, \ldots, z_r)Q(z_{r+1}, \ldots, z_{r+r'})$.
It is easy to check that the maps $m_{r,r'}$ endow the space
$\mathbf{S}=\C 1 \oplus \bigoplus_{r \geqslant 1} \mathbf{S}_r$
with the structure of an associative algebra.
The product in $\mathbf{S}$ may be explicitly written
as the shuffle operation
$$h(z_1, \ldots, z_r) \cdot f(z_1, \ldots, z_{r'})=
\frac{1}{r!r'!} 
\Sym_{r+r'} \bigg( \hspace{-.1in}\prod_{\substack{1 \leqslant i \leqslant r\\ r
+1 \leqslant j \leqslant r+r'}} \hspace{-.15in}g(z_i/z_j)\cdot 
\;h(z_1, \ldots, z_r) f(z_{r+1}, \ldots, z_{r+r'})\bigg).$$
Note that by construction the algebra $\mathbf{S}$ is generated by the subspace
$\mathbf{S}_1$. We may replace the ground field $\C$ by
$\mathbf{K}$ and define $\mathbf{S}$ for any $g(z)\in \mathbf{K}(z)$.
From now on we fix
\begin{equation}
g(z)=\frac{(1-tz)(1-qz)}{(1-z)(1-tqz)}.
\end{equation}
Recall that the action of $\U^0$ on $\U^>$ by Hecke operators factors
yields, for each $r$, an action
$$\bullet~: \Kb[z_1^{\pm 1}, \ldots, z_r^{\pm 1}]^{\mathfrak{S}_r}
\otimes \U^>[r] \to \U^>[r].$$

\vspace{.1in}

\begin{theo}\label{T:shuffle} The assignment
$\bt_{1,l} \mapsto z_1^l $ for $l \in \Z^*$ extends to a graded
algebra isomorphism $\Upsilon~: \U^> \to
\mathbf{S}$. Moreover, for any $P(z_1, \ldots, z_r) \in \Kb[z_1^{\pm
1}, \ldots, z_r^{\pm 1}]^{\mathfrak{S}_r}$ and any $u \in \U^>[r]$
we have
\begin{equation}\label{E:shuffle2}
\Upsilon(P(z_1, \ldots, z_r) \bullet u)=P(z_1, \ldots, z_r) \Upsilon(u).
\end{equation}
\end{theo}

As mentioned above, the above theorem may be proved by considering the action
of $\U^>$ on $\mathbf{L}_{\Kb}$ and expanding it in terms of the virtual classes. However,
since a more general theorem (giving a shuffle presentation for the Hall algebra
of any smooth projective curve) appears in \cite{SV3} we omit the proof here.

\vspace{.2in}

\paragraph{\textbf{10.2.}}
Consider the generating series
$${E}^{+}(z)=\sum_{p \in \Z} \bt_{1,p}z^p,$$
Kapranov showed in \cite{Kap}, Theorem~3.3, that the following relation holds~:
\begin{equation}\label{E:Kap1}
\zeta(z_2/z_1) {E}^{+}(z_1){E}^{+}(z_2)=\zeta(z_1/z_2)
{E}^{+}(z_2){E}^{+}(z_1),
\quad \zeta(z)=\frac{(1-\s z)(1-\bs z)}{(1-z)(1-\s\bs z)}=g(z^{-1}).
\end{equation}
Equations (\ref{E:Kap1}) is the so-called functional equation for Eisenstein
series. Let $\mathfrak{E}^>$ be the associative algebra generated by some
elements $\bt_{1,p}, \;p \in \Z$ subject to
(\ref{E:Kap1}). There is a surjective map $\pi^>~:\mathfrak{E}^> \to \U^>$.
It is known however, that this map is NOT an isomorphism.
Its kernel is presumably described by some higher rank
analogs of (\ref{E:Kap1}).
One can complete the algebra $\mathfrak{E}^>$ by adding new generators
$$E^-(z)=\sum_{p \in \Z} \bt_{-1,p}z^p,\quad
\psi^{\pm}(z)=\exp\big( \sum \bt_{0,\pm r}/\a_r\big),$$
and certain new relations,  see, e.g., \cite{FT}.
One gets in this way the so-called \textit{Ding-Iohara} algebra $\mathfrak{E}$
associated to $\zeta(z)$. There is a surjective (but not injective)
algebra map $\pi~:\mathfrak{E} \to \U_c$.

\vspace{.2in}

\paragraph{\textbf{10.3.}} In \cite{FT} the authors define the algebras
$\mathfrak E$, $\mathbf S$ and construct two representations
$$\rho~: \mathfrak{E} \to \text{End}_\Kb(\mathbf{L}_{\Kb}), \qquad
\rho': \mathbf{S}\to \text{End}_\Kb(\mathbf{L}_{\Kb}).$$
They prove that $\rho$, $\rho'$ are compatible in the sense that  we have
$\rho|_{\mathfrak{E}^>}=\rho'\circ\Upsilon \circ \pi|_{\mathfrak{E}^>}$.
The representation $\rho'$ coincides, via $\Upsilon$, with our representation
$\varphi~:\U_c \to \text{End}_\Kb(\mathbf{L}_{\Kb})$ in Proposition 4.10, i.e.,
we have $\rho'=\varphi \circ \Upsilon^{-1}$ over $\mathbf S$.

Our results imply also that
the representation $\rho$ factors through our $\varphi$,
i.e., we have $\rho=\phi \circ \pi$, that $\varphi$ is faithful,
that $\U_c$ is related to the Double Affine Hecke algebra, that
$\varphi$ can be expressed in terms of Macdonald's difference operator,
and that $\U^>$ is related to the $K$-theory of the commuting variety.
Note also that the subalgebra of $\U^>$ generated
by the virtual fundamental classes
is probably the same as the commutative subalgebra of $\mathbf S$
studied in \cite{FT} (we have not checked this).

\newpage

\centerline{APPENDICES}

\vspace{.2in}

\appendix

\section{Combinatorics}

\paragraph{\textbf{A.1. Proof of relation (\ref{E:5})}}

\vspace{.1in}

We begin with the following lemma.

\vspace{.1in}

\begin{lem}\label{L:tyt} For any $k, l \in \Z$ the class
$[\f_{1,l}, \f_{-1,k}]$ is supported on the diagonal of $\Hi\times \Hi$.
\end{lem}
\begin{proof} Let $\lambda$ be any partition and $s$ be an addable box of $\lambda$.
Let $\mu$ be the partition such that $\mu \backslash \lambda=s$. Let
$s'$ be a removable box of $\mu$ different from $s$. It is also a
removable box of $\lambda$. Define partitions $\nu$ and $\sigma$ by
$\mu \backslash \nu=s'$ and $\lambda \backslash \sigma=s'$. We have
$$\langle \nu , \f_{1,l}\f_{-1,k} \cdot \lambda \rangle=
\langle \nu, \f_{1,l} \cdot \sigma \rangle\,
\langle \sigma, \f_{-1,k} \cdot \lambda \rangle,$$
$$\langle \nu , \f_{-1,k}\f_{1,l} \cdot \lambda \rangle=
\langle \nu, \f_{-1,k} \cdot \mu \rangle\,
\langle \mu, \f_{1,l} \cdot \lambda \rangle.$$
Lemma~\ref{L:gyt} yields
\begin{equation}
\begin{split}
\langle \sigma , \f_{-1,k} \cdot \lambda \rangle \,
\langle \nu, \f_{-1,k} \cdot \mu \rangle ^{-1}
&=L_{\sigma,\nu}(q,t) L_{\lambda,\mu}(q,t)^{-1}\\
&=\frac{(t^{l_{\lambda}(u)}-q^{a_{\lambda}(u)})
(t^{l_{\lambda}(u)+1}-q^{a_{\lambda}(u)+1})}
{(t^{l_{\lambda}(u)}-q^{a_{\lambda}(u)+1})
(t^{l_{\lambda}(u)+1}-q^{a_{\lambda}(u)})}.
\end{split}
\end{equation}
where $u$ is the (unique) intersection point of $C_s \cup R_s$
and $C_t \cup R_t$. Similarly,
\begin{equation}
\begin{split}
\langle \nu , \f_{1,l} \cdot \sigma \rangle \,
\langle \mu, \f_{1,l} \cdot \lambda \rangle ^{-1}
&=L_{\lambda,\sigma}(q,t) L_{\mu,\nu}(q,t)^{-1}\\
&=\frac{(t^{l_{\lambda}(u)}-q^{a_{\lambda}(u)+1})
(t^{l_{\lambda}(u)+1}-q^{a_{\lambda}(u)})}
{(t^{l_{\lambda}(u)}-q^{a_{\lambda}(u)})
(t^{l_{\lambda}(u)+1}-q^{a_{\lambda}(u)+1})}.
\end{split}
\end{equation}
Hence $\langle \nu, \f_{1,l}\f_{-1,k} \cdot \lambda \rangle=\langle
\nu, \f_{-1,k}\f_{1,l} \cdot \lambda \rangle$ for all $\nu \neq
\lambda$. We are done.
\end{proof}

\vspace{.1in}

The precise computation of the quantities $\langle \lambda ,
[\f_{1,l},\f_{-1,k}] \cdot \lambda\rangle$ is more involved. For this it is
convenient to use the following changes of
variables, see \cite{Garsia}, Section~2. Let $\lambda$ be our fixed
partition. We label the removable boxes of $\lambda$ by $A_1, A_2,
\ldots, A_r$ from left to right, and we set $\alpha_k=y(A_k),
\beta_k=x(A_k)$. We introduce new variables $x_1, \ldots, x_r$ and
$u_0, \ldots, u_r$ by
$$x_k=t^{\alpha_k}q^{\beta_k}, \qquad u_l=t^{\alpha_{l+1}}q^{\beta_l}$$
where by convention $\beta_0=\alpha_{r+1}=-1$. We also let $C_0,
\ldots, C_{r}$ stand for the addable boxes of $\lambda$, so that
$y(C_i)=\alpha_{i+1}+1$, $x(C_i)=\beta_i +1$. Here is an example
with $\lambda=(10,9^3,6,3^2)$

\vspace{.2in}

\centerline{
\begin{picture}(200,110)
\put(0,20){\line(0,1){70}}
\put(10,90){\line(0,-1){70}}
\put(0,90){\line(1,0){30}}
\put(20,90){\line(0,-1){70}}
\put(30,90){\line(0,-1){70}} \put(0,80){\line(1,0){30}}
\put(40,70){\line(0,-1){50}} \put(0,70){\line(1,0){60}}
\put(50,70){\line(0,-1){50}} \put(60,70){\line(0,-1){50}}
\put(0,60){\line(1,0){90}} \put(70,60){\line(0,-1){40}}
\put(80,60){\line(0,-1){40}} \put(90,60){\line(0,-1){40}}
\put(0,50){\line(1,0){90}} \put(0,40){\line(1,0){90}}
\put(0,30){\line(1,0){100}}
\put(0,20){\line(1,0){100}} \put(100,30){\line(0,-1){10}}
\put(20,83){\tiny{${A_1}$}}
\put(50,63){\tiny{${A_2}$}}
\put(80,53){\tiny{${A_3}$}}
\put(90,23){\tiny{${A_4}$}}
\put(1,93){\tiny{${C_0}$}}
\put(31,73){\tiny{${C_1}$}}
\put(61,63){\tiny{${C_2}$}}
\put(91,33){\tiny{${C_3}$}}
\put(101,23){\tiny{${C_4}$}}
\multiput(0,90)(0,2){5}{\line(0,1){1}}
\multiput(0,100)(2,0){5}{\line(1,0){1}}
\multiput(10,90)(0,2){5}{\line(0,1){1}}
\multiput(30,80)(2,0){5}{\line(1,0){1}}
\multiput(40,70)(0,2){5}{\line(0,1){1}}
\multiput(60,70)(2,0){5}{\line(1,0){1}}
\multiput(70,60)(0,2){5}{\line(0,1){1}}
\multiput(90,40)(2,0){5}{\line(1,0){1}}
\multiput(100,30)(0,2){5}{\line(0,1){1}}
\multiput(110,20)(0,2){5}{\line(0,1){1}}
\multiput(100,30)(2,0){5}{\line(1,0){1}}
\multiput(100,20)(2,0){5}{\line(1,0 ){1}}
\put(160,80){{$\alpha_1=6,$}} \put(210,80){{$\beta_1=2$}}
\put(160,65){{$\alpha_1=4,$}} \put(210,65){{$\beta_1=5$}}
\put(160,50){{$\alpha_1=3,$}} \put(210,50){{$\beta_1=8$}}
\put(160,35){{$\alpha_1=0,$}} \put(210,35){{$\beta_1=9$}}
\end{picture}}
\vspace{.05in}
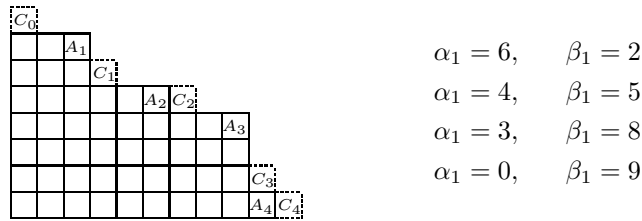
\centerline{\textbf{Figure~3.} Garsia and Tesler's change of variables.}

\vspace{.2in}

\begin{lem}[\cite{Garsia}]\label{L:GT} The following formulas hold
\begin{enumerate}
\item[(a)] for any $m \in \Z$ we have
\begin{equation}\label{E:GT1}
\sum_i x_i^m-\sum_j u_j^m=(1-t^{-m})(1-q^{-m})B^{m}_{\lambda}(q,t)
-t^{-m}q^{-m},\quad
B^{m}_{\lambda}(q,t)=\sum_{s\in\lambda}q^{mx(s)}t^{my(s)},
\end{equation}
\item[(b)] if $\mu\subset \lambda$ is such that
$\lambda \backslash \mu= A_i$ then we have
\begin{equation}\label{E:GT2}
L_{\lambda,\mu}(q,t)=\frac{tqx_i^{-1}}{(1-t)(1-q)} \prod_{j=0}^r
(u_j-x_i)\prod_{\substack{j=1\\ j \neq i}}^r(x_j-x_i)^{-1},
\end{equation}
\item[(c)] if $\mu\supset \lambda$ is such that
$\mu \backslash \lambda= C_i$ then we have
\begin{equation}\label{E:GT3}
L_{\lambda,\mu}(q,t)=\frac{u_i^{-1}}{qt} \prod_{j=1}^r (u_i-x_j)\prod_{\substack{j=0\\ j \neq i}}^r(u_i-u_j)^{-1}.
\end{equation}
\end{enumerate}
\end{lem}

\vspace{.1in}

Using the above (\ref{E:GT2}) and (\ref{E:GT3}) we obtain,
after a little arithmetic
\begin{equation}\label{E:kil}
\langle \mu, \f_{1,l} \cdot \lambda \rangle
=\frac{1}{(1-q)(1-t)}(qt)^{l}u_i^{l} \cdot \prod_{j=1}^r (u_i-x_j)
\prod_{\substack{j=0 \\j \neq i}}^r(u_i-u_j)^{-1}
\ \text{if}\ \mu \supset \lambda\ \text{and}\  \mu \backslash \lambda =C_i,
\end{equation}
\begin{equation}\label{E:kkil}
\langle \nu, \f_{-1,k} \cdot \lambda \rangle =\frac{tq}{(1-q)(1-t)}
x_i^{k-1} \cdot \prod_{j=0}^r (u_j-x_i) \prod_{\substack{j=1 \\ j
\neq i}}^r (x_j-x_i)^{-1}
\ \text{if}\ \nu\subset\lambda\ \text{and}\ \lambda \backslash \nu=A_i,
\end{equation}
\begin{equation}\label{E:kilo}
\begin{split}
\langle \lambda, \f_{-1,k}\f_{1,l} \cdot \lambda \rangle &=\frac{1}{(1-q)(1-t)}\sum_{i=0}^r \;(qt)^{l+k+1}u_i^{l+k+1} \cdot \prod_{j=1}^r \frac{u_i-x_j}{tqu_i-x_j} \cdot \prod_{\substack{j=0 \\ j \neq i}}^r \frac{tq u_i-u_j}{u_i-u_j},\\
\langle \lambda, \f_{1,l}\f_{-1,k} \cdot \lambda \rangle
&=\frac{tq}{(1-q)(1-t)}\sum_{i=1}^r \;x_i^{l+k+1} \cdot
\prod_{j=0}^r \frac{u_j-x_i}{tqu_j-x_i} \cdot \prod_{\substack{j=1
\\ j \neq i}}^r \frac{tq x_j-x_i}{x_j-x_i}.
\end{split}
\end{equation}
In particular, the commutator $[\f_{-1,k},\f_{1,l}]$, 
being supported on the diagonal, only depends on $k+l$. 
For $k+l=m$ we put
$$c_{m}=[\f_{-1,k},\f_{1,l-1}].$$
We now establish a formula for the generating series of $(c_m)$.

\vspace{.1in}

\begin{prop}\label{L:GTRE} We have

(a) 
$(1-q)(1-t)\sum_{m\geqslant 0} c_{m} s^m=$
$$=\frac{1}{(1-qt)} \Big[ 1 -qt\; \exp\Big(\sum_{n \geqslant 1}
(1-t^{n}q^{n}) \big( (1-t^{-n})(1-q^{-n})\f_{0,n}
-t^{-n}q^{-n}\big) \frac{s^n}{n}\Big)\Big],
$$

(b)
$(1-q)(1-t)\sum_{m\geqslant 0} c_{-m} s^m=$
$$=\frac{1}{(1-qt)} \Big[-qt + \exp\Big(\sum_{n \geqslant 1}
(1-t^{-n}q^{-n}) \big( (1-t^n)(1-q^n)\f_{0,-n} -t^{n}q^{n}\big)
\frac{s^n}{n}\Big)\Big].$$
\end{prop}
\begin{proof} We sketch the proof of the first equality, the second is similar.
Using the change of variables introduced above,
(\ref{E:GT1}) and (\ref{E:kilo}),
we are reduced to showing the following identity
\begin{equation}\label{E:kop}
 \sum_{i=0}^r \frac{1}{1-tqu_is} T_i -\sum_{i=1}^r \frac{tq}{1-x_is}S_i  =\frac{1}{1-tq} -\frac{tq}{1-tq}\prod_{i=1}^r \frac{1-tq x_i s}{1-x_is} \prod_{i=0}^r \frac{1-u_is}{1-tqu_is}
\end{equation}
where we have set
$$T_i=\prod_{j=1}^r \frac{u_i-x_j}{tqu_i-x_j} \cdot \prod_{\substack{j=0 \\ j \neq i}}^r \frac{tq u_i-u_j}{u_i-u_j},\qquad
S_i=\prod_{j=0}^r \frac{u_j-x_i}{tqu_j-x_i} \cdot \prod_{\substack{j=1 \\ j \neq i}}^r \frac{tq x_j-x_i}{x_j-x_i}.$$
The proof of (\ref{E:kop}) is standard. Namely, by direct computations we check that the left hand side and the right hand side 
have the same poles and residues. Since both are rational functions of degree zero, they differ by at most a constant.
To compute this constant we set $s=0$. This reduces (\ref{E:kop}) to
\begin{equation}
\sum_{i=0}^r \prod_{j=1}^r \frac{u_i-x_j}{tqu_i-x_j} \cdot
\prod_{\substack{j=0 \\ j \neq i}}^r \frac{tq
u_i-u_j}{u_i-u_j}-tq\sum_{i=0}^r\prod_{j=0}^r
\frac{u_j-x_i}{tqu_j-x_i} \cdot \prod_{\substack{j=1 \\ j \neq i}}^r
\frac{tq x_j-x_i}{x_j-x_i}=1
\end{equation}
which is itself a corollary of the following formula : for any variables
$\xi_1, \ldots, \xi_n$,  $p$ we have
$$\sum_{i=1}^n \left(\prod_{j \neq i} \frac{p \xi_i-\xi_j}{\xi_i-\xi_j}\right)=1+p + \cdots + p^{n-1}.$$
We leave the details to the reader. \end{proof}

\vspace{.1in}

It is straightforward to check (\ref{E:5}) from the above proposition and (\ref{E:hrd1}--\ref{E:hrd2}).

\vspace{.4in}

\paragraph{\textbf{A.2. Proof of Claim \ref{claim:5.4}}}
Let us denote by $A_{r,s}=A_{r,s}(x_i,x'_i,u_i,u'_i)$ the right hand side
of \eqref{E:vert8}. This is a rational function of degree one, whose poles 
are of order at most one and are located along the hyperplanes 
$x_i=x_j$, $u'_i=u'_j$ and $u'_i=x_j$ for $i \neq j$. A simple calculation 
shows that the residues of $A_{r,s}$ vanish along all these hyperplanes, hence 
$A_{r,s}$ is in fact a polynomial in the variables $x_i,$ $x'_i,$ $u_i,$ $u'_i$
which we may write as
$$A_{r,s}=\sum_i \alpha_i x_i + \sum_i \alpha'_i x'_i + 
\sum_i \beta_i u_i + \sum_i \beta'_i u'_i + \gamma$$
for some scalars $\alpha_i,$ $\alpha'_i,$ $\beta_i,$ $\beta'_i$ and 
$\gamma$. By considering the leading coefficients of $A_{r,s}$ in $x_i$ and 
$u'_i$ one sees that $\alpha_i=\beta'_i=-1$ for all $i$. Next, observe that the
specialization of $A_{r,s}$ at $u_r=x_r$ is equal to $A_{r-1,s}$ and that 
similarly the specialization of $A_{r,s}$ at $x'_s=u'_s$ is equal to 
$A_{r,s-1}$. It follows that $\beta_i=-\alpha_i=1$ and $\alpha'_i=-\beta'_i=1$,
and that $\gamma=A_{0,0}=0$. The claim is proved.

\vspace{.2in}

\section{Hecke correspondences in higher rank}

\vspace{.1in}

Fix a pair of vector space $E_1$, $E_2$ of dimension $n,n+1$
respectively. We define a $T$-equivariant complex of vector bundles
over $M_{r,n}\times M_{r,n+1}$ by
\begin{equation}
\mathrm{Hom}(E_1,E_2)\ {\buildrel\sigma\over\to}\
\begin{gathered}
q\mathrm{Hom}(E_1,E_2)\oplus t\mathrm{Hom}(E_1,E_2)\cr \oplus\cr
qt\mathrm{Hom}(E_1,\C^r)\oplus\mathrm{Hom}(\C^r,E_2)
\end{gathered}\ {\buildrel\tau\over\to}\
qt\mathrm{Hom}(E_1,E_2)\oplus qt,
\end{equation}
where the maps $\sigma$, $\tau$ are defined by
$$\sigma(\xi)=\begin{pmatrix}\xi a_1-a_2\xi\cr\xi
b_1-b_2\xi\cr-\varphi_2\xi\cr\xi v_1\end{pmatrix},\quad
\tau\begin{pmatrix}c\cr d\cr\Phi\cr
V\end{pmatrix}=([a,d]+[c,b]+V\varphi_1+v_2\Phi)\oplus(\mathrm{tr}_{E_1}(v_1
\Phi)\oplus\mathrm{tr}_{E_2}(V\varphi_2)),$$ with
$$[a,d]=a_2d-da_1,\quad[c,b]=cb_1-b_2c.$$
It is well-known that the unique non-zero cohomology group is
$\mathcal V=\mathcal{K}er(\tau)/\mathcal{I}m(\sigma)$ and that the latter is a
$T$-equivariant vector bundle of rank $r(2n+1)-1$ which admits a
section which vanishes precisely on the Hecke correspondence
$M_{r,n,n+1}$. Therefore, for each $\mu\subset\lambda$ we have
$$N_{\mu,\lambda}=[\mathcal V^*|_{I_{\mu,\lambda}}]. $$

\vspace{.2in}

\section{Complements on base change}

\vspace{.1in}

Consider the following Cartesian square of smooth connected varieties
$$\xymatrix{ X'\ar[d]_g&Y'\ar[l]_-{f'}  \ar[d]_-{g'} \\ X&\ar[l]_{f}  Y.}$$
Observe that the maps $f,f',g,g'$ are regular. So the derived pull-back in
K-theory are well-defined.
Let $x,y,x',y'$ be the dimensions of $X,Y,X',Y'$.
The aim of this appendix is to prove the following.

\begin{prop}\label{D:1} Assume that $g$ is proper, that $f'\times g'$
is a closed embedding $Y'\to X'\times Y$ and that
$x+y'=x'+y$. Then $g'$ is proper and we have
$Lf^*\circ Rg_*=Rg'_*\circ L(f')^*$,
an equality of group homomorphisms
$K(X')\to K(Y)$.
\end{prop}

\begin{proof}
First, observe that $g'$ is proper because the base-change of a proper map is
again proper. Now we prove the proposition in two steps.

$(a)$ Assume that $g$, $g'$ are closed (regular) embeddings.
The hypothesis on the dimensions implies that
$\mathrm{Tor}_k^{\mathcal O_X}(\mathcal O_{X'},\mathcal O_Y)=0$ for each $k>0$,
see \cite{T}, Lemma 3.2.
Therefore $X'$ and $Y$ are tor-independent over $X$, see
\cite{SGA6}, Definition~1.5.

$(b)$ Now, assume that $g$, $g'$ are any regular morphisms.
Consider the following diagram with Cartesian squares
$$\xymatrix{ X'\ar[d]_g&X'\times Y\ar[l]_-{p_1}\ar[d]_{g\times 1}
&Y'\ar[l]_-{f'\times g'}  \ar[d]_-{g'} \\
X&X\times Y\ar[l]_-{p_1}&\ar[l]_-{f\times 1}  Y.}$$
Here $p_1$ means the projection along the second factor.
Since $p_1$ is a smooth morphism the varieties $X'$ and $X\times Y$ are
tor-independent over $X$.
Since $f\times 1$ and $f'\times g'$ are closed
(regular) embeddings the varieties $X'\times Y$ and $Y$ are
tor-independent over $X\times Y$
by part $(a)$.
Thus the varieties $X'$ and $Y$ are tor-independent over $X$
by \cite{SGA6}, Lemma III.1.5.1.
Thus \cite{SGA6}, Proposition IV.3.1.1, yields
$$Lf^*\circ Rg_*=Rg'_*\circ L(f')^*.$$

\end{proof}

\vskip1cm

\vspace{.2in}

\centerline{\textsc{Acknowledgements}}

\vspace{.1in}

We would like to thank A. Braverman, H. Nakajima, W. Wang and especially A. Oblomkov for useful discussions and correspondence. We are also grateful to M. Finkelberg and A. 
Tsymbaliuk for communicating on \cite{FT}. This work was started while O.S. was visiting the Institute for Advanced Studies at Princeton, which he would like to thank for the 
stimulating research environment.

\newpage

\newpage

\centerline{\textsc{Index of notations}}

\vspace{.4in}

\begin{tabular}{l @{\hspace{.3in}} r |  l @{\hspace{.3in}} r}
{\textbf{Rings}} & &{\textbf{Varieties, Bundles}} &\\
&&&\\
$R=R_T=\C[q^{\pm 1}, t^{\pm 1}]$ & 2.2 & $\text{Hilb}_n$ & 2.1\\
$K=K_T=\C(q,t)$ & 3.2 & $\boldsymbol{\tau}_n, \boldsymbol{\tau}_{\lambda}$ & 2.3\\
$\Kb=\C(\s^{1/2},\bs^{1/2})=\C(q^{1/2},t^{1/2})$ & 1.1, 1.3 & $T_{\mu,\lambda}, T_{\lambda}, N_{\mu,\lambda}$ & 2.4\\
$\s=q^{-1}, \bs=t^{-1}$ & 1.1, 1.3 & $Z_{n,k}$ & 2.4\\
& & $\boldsymbol{\tau}_{n,k}, \boldsymbol{\tau}_{\lambda,\mu}$ & 2.5\\
{\textbf{Algebras}} & & $C_E$ & 7.3\\
& & &\\
$\widehat{\U}, \widehat{\U}^{\pm}, \widehat{\U}^0, \widehat{\U}^>, \widehat{\U}^<, \widehat{\U}^{\geqslant}, \widehat{\U}^{\leqslant}$ & 1.1 & {\textbf{Others}} &\\
$\U,\U^{\pm}, \U^0, \U^>\U^<, \U^{\geqslant}, \U^{\leqslant}$ & 1.1& & \\
$\ddot{\H}_n, \mathbf{S}\ddot{\H}_n, \SH^+_{n}, \SH^+_{\infty}$ & 1.3 & $\Zb, \Zb^*, \Zb^>, \Zb^<, \Zb^{\geqslant}, \Zb^{\leqslant}, \Zb^{\pm}, \Zb^0$ & 1.1\\
$\SH^>_{\infty}, \SH^{\geqslant}_{\infty}$ & 1.4 & $\Zb^{-1}, \Zb^{0}, \Zb^{1}$ & 1.1\\
$\H_K, \H_{\Kb}$ & 3.4 & $\Psi_l$ (Adams operations) & 3.4\\
$\H_{\Kb}^>, \H_{\Kb}^<, \H_{\Kb}^0$ & 4.6 & $\bt_{r,d}$ & 1.1\\
$\H_{\Kb}^{\geqslant}, \H_{\Kb}^{\leqslant}$ & 4.4 & $\mathbf{h}_{r,d}$ & 3.4\\
$\mathbf{E}_K, \mathbf{E}_{\Kb}$ & 3.2 & $\Vcb_r, \Lambda(\Vcb_r)$ (virtual class) & 5.2\\
$\bar\Cb_R, \bar\Cb_{\Kb}$ & 7.4 & $\gamma_t \in \text{End}(\Lb_{\Kb})$ (plethysm) & 3.3\\
$\overline\SCb_R, \overline\SCb_{\Kb}$ & 7.8 & $\omega~: \U^{\geqslant} \to \U^>$ & 6.1\\
$\SCb_R, \SCb_{\Kb}$ & 7.9 & $\pi_r~: \U^0 \to \Kb[z_1^{\pm 1}, \ldots, z_r^{\pm 1}]^{\mathfrak{S}_r}$ & 6.1\\
$\mathbf{S}=\mathbf{S}_{g(z)}$ & 9.1 & $\pi_E~: \H^0_{\Kb} \to R_{T \times GL(E)} \otimes_R \Kb$ & 7.7\\
& & $\bullet~: \U^0 \otimes \U^> \to \U^>$ & 6.1\\
{\textbf{Spaces}} & & $\bullet~: \H_{\Kb}^0 \otimes \text{End}(\mathbf{L}_{\Kb}) \to \text{End}(\mathbf{L}_{\Kb})$ & 6.4\\
& & $C^\mu, {C'}^\mu$ (Casimir operators) & 9.2\\
$\boldsymbol{\Lambda}_{\Kb}$, $\boldsymbol{\Lambda}^n_{\Kb}$ & 0.4 & $\nabla$ (nabla operator) & 9.2\\
$\mathbf{L}_R,\mathbf{L}_K, \mathbf{L}_{\Kb}$ & 3.2 & &\\
$\mathbf{M}_{\Kb} \simeq \mathbf{L}_{\Kb}$ & 7.5 & &\\
& & &\\
{\textbf{Maps, Representations}} & & &\\
& & &\\
$\widetilde{\varphi}_\infty~: \SH^{\geqslant}_{\infty} \to \text{End}(\Lb_{\Kb})$ & 1.4 & &\\
$\widetilde{\varphi}~: \U^{\geqslant} \to \text{End}(\Lb_{\Kb})$ & 1.4 & &\\
$\varphi~: \U_c \to \text{End}(\Lb_{\Kb})$ & 4.7 & &\\
$\psi~: \H_{\Kb} \to  \text{End}(\Lb_{\Kb})$ & 3.2 & &\\
$\rho~: \bar\Cb_{\Kb} \to  \text{End}(\mathbf{M}_{\Kb})=\text{End}(\Lb_{\Kb})$ & 7.6 & &\\
$\Omega~: \U_c \to \H_{\Kb}$ & 3.4 & &\\
$ \Phi^+_{\infty}~: \U^+ \to \SH^+_{\infty}$ & 1.3 & &\\
$\Xi~: \SCb_{\Kb} \to \H_{\Kb}^>$ & 7.9 & &\\
$\Gamma~:\SCb_{\Kb} \to \U^>$ & 7.9 & &\\
$\Upsilon~: \U^> \to \mathbf{S}$ & 9.1 & &\\
\end{tabular}

\newpage

\small{}

\vspace{4mm}

\noindent
O. Schiffmann, \texttt{olive@math.jussieu.fr},\\
D\'epartement de Math\'ematiques, Universit\'e de Paris 6, 175 rue du Chevaleret, 75013 Paris, FRANCE.

\vspace{.1in}

\noindent
E. Vasserot, \texttt{vasserot@math.jussieu.fr},\\
D\'epartement de Math\'ematiques, Universit\'e de Paris 7, 175 rue du Chevaleret, 75013 Paris, FRANCE.


\begin{thebibliography}{99}

\bibitem{BS}
I. Burban, O. Schiffmann, \emph{On the Hall algebra of an elliptic curve, I.}, preprint math.AG/0505148, (2005).

\bibitem{CO}
E. Carlsson, A. Okounkov, \emph{Exts and vertex operators}, preprint arXiv:0801.2565 (2008).

\bibitem{Cheah}
J. Cheah, \emph{Cellular decompositions for nested Hilbert schemes of points},
 Pacific Journal of Mathematics, Vol. \textbf{183}, No. 1, (1998), 39--90.

\bibitem{C}
I. Cherednik, \emph{Double affine Hecke algebras},  Cambridge University Press, (2004).

\bibitem{CG}
N. Chriss, V. Ginzburg, \emph{Representation theory and complex geometry}, Birkha\" user (1996).

\bibitem{ES}
G. Ellingsrud, S. A. Str\o mme , \emph{On the homology of the Hilbert scheme
of points in the plane}, Invent. Math. \textbf{87} (1987), p. 343--352.

\bibitem{FO}
B. Feigin, A. Odesskii, \emph{Vector bundles on an elliptic curve and Sklyanin algebras}, Topics in quantum groups and finite-type invariants, 65--84, Amer. Math. Soc. Transl. Ser. 2, \textbf{185}, (1998).

\bibitem{FT}
B. Feigin, A. Tsymbaliuk, \emph{Heisenberg action in the equivariant $K$-theory of Hilbert schemes via Shuffle Algebra}, preprint arXiv:0904.1679 (2009).


\bibitem{FGV}
E. Frenkel, D. Gaitsgory, K. Vilonen, \emph{On the geometric Langlands conjecture},
J. Amer. Math. Soc. \textbf{15} (2002), no. 2, 367--417.

\bibitem{Garsia}
A. Garsia, G. Tesler, \emph{Plethystic formulas for Macdonald $q,t$-Kostka coefficients}, Adv. Math. \textbf{123}, no. 2, (1996), 144--222.

\bibitem{Go}
L. G\"ottsche,  \emph{Hilbert schemes of zero-dimensional subschemes of smooth varieties}, Lecture Notes in Mathematics, \textbf{1572}, (1994).

\bibitem{Groj}
I. Grojnowski, \emph{ Instantons and affine algebras I: the Hilbert scheme and vertex operators}, Math. Res. Letters \textbf{3} (1996), 275--291.

\bibitem{Haiman}
M. Haiman, \emph{Notes on Macdonald polynomials and the geometry of Hilbert Schemes}, Symmetric Functions 2001, Kluwer (2002), 1--64.

\bibitem{Kap}
M. Kapranov, \emph{Eisenstein series and quantum affine algebras},
Algebraic geometry, \textbf{7}. J. Math. Sci. (New York) \textbf{84} (1997), no. 5, 1311--1360.

\bibitem{Lehn}
M. Lehn, \emph{Chern classes of tautological sheaves on Hilbert schemes of points on surfaces}, Invent. Math. \textbf{136} (1999), no. 1, 157--207.

\bibitem{SL}
M. Lehn, C. Sorger, \emph{Symmetric groups and the cup product on the cohomology of Hilbert schemes}, Duke Math. J. \textbf{110} (2001), no. 2, 345--357.

\bibitem{L}
G. Lusztig,
\emph{Quivers, perverse sheaves, and quantized enveloping algebras},
J. Amer. Math. Soc., \textbf{4} (1991), 365--421.

\bibitem{Mac}
I.G. Macdonald , \emph{Symmetric functions and Hall polynomials, second edition}, Oxford Math. Mon., (1995).

\bibitem{MO}
D. Maulik, A. Okounkov, \emph{Nested Hilbert schemes and symmetric functions}, unpublished.

\bibitem{Nak1}
H. Nakajima, \emph{Heisenberg algebra and Hilbert schemes of points on projective surfaces},  Ann. of Math. (2) \textbf{145} (1997), no. 2, 379--388.

\bibitem{NY}
H. Nakajima, K. Yoshioka, \emph{Instanton counting on blowup. I.
4-dimensional pure gauge theory}, Invent. Math. \textbf{162} (2005),
313--355.

\bibitem{QinWang1}
W.-P. Li, Z. Qin, W. Wang, \emph{Vertex algebras and the cohomology ring structure of Hilbert schemes of points on surfaces}, Math. Ann. \textbf{324} (2002), no. 1, 105--133.

\bibitem{OP}
A. Okounkov, R. Pandharipande, \emph{The quantum differential equation of the Hilbert scheme of points in the plane}, preprint arXiv:0906.3587 (2009).

\bibitem{SGA6}
P. Berthelot, A. Grothendieck, A. Illusie,
\emph{Th\'eorie des intersections et Th\'eor\`eme de Riemmann-Roch},
Lect. Notes Math. \textbf{162} (1971).

\bibitem{Scano}
O. Schiffmann, \emph{On the Hall algebra of an elliptic curve, II.}, preprint arXiv:math/0508553 (2005).

\bibitem{SLectures2}
O. Schiffmann, \emph{Lectures on canonical and crystal bases of Hall algebras.}, in preparation.

\bibitem{SV}
O. Schiffmann, E. Vasserot, \emph{The elliptic Hall algebra, Cherednick Hecke algebras and Macdonald polynomials}, preprint arXiv:0802.4001 (2008).

\bibitem{SV5}
O. Schiffmann, E. Vasserot, \emph{Degenerate double affine Hecke algebras and the equivariant
cohomology of the moduli space of torsion free sheaves over $\Bbb A^2$}, in preparation.

\bibitem{SV3}
O. Schiffmann, E. Vasserot, 

\bibitem{T}
R.W. Thomason, \emph{Les K-groupes d'un sch\'ema \'eclat\'e et une
formule d'intersection exc\'edentaire},
Invent. Math. \textbf{112} (1993),
195--215.

\bibitem{V}
E. Vasserot, \emph{Sur l'anneau de cohomologie du sch\'ema de Hilbert de 
$\bold C\sp 2$},. C. R. Acad. Sci. Paris S\'er. I Math. \textbf{332} (2001), 
no. 1, 7--12.

\bibitem{VV}
M. Varagnolo, E. Vasserot, \emph{On the K-theory of the cyclic quiver variety.}  Internat. Math. Res. Notices (1999), no. \textbf{18}, 1005--1028.

\end{thebibliography}
\end{document}